\documentclass[11pt,reqno]{amsart}

\RequirePackage{amsthm}   

\usepackage[T1]{fontenc}
\usepackage[utf8]{inputenc} 
\usepackage{ifthen}
\usepackage{xstring}
\usepackage{xparse}
\usepackage{xargs}
\usepackage{xspace} 
\usepackage{color}

\usepackage{libertine}
\usepackage[british]{quotmark}
\usepackage{mathtools,textcomp} 
\usepackage{mathalfa} 

\usepackage[libertine]{newtxmath} 

\makeatletter
\DeclareDocumentCommand{\mathscr}{m}{\mbox{{\usefont{U}{eus}{m}{n} #1}}}

\makeatother

\makeatletter
\DeclareFontFamily{T1}{pzc}{}
\DeclareFontShape{T1}{pzc}{m}{it}{<->s *[1.2] pzcmi7t}{} 
\DeclareMathAlphabet{\mathpzc}{T1}{pzc}{m}{it}

\DeclareFontFamily{T1}{Pzc}{}
\DeclareFontShape{T1}{Pzc}{m}{it}{<->s*[1.4] pzcmi7t}{} 
\DeclareMathAlphabet{\mathPzc}{T1}{Pzc}{m}{it}

\makeatother

\usepackage{pifont} 
\usepackage{bullcntr} 
\usepackage{stmaryrd} 
\usepackage{upgreek} 
\usepackage{ushort} 
\usepackage{accents} 

\usepackage{pgf,pgffor} %
\usepgfmodule[matrix]
\usepackage{tikz}
\usetikzlibrary{%
	calc, %
	math,%
	arrows,
	arrows.meta,%
	trees,%
	shapes.misc,
	shapes.arrows,%
	shapes.geometric,%
	shapes.symbols,
	chains,%
	matrix,%
	positioning,
	scopes,%
	decorations.markings,%
	decorations.pathreplacing,%
	decorations.pathmorphing,
	shadows,%
	backgrounds,%
	graphs,%
	quotes,%
	fit%
}

\newcommand*{\midsloppy}{%
	\tolerance 5000%
	\hbadness 4000
	\emergencystretch 1.5em%
	\hfuzz .1pt 
	\vfuzz\hfuzz}

\makeatletter
\newcommand{\startcent}[1][-0]{\begin{center} \vspace{#1\baselineskip} }
	\newcommand{\stopcent}[1][-0]{\end{center}\vspace{#1\baselineskip}}
\makeatother

\makeatletter
\NewDocumentCommand{\startalignl}{O{}}{
	\pgfkeys{/vskip/.initial=-0.25,/wscale/.initial=0.9,#1}
	\begin{center}
		\vspace{\pgfkeysvalueof{/vskip}\baselineskip}
		\begin{minipage}[t]{\pgfkeysvalueof{/wscale}\textwidth}}
		\makeatother
		
		\makeatletter
		\NewDocumentCommand{\stopalignl}{O{}}{ 
			\pgfkeys{/vskip/.initial=-0.25,#1}
	\end{minipage}\end{center}
	\vspace{\pgfkeysvalueof{/vskip}\baselineskip}} 
\makeatother    

\usepackage{enumitem}  
\newlist{Enum}{enumerate}{3}
\setlist[Enum]{itemsep=0.2em,labelsep=0.2em,
	leftmargin=*,listparindent=0.0em,font=\sffamily\mdseries\upshape,align=left}
\setlist[Enum,1]{label=\arabic*.,ref=\arabic*}
\setlist[Enum,2]{label= \theEnumi.\arabic*.,ref=\arabic*}
\setlist[Enum,3]{label= \theEnumi.\theEnumii.\arabic*.}

\newlist{ItemList}{itemize}{3} 
\setlist[ItemList]{itemsep=0.2em,labelsep=0.2em,
	leftmargin=*,listparindent=0.0em,font=\sffamily\mdseries\upshape,align=left}
\setlist[ItemList,1]{label={\scriptsize{\ding{109}}}}
\setlist[ItemList,2]{label={\scriptsize{\ding{71}}}}
\setlist[ItemList,3]{label={\scriptsize{\ding{73}}}}

\newcommand{\startitem}{\begin{ItemList}}  
	\newcommand{\stopitem}{\end{ItemList}}

\NewDocumentCommand{\itemnob}{o}{\IfNoValueTF{#1}{
		\item}{\item[\mth{#1\text{\textdegree}\rparen}]}} 
 
\NewDocumentCommand{\nobp}{D(){}}{\mth{#1\text{\textdegree}\rparen}} 

\def\rmnum(#1){\MakeUppercase{\romannumeral #1}}
\NewDocumentCommand{\rmNb}{D(){}}{\mth{(\text{{\sc{#1}}})}}

\NewDocumentCommand{\sectitem}{ o g}{\medskip\noindent{\IfNoValueTF{#1}{}{#1.~~}
		\IfNoValueTF{#2}{}{{\slshape #2}.\,}}}

\NewDocumentCommand{\fsize}{oO{1}m}{{\IfNoValueTF{#1}{#3}{%
			\ifmmode{\text{\fontsize{#1}{#2}\selectfont\ensuremath{#3}}}%
			\else{\fontsize{#1}{#2\baselineskip}\selectfont #3}\fi}}}

\newcommand*{\bsym}[1]{\boldsymbol{#1}}
\newcommand*{\mc}[1]{\mathcal{#1}}
\newcommand*{\mscr}[1]{\mathscr{#1}} 
\newcommand*{\mpzc}[1]{\mathpzc{#1}}
\newcommand*{\mb}[1]{\mathbb{#1}}
\newcommand*{\mrm}[1]{\mathrm{#1}}

\newcommand*{\msf}[1]{\ensuremath{\mathsf{#1}}}
\DeclareDocumentCommand{\pmsf}{D(){}}{{\ensuremath{\mathsf{(#1)}}}}
\newcommand*{\mfr}[1]{\mathfrak{#1}} 

\DeclareDocumentCommand{\mth}{m}{\ensuremath{#1}} 
\DeclareDocumentCommand{\smth}{m}{\ensuremath{{\scriptstyle{#1}}}}
\DeclareDocumentCommand{\ssmth}{m}{\ensuremath{{\scriptscriptstyle{#1}}}}

\DeclareDocumentCommand{\mthbox}{m}{\mbox{\ensuremath{#1}}}

\NewDocumentCommand{\msp}{O {3.6}}{\ensuremath{\mspace{#1mu}}}
\newcommand*{\nsc}{\ensuremath{\mspace{-0.5mu},\mspace{-0.64mu}}}

\def\kwd#1{\mth{\bsym{\mrm{#1}}}}

\makeatletter

\DeclareDocumentCommand{\bar}{oD(){}}{\IfNoValueTF{#1}{\overline{#2}}{
		\overline{#2\,}\sp{\hskip-0.4pt \scriptstyle{#1}}}}
\makeatother 

\NewDocumentCommand{\ubar}{D(){}}{{\ushortw{#1}}} 
\newcommand*{\wt}[1]{\widetilde{#1}} 
\newcommand*{\sprime}{\ssmth{\prime}} 

\DeclareDocumentCommand{\Cramped}{D<>{} O{} m}{
	{ 
		\IfEqCase*{#2}{
			{}{\cramped[#1]{#3}}%
			{clap}{\crampedclap[#1]{#3}}%
			{llap}{\crampedllap[#1]{#3}}%
			{rlap}{\crampedrlap[#1]{#3}}%
		}[\textrm{Wrong [opt]!! must be [],[clap], [llap] or [rlap]}] 
} }

\makeatletter
\let\OrigSb\sb
\DeclareDocumentCommand{\sb}{O{} m}{\OrigSb{\Cramped<\scriptstyle>[#1]{\mspace{-0.72mu}#2}}} 
\DeclareDocumentCommand{\ssb}{O{} m}{\OrigSb{\Cramped<\scriptscriptstyle>[#1]{\mspace{-0.54mu}#2}}}
\newcommandx*{\sbsb}[1]{{\sb{\sb{#1}}}}
\DeclareDocumentCommand{\sbp}{u{(}u{)}}{\OrigSb{
		\cramped[\scriptstyle]{\mspace{-0.72mu}#1(#2)}}} 
\DeclareDocumentCommand{\sbseq}{O{}D(){}}{\OrigSb{\mspace{-0.72mu}
		\Cramped<\scriptstyle>[#1]{\cseq(#2)}}}     
\makeatother

\makeatletter
\let\OrigSp\sp
\DeclareDocumentCommand{\sp}{O{} m}{\OrigSp{
		\Cramped<\scriptstyle>[#1]{\mspace{-0.5mu}#2}}} 
\DeclareDocumentCommand{\ssp}{O{} m}{\OrigSp{\Cramped<\scriptscriptstyle>[#1]{\mspace{-0.5mu}#2}}}
\DeclareDocumentCommand{\spp}{u{(}u{)}}{\OrigSp{
		\cramped[\scriptstyle]{\mspace{-0.72mu}#1(#2)}}}
\makeatother
\makeatletter
\DeclareDocumentCommand{\sbsp}{O{}}{
	\foreach \x [count=\cpt]  in {#1} {%
		\ifnum\cpt=1%
		\sb{\scriptstyle{\x}}%
		\else\ifnum\cpt=2%
		\sp{\scriptstyle{\x}}%
		\else%
		\x\fi\fi}}
\makeatother

\NewDocumentCommand{\sbspscript}{o m g}{\IfNoValueTF{#1}{
		\IfNoValueTF{#3}{\sb{#2}}{\sb{#2}\sp{#3}}}{
		\IfNoValueTF{#3}{\sb{#1\nsc #2}}{\sb{#1\nsc#2}\sp{#3}}}}

\makeatletter
\NewDocumentCommand{\MultiSubSupscript}{m G{} G{} g g}{
	\IfNoValueTF{#5}{ \IfNoValueTF{#4}{
			{}\OrigSb{\cramped[\scriptstyle]{#2}\mspace{-0.8mu}}#1\OrigSb{\mspace{-0.72mu}
				\cramped[\scriptstyle]{#3}}}{{}\OrigSb{\cramped[\scriptstyle]{#2}\mspace{-0.8mu}}#1\OrigSb{\mspace{-0.72mu}
				\cramped[\scriptstyle]{#3}}\OrigSp{\mspace{-0.72mu}
				\cramped[\scriptstyle]{#4}}}}{
		{}\OrigSb{\cramped[\scriptstyle]{#2}\mspace{-0.8mu}}\OrigSp{\mspace{-0.8mu}
			\cramped[\scriptstyle]{#4}}#1\OrigSb{\mspace{-0.72mu}
			\cramped[\scriptstyle]{#3}}\OrigSp{\mspace{-0.72mu}
			\cramped[\scriptstyle]{#5}}  
}}
\NewDocumentCommand{\msbsp}{>{\SplitList{ , }} O{}}{
	\mth{\MultiSubSupscript#1}}

\makeatother

\DeclarePairedDelimiter{\Paren}{\lparen}{\rparen}
\DeclareDocumentCommand{\paren}{D(){}}{\Paren[\big]{#1}}
\DeclareDocumentCommand{\Bigparen}{D(){}}{\Paren[\Big]{#1}}
\DeclareDocumentCommand{\biggparen}{D(){}}{\Paren[\bigg]{#1}}
\DeclareDocumentCommand{\Biggparen}{D(){}}{\Paren[\Bigg]{#1}}

\DeclarePairedDelimiter{\lrangle}{\langle}{\rangle} 
\NewDocumentCommand{\tup}{D(){}}{\cramped{\lrangle{#1}}} 

\newcommand*{\N}{\mb{N}} 

\newcommand*{\Z}{\mb{Z}}

\NewDocumentCommand{\len}{o}{\IfNoValueTF{#1}{{\displaystyle{\ell}}}{{\displaystyle{\ell}}\sb{#1}}}

\def\udim{\msf{udim}}

\DeclarePairedDelimiter{\set}{\lbrace}{\rbrace}
\DeclareDocumentCommand{\vnoth}{o}{\IfNoValueTF{#1}{\varnothing}{{\varnothing\sb{[#1]}}}} 

\NewDocumentCommand{\Interv}{G{}g}{ \IfNoValueTF{#2}{\llbracket1\hspace{0.1em},\hspace{0.1em}#1\rrbracket}{\llbracket#1\hspace{0.1em},\hspace{0.1em}#2\rrbracket}}
\NewDocumentCommand{\irg}{>{\SplitList{ , }} O{n}}{\mth{\Interv#1}}

\makeatletter
\DeclareDocumentCommand{\with}{}{\hskip3pt \colon \hskip3pt}
\makeatother

\DeclareDocumentCommand{\cseq}{D(){}}{%
	\ensuremath{\foreach \x [count=\cpt]  in {#1} {%
			\ifnum\cpt=1%
			{\x}
			\else%
			\nsc\x\fi}}}

\NewDocumentCommand{\aarg}{O{-1.2}O{0.4}}{{\color[gray]{0.75}\rule[#1pt]{#2em}{1pt}}
}

\makeatletter

\DeclareDocumentCommand{\setminus}{O{0.36}}{\ensuremath{\raisebox{#1ex}{\scalebox{0.8}[0.625]{\mth{\backslash}}}}}
\makeatother

\DeclareMathOperator*{\tsum}{\textstyle{\sum}}
\NewDocumentCommand{\dsum}{o}{\IfNoValueTF{#1}{\sum}
	{\sum\limits\sbsp[#1]}}

\DeclarePairedDelimiter{\Abs}{\lvert}{\rvert}
\NewDocumentCommand{\abs}{o D(){}}{\IfNoValueTF{#1}{\Abs{#2}}{\Abs[#1]{#2}}}
\makeatletter
\let\origdeg\deg
\DeclareDocumentCommand{\deg}{o}{\IfNoValueTF{#1}{\origdeg}{\origdeg\sb{#1}}}
\NewDocumentCommand{\mdeg}{o}{\IfNoValueTF{#1}{\mrm{mdeg}}{\mrm{mdeg}\sb{#1}}}
\NewDocumentCommand{\Deg}{o}{\IfNoValueTF{#1}{\mrm{D}}{\mrm{Deg}\sb{#1}}}
\makeatother

\DeclareDocumentCommand{\rchi}{O{}}{\raisebox{0.18em}{\ensuremath{\upchi}}\sb{#1}} 

\makeatletter
\let\Origalpha\alpha
\DeclareDocumentCommand{\alpha}{o}{\IfNoValueTF{#1}{\Origalpha}{
		\Origalpha\sbsp[#1]}}

\DeclareDocumentCommand{\alphab}{o}{\IfNoValueTF{#1}{\bar(\Origalpha)}{
		\bar(\Origalpha)\sbsp[#1]}}

\let\Origbeta\beta
\DeclareDocumentCommand{\beta}{o}{\IfNoValueTF{#1}{\Origbeta}{
		\Origbeta\sbsp[#1]}}

\let\Origtau\tau
\DeclareDocumentCommand{\tau}{o}{\IfNoValueTF{#1}{\Origtau}{
		\Origtau\sbsp[#1]}}

\let\Origdelta\delta
\DeclareDocumentCommand{\delta}{o}{\IfNoValueTF{#1}{\Origdelta}{
		\Origdelta\sbsp[#1]}}

\DeclareDocumentCommand{\deltab}{o}{\IfNoValueTF{#1}{\bar(\Origdelta)}{
		\bar(\Origdelta)\sbsp[#1]}}

\let\OrigDelta\Delta
\DeclareDocumentCommand{\Delta}{o}{\IfNoValueTF{#1}{\OrigDelta}{
		\OrigDelta\sbsp[#1]}}

\DeclareDocumentCommand{\vDelta}{o}{\IfNoValueTF{#1}{\varDelta}{
		\varDelta\sbsp[#1]}}

\let\Origgamma\gamma
\DeclareDocumentCommand{\gamma}{o}{\IfNoValueTF{#1}{\Origgamma}{
		\Origgamma\sbsp[#1]}}

\let\OrigGamma\Gamma
\DeclareDocumentCommand{\Gamma}{o}{\IfNoValueTF{#1}{\OrigGamma}{
		\OrigGamma\sbsp[#1]}}

\let\Origkappa\kappa
\DeclareDocumentCommand{\kappa}{o}{\IfNoValueTF{#1}{\Origkappa}{
		\Origkappa\sbsp[#1]}}

\let\Origphi\phi
\DeclareDocumentCommand{\phi}{o}{\IfNoValueTF{#1}{\Origphi}{
		\Origphi\sbsp[#1]}}

\let\OrigPhi\Phi
\DeclareDocumentCommand{\Phi}{o}{\IfNoValueTF{#1}{\OrigPhi}{
		\OrigPhi\sbsp[#1]}}

\let\Origpsi\psi
\DeclareDocumentCommand{\psi}{o}{\IfNoValueTF{#1}{\Origpsi}{
		\Origpsi\sbsp[#1]}}

\let\OrigPsi\Psi
\DeclareDocumentCommand{\Psi}{o}{\IfNoValueTF{#1}{\OrigPsi}{
		\OrigPsi\sbsp[#1]}}

\let\Origxi\xi
\DeclareDocumentCommand{\xi}{o}{\IfNoValueTF{#1}{\Origxi}{
		\Origxi\sbsp[#1]}}

\let\Origzeta\zeta
\DeclareDocumentCommand{\zeta}{o}{\IfNoValueTF{#1}{\Origzeta}{
		\Origzeta\sbsp[#1]}}

\let\Origmu\mu
\DeclareDocumentCommand{\mu}{o}{\IfNoValueTF{#1}{\Origmu}{
		\Origmu\sbsp[#1]}}

\let\Orignu\nu
\DeclareDocumentCommand{\nu}{o}{\IfNoValueTF{#1}{\Orignu}{
		\Orignu\sbsp[#1]}}

\let\Origomega\omega
\DeclareDocumentCommand{\omega}{o}{\IfNoValueTF{#1}{\Origomega}{
		\Origomega\sbsp[#1]}}

\let\OrigOmega\Omega
\DeclareDocumentCommand{\Omega}{o}{\IfNoValueTF{#1}{\OrigOmega}{
		\OrigOmega\sbsp[#1]}}

\let\Origlambda\lambda
\DeclareDocumentCommand{\lambda}{o}{\IfNoValueTF{#1}{\Origlambda}{
		\Origlambda\sbsp[#1]}}

\let\OrigLambda\Lambda
\DeclareDocumentCommand{\Lambda}{o}{\IfNoValueTF{#1}{\OrigLambda}{
		\OrigLambda\sbsp[#1]}}

\let\Origpi\pi
\DeclareDocumentCommand{\pi}{o}{\IfNoValueTF{#1}{\Origpi}{
		\Origpi\sbsp[#1]}}

\let\OrigPi\Pi
\DeclareDocumentCommand{\Pi}{o}{\IfNoValueTF{#1}{\OrigPi}{
		\OrigPi\sbsp[#1]}}

\let\Origeta\eta
\DeclareDocumentCommand{\eta}{o}{\IfNoValueTF{#1}{\Origeta}{
		\Origeta\sbsp[#1]}}

\let\Origtheta\theta
\DeclareDocumentCommand{\theta}{o}{\IfNoValueTF{#1}{\Origtheta}{
		\Origtheta\sbsp[#1]}}

\let\OrigTheta\Theta
\DeclareDocumentCommand{\Theta}{o}{\IfNoValueTF{#1}{\OrigTheta}{
		\OrigTheta\sbsp[#1]}}

\let\Origsigma\sigma
\DeclareDocumentCommand{\sigma}{o}{\IfNoValueTF{#1}{\Origsigma}{
		\Origsigma\sbsp[#1]}}

\let\OrigSigma\Sigma
\DeclareDocumentCommand{\Sigma}{o}{\IfNoValueTF{#1}{\OrigSigma}{
		\OrigSigma\sbsp[#1]}}
\makeatother

\makeatletter	
\let\Origupsilon\upsilon
\DeclareDocumentCommand{\upsilon}{o}{\IfNoValueTF{#1}{\Origupsilon}{
		\Origupsilon\sbsp[#1]}}

\let\Origupsilonup\upsilonup
\DeclareDocumentCommand{\upsilonup}{o}{\IfNoValueTF{#1}{\Origupsilonup}{
		\Origupsilonup\sbsp[#1]}}
\makeatother

\DeclareDocumentCommand{\eps}{o}{\IfNoValueTF{#1}{\varepsilon}{
		\varepsilon\sbsp[#1]}}
\DeclareDocumentCommand{\upeps}{o}{\IfNoValueTF{#1}{\upvarepsilon}{
		\upvarepsilon\sbsp[#1]}}
\DeclareDocumentCommand{\vtheta}{o}{\IfNoValueTF{#1}{\vartheta}{
		\vartheta\sbsp[#1]}}
\DeclareDocumentCommand{\upvtheta}{o}{\IfNoValueTF{#1}{\upvartheta}{
		\upvartheta\sbsp[#1]}}
\DeclareDocumentCommand{\vphi}{o}{\IfNoValueTF{#1}{\varphi}{
		\varphi\sbsp[#1]}}
\DeclareDocumentCommand{\upvphi}{o}{\IfNoValueTF{#1}{\upvarphi}{
		\upvarphi\sbsp[#1]}}

\makeatletter
\let\Origrho\rho
\NewDocumentCommand{\processrho}{m g}{
	\IfNoValueTF{#2}{\Origrho\sb{\raisebox{-0.16ex}{\hskip-0.8pt\smth{#1}}}}{\Origrho\sb{#1}\sp{#2}}}
\DeclareDocumentCommand{\rho}{>{\SplitList{ , }} o}{\IfNoValueTF{#1}{\Origrho}{
		\processrho#1}}
\makeatother

\makeatletter
\let\Origuprho\uprho
\DeclareDocumentCommand{\uprho}{o}{\IfNoValueTF{#1}{\Origuprho}{
		\Origuprho\sbsp[#1]}}
\makeatother

\def\ring{R} 
\def\ringS{S}
\def\extringS{\mscr{S}} 
\NewDocumentCommand{\Ring}{ o }{\IfNoValueTF{#1}{\bsym{\msf{R}}}{
		\bsym{\msf{R}}\sb{#1}}}
\NewDocumentCommand{\RingS}{ o }{\IfNoValueTF{#1}{\mb{S}}{\mb{S}\sb{#1}}} 
\NewDocumentCommand{\ringZ}{ o }{\IfNoValueTF{#1}{\msf{Z}}{\msf{S}\sb{#1}}} 

\newcommand{\ACC}{{\sl ascending chain condition }}
\newcommand{\DCC}{{\slshape descending chain condition }}

\NewDocumentCommand{\bvar}{ D(){x}}{\bsym{#1}} 
\NewDocumentCommand{\term}{}{\msf{term}}  
\NewDocumentCommand{\varX}{>{\SplitList{:}} o }{\IfNoValueTF{#1}{\bsym{X}}{\bsym{X}\sbspscript#1}}
\NewDocumentCommand{\varY}{>{\SplitList{:}} o }{\IfNoValueTF{#1}{\bsym{Y}}{\bsym{Y}\sbspscript#1}} 

\NewDocumentCommand{\coef}{o}{\IfNoValueTF{#1}{\kwd{coef}}{\kwd{coef}\sb{#1}}}
\NewDocumentCommand{\lc}{o}{\IfNoValueTF{#1}{\kwd{lc}}{\kwd{lc}\sb{#1}}}
\NewDocumentCommand{\rlc}{o}{\IfNoValueTF{#1}{\kwd{rlc}}{\kwd{rlc}\sb{#1}}}
\NewDocumentCommand{\lcc}{o}{\IfNoValueTF{#1}{{\bsym{\mfr{c}}}}{{\bsym{\mfr{c}}}\sb{#1}}}
\NewDocumentCommand{\lt}{o}{\IfNoValueTF{#1}{\kwd{T}}{\kwd{T}\ssb{#1}}}

\NewDocumentCommand{\lm}{o}{\IfNoValueTF{#1}{\kwd{M}}{\kwd{M}\ssb{#1}}}
\NewDocumentCommand{\lex}{o}{\IfNoValueTF{#1}{\mathrel{<\sb{\kwd{lex}}}}{
		\mathrel{<\sb{\kwd{lex}\nsc #1}}}}
\NewDocumentCommand{\rlex}{o}{\IfNoValueTF{#1}{\mathrel{<\sb{\kwd{revlex}}}}{
		\mathrel{<\sb{\kwd{revlex}\nsc #1}}}}
\NewDocumentCommand{\rlexeq}{o}{\IfNoValueTF{#1}{\mathrel{\leq \sb{\kwd{revlex}}}}{
		\mathrel{\leq\sb{\kwd{revlex}\nsc #1}}}}

\NewDocumentCommand{\T}{>{\SplitList{:}} o }{\IfNoValueTF{#1}{\mb{T}}{\mb{T}\sbspscript#1}}
\NewDocumentCommand{\mcT}{ >{\SplitList{:}} o }{\IfNoValueTF{#1}{\mpzc{T}}{\mpzc{T}\sbspscript#1}}
\NewDocumentCommand{\cT}{ >{\SplitList{:}} o }{\IfNoValueTF{#1}{\mpzc{T}\sb{\mpzc{c}}}{\mpzc{T}\sbspscript[\sb{\mpzc{c}}]#1}}

\NewDocumentCommand{\spT}{ o }{\IfNoValueTF{#1}{\mb{T}\sprime}{\mb{T}\sbsp[#1,\prime]}}

\NewDocumentCommand{\essT}{>{\SplitList{:}} o }{\IfNoValueTF{#1}{\bsym{\mrm{T}}}{\bsym{\mrm{T}}\sbspscript#1}}
\NewDocumentCommand{\sT}{>{\SplitList{:}} o }{\IfNoValueTF{#1}{\msf{T}}{\msf{T}\sbspscript#1}}

\NewDocumentCommand{\q}{o}{\IfNoValueTF{#1}{\mfr{q}}{\mfr{q}\sb{#1}}}
\NewDocumentCommand{\qb}{o}{\IfNoValueTF{#1}{\bar(\mfr{q})}{\bar(\mfr{q})\sb{#1}}}

\def\ldenom(#1){#1\sb{\lft}}

\NewDocumentCommand{\reg}{o}{\IfNoValueTF{#1}{\mpzc{\msf{S}}}{\mpzc{\msf{S}}\sb{#1}}}
\NewDocumentCommand{\regC}{o}{\IfNoValueTF{#1}{\mscr{C}}{\mscr{C}\sb{#1}}}

\NewDocumentCommand{\Zc}{o}{\IfNoValueTF{#1}{\msf{Z}}{\msf{Z}\sb{#1}}}
\NewDocumentCommand{\zc}{o}{\IfNoValueTF{#1}{\mpzc{z}}{\mpzc{z}\sb{#1}}}

\DeclareDocumentCommand{\ide}{o}{\IfNoValueTF{#1}{\mpzc{1}}{\mpzc{1}\sb{#1}}}
\DeclareDocumentCommand{\idm}{d()}{\IfNoValueTF{#1}{\mathpzc{1}}{{\mathpzc{1}\sb{#1}}}}
\DeclareDocumentCommand{\id}{d()}{\IfNoValueTF{#1}{\mathsf{id}}{{\mathsf{id}\sb{#1}}}}

\NewDocumentCommand{\ann}{o}{\IfNoValueTF{#1}{\mathPzc{a}}{\mathPzc{a}\sb{#1}}}
\NewDocumentCommand{\annr}{o}{\IfNoValueTF{#1}{\mathPzc{a}'}{\mathPzc{a}\sb{#1}'}}
\NewDocumentCommand{\rann}{o}{\IfNoValueTF{#1}{\msf{r}\hskip-0.2em.\hskip-0.1em\mathPzc{ a}}{\msf{r}\hskip-0.2em.\hskip-0.1em\mathPzc{a}\sb{#1}}}

\def\sing{\bsym{\zetaup}}

\tikzset{morph/.style={>=stealth',every node/.style={inner sep=1pt,outer sep=1pt}},%
	nattransf/.style={-{stealth'}}, 
	diagm/.style={>=stealth',every node/.style={inner sep=2pt,outer sep=0.2pt},shorten <=-1.5pt, shorten >=-1.5pt},%
	grph/.style={>=latex',every node/.style={inner sep=1.5pt,outer sep=0.2pt},shorten <=-1.5pt, shorten >=-1.5pt},%
	rowstyle/.style={},
	every on chain/.append style={join},%
	every join/.style={->},%
	/seqStyle/.style ={>=stealth'}
}

\pgfkeys{
	/tipscale/.initial=0.8,
	/Arrowstyle/.code=\pgfkeys{/Arrowstyle/.initial=#1},/Arrowstyle/.default=morph,
	/Dir/.code=\pgfkeys{/Dir/.initial=#1},/Dir/.default=->,
	/Dist/.code=\pgfkeys{/Dist/.initial=#1},/Dist/.default=2.5,
	/Rowsep/.code=\pgfkeys{/Rowsep/.initial=#1},/Rowsep/.default=2,
	/Bline/.code=\pgfkeys{/Bline/.initial=#1}, /Bline/.default=-0.6ex,
	/Labsep/.code=\pgfkeys{/Labsep/.initial=#1},/Labsep/.default=0,
	/isymb/.code=\pgfkeys{/isymb/.initial=#1},/isymb/.default={\ensuremath{\sim}}, 
	/isymbb/.code=\pgfkeys{/isymb=#1,/isymbside/.initial=below},%
	/asymb/.code=\pgfkeys{/symb=#1,/symbside/.initial=above}, %
	/Arrowstyle,/Dir,/Dist,/Rowsep,/Bline,/Labsep
}

\makeatletter

\DeclareDocumentCommand{\to}{O{}}{
	\pgfkeys{/arrowstyle/.initial=\pgfkeysvalueof{/Arrowstyle},/style/.initial=,/bline/.initial=\pgfkeysvalueof{/Bline},
		/dir/.initial=\pgfkeysvalueof{/Dir},/dscale/.initial=1,/lab/.initial=,/speciallab/.initial=,
		/labsep/.initial=\pgfkeysvalueof{/Labsep},/side/.initial=above,/isymb/.initial=,/isymbside/.initial=above,
		/isymbsep/.initial=-1,/symb/.initial=,/symbsep/.initial=0.72,/symbside/.initial=below,/symbendsep/.initial=0.1,
		/lshift/.initial=0,/rshift/.initial=0,/lsep/.initial=1.8,/rsep/.initial=1.8,/slabsep/.initial=1,#1}
	\msp[\pgfkeysvalueof{/lsep}]
	\tikz[\pgfkeysvalueof{/arrowstyle},baseline=\pgfkeysvalueof{/bline},
	node distance=\pgfkeysvalueof{/dscale}*\pgfkeysvalueof{/Dist}em,trim left=(ito.east),trim right=(jto.west)]{
		\node(ito) [left]{}; 
		\node(jto) [right=of ito]  {}; 
		\path[\pgfkeysvalueof{/dir}] ($(ito.east)+(-\pgfkeysvalueof{/lshift}em,0em)$) 
		edge [\pgfkeysvalueof{/style}] node  [midway,\pgfkeysvalueof{/side}=\pgfkeysvalueof{/labsep}pt]{\ensuremath{\pgfkeysvalueof{/lab}}}
		node[\pgfkeysvalueof{/isymbside}=\pgfkeysvalueof{/isymbsep}pt]{\ensuremath{\scriptscriptstyle{
					\pgfkeysvalueof{/isymb}}}} ($(jto.west)+(\pgfkeysvalueof{/rshift}em,0em)$);
		\node (speciallabnode)  [inner sep=\pgfkeysvalueof{/slabsep}pt,above] at 
		($(ito.east)!0.5! (jto.west) $) {};
		\pgftext[at=(speciallabnode.north),bottom]{\ensuremath{\pgfkeysvalueof{/speciallab}}};
		\node (symb) [inner sep=\pgfkeysvalueof{/symbsep},\pgfkeysvalueof{/symbside}] at 
		($ (-\pgfkeysvalueof{/symbendsep}em,0) + (jto.west) $) {\ensuremath{\scriptstyle{\pgfkeysvalueof{/symb}}}};
	} \msp[\pgfkeysvalueof{/lsep}]}
\makeatother

\newcommand*{\Seq}[2][]{\pgfkeys{#1} 
	\ensuremath{\makebox{\ensuremath{#2}} 
	}\pgfkeys{/Arrowstyle,/Dir,/Dist,/Rowsep,/Bline,/Labsep}}

\theoremstyle{thmstyleone}%
\newtheorem{thm}{Theorem}[section] 
\newtheorem{prop}[thm]{Proposition}
\newtheorem{cor}[thm]{Corollary}
\newtheorem{lem}[thm]{Lemma}
\newtheorem*{thmn}{Theorem} 
\newtheorem*{stmn}{} 

\theoremstyle{thmstyletwo}%
\newtheorem{expl}[thm]{\bfseries Example}
\newtheorem{rem}[thm]{\bfseries Remark}

\theoremstyle{thmstylethree}%
\newtheorem{defn}[thm]{Definition} 
\newtheorem{prob}{Problem}[section]
\newtheorem{quest}[prob]{Question}

\numberwithin{equation}{section}

\newcommand*{\startthm}{\begin{thm}}
	\newcommand*{\stopthm}{\end{thm}}

\newcommand*{\startthmn}{\begin{thmn}}
	\newcommand*{\stopthmn}{\end{thmn}}

\newcommand*{\startcor}{\begin{cor}}
	\newcommand*{\stopcor}{\end{cor}}
\newcommand*{\startlem}{\begin{lem}}
	\newcommand*{\stoplem}{\end{lem}}

\newcommand*{\startprop}{\begin{prop}}
	\newcommand*{\stopprop}{\end{prop}} 

\newcommand*{\startstmn}{\begin{stmn}}
	\newcommand*{\stopstmn}{\end{stmn}}

\newcommand*{\startdef}{\begin{defn}}
	\newcommand*{\stopdef}{\end{defn}}

\newcommand*{\startexpl}{\begin{expl}}
	\newcommand*{\stopexpl}{\end{expl}}

\newcommand*{\startrem}{\begin{rem}}
	\newcommand*{\stoprem}{\end{rem}}

\newcommand*{\startnote}{\begin{rem}} 
	\newcommand*{\stopnote}{\end{rem}}

\newcommand*{\startquest}{\begin{quest}}
	\newcommand*{\stopquest}{\end{quest}}

\newcommand*{\startprob}{\begin{prob}}
	\newcommand*{\stopprob}{\end{prob}}

\NewDocumentCommand{\startpr}{}{\begin{proof}}
	\NewDocumentCommand{\stoppr}{}{\nobreak
\end{proof}}

\makeatletter
\DeclareDocumentCommand{\sbdsp}{m m g}{\IfNoValueTF{#3}{\sb{#1}\sp{#2}}{
		\sb{#1\nsc #2}\sp{#3}}}
\DeclareDocumentCommand{\bcoef}{>{\SplitList{ , }} D (){} }{\msf{C}\sbdsp#1}
\makeatother

\NewDocumentCommand{\Sn}{o}{\IfNoValueTF{#1}{\mrm{S}}{\mrm{S}\sb{#1}}}
\NewDocumentCommand{\bSn}{o}{\IfNoValueTF{#1}{\bar(\mrm{S})}{\bar(\mrm{S})\sb{#1}}}  
\NewDocumentCommand{\mrk}{O{}}{\bsym{m}\sb{#1}}  

\NewDocumentCommand{\upsilonx}{O{x}}{\bsym{\upsilon}\sb{#1}}  
\NewDocumentCommand{\splterm}{O{}}{\lambdaup\sbsp[#1]}  
\NewDocumentCommand{\ordn}{o}{\IfNoValueTF{#1}{\mfr{n}}{\mfr{n}\sb{#1}}}  

\NewDocumentCommand{\eext}{D(){\aarg} o}{\IfNoValueTF{#2}{\wt{#1}}{\wt{#1}\sb{#2}}}
\def\eextA{\eext(A)} 
    
\NewDocumentCommand{\isubext}{D(){\aarg} o}{\IfNoValueTF{#2}{#1}{#1\sb{(#2)}}}

\NewDocumentCommand{\dd}{o}{\IfNoValueTF{#1}{\msf{d}}{\msf{d}\sb{#1}}}  
\def\esspecial{essentially special }
\def\anesspecial{an essentially special }
\def\Esspecial{Essentially special }

\def\special{special }  

\let\citep\cite

\begin{document}
 \fontfamily{ptm}\fontseries{m}\selectfont
 	
	\title[Invariance of uniform dimension in subextensions]{On the uniform dimension of subextensions in skew polynomial rings} 
 	
	\author[B. Nguefack]{BertrandNguefack}
	\address{University of Yaounde~\MakeUppercase{\romannumeral 1},  P.O. Box 812, Yaounde, Cameroon}
	
	\email{Bertrand.Nguefack@facsciences-uy1.cm}






%
%
		
\begin{abstract}
This work  investigates the invariance of the non-necessarily finite  uniform dimension and related concepts 
for  subextensions in  skew polynomial rings \mbox{$ \mathbb{S}=R[ \mathbf{\mathrm{X}}; \mathbf{\alpha} , \mathbf{\delta} ]$}  of bijective type over a well-ordered set of variables. 
When the coefficient ring has enough uniform left ideals,  in the   commuting variables case we show that classical results on this topic for polynomial rings extend  to  subextensions of  
skew Laurent polynomial rings \mbox{$ \mathbb{S}=R[ \mathbf{\mathrm{X}} \sp{\pm1}; \mathbf{\alpha}]$}, generated over $R$   by any family of (standard) terms. The situation in the non-commuting variables context is more complex; easily formed polynomial-like subrings can  behave very oddly from the ambient ring.
 We provide easy examples of  a (semi)prime left Goldie skew polynomial ring of bijective type containing a monoid subring isomorphic to a free non-commutative polynomial ring. We then study the so-called subclass of \emph{essentially special subextensions} and  obtain for them the preservation of the uniform dimension and related concepts.   
\end{abstract}

\maketitle
	
\midsloppy


\section{Introduction}
Throughout this paper, every ring is associative and has an identity;   $\ring$ denotes a ground  ring (of coefficients);  a ring extension  of $\ring$ is alternatively termed as an \emph{$\ring$-ring}. An $\ring$-ring morphism is any ring 
map between $\ring$-rings that restricts to the identity map on the ground  ring;  a subextension of  an $\ring$-ring $\ringS$  is generally any  $\ring$-subring  of $\ringS$.   
The uniform (or Goldie) dimension $\udim(M)$ of a module $M$ (over some ring) is the supremum of those cardinal numbers $\mfr{N}$  such that $M$  contains an $\mfr{N}$-indexed direct sum of non-zero submodules \citep{Goodearl.Warfield2004,Dinh-et-al2006}. Except otherwise specified,  the uniform dimension of a ring $\ringS$ mean its  left uniform dimension denoted  and given by  $\udim(\ringS)=\udim({\sb{\ringS}S})$.  
Left Goldie conditions for a ring mean the  finiteness of the uniform dimension together with  the \ACC on left annihilators.  
It is well-known that a subring  
needs not have the same uniform dimension as 
the ambient ring. 
Starting with a classical paper by~\cite{Shock1972} on polynomial rings with finite right uniform dimension, the  invariance of the uniform dimension and  the preservation of Goldie conditions under ring extensions   has been   an active  research topic till date.  Due to the importance of polynomial rings,   and more generally, of skew polynomial rings  and related polynomial-like rings,  
Shock's result was immediately extended to group rings by~\cite{Wilkerson1973}, and   latter on, to (finitely) iterated Ore extensions of bijective type by \cite{Matczuk1995};  the case of iterated Ore extensions of injective type over semiprime left Goldie ring was dealt with by  \cite{LeroyMatczuk1995}, and that of skew PBW extensions is started by \cite{ArmandoReyes2014}. It is a well-established fact that the uniform dimension and Goldie conditions have important applications in ring theory itself \citep{Jategaonkar1972,CamCozz1973,Shock1974,CamGur1986,McConnellRobson1987,Kerr1990,Faith1996,Huynh-et-al2003} and references therein, as well as in  module theory \citep{Camillo1977,HJL1992,LeroyMatczuk2004}.      

Polynomial-like subrings of skew polynomial rings can behave very strangely from ordinary skew polynomial rings; however they arise naturally and in various specialized shapes in algebra and geometry. 
The present work  considers an ambient  $\ring$-ring $\RingS$  together with its subextensions, where  $\RingS$ is  
\startitem 
\item either a multivariate skew  Laurent polynomial ring  $ \ring[\varX\sp{\pm 1}; \bsym{\alpha}]$ of bijective type with  commuting variables,
\item  or a multivariate skew polynomial ring $\ring[\varX;\bsym{\alpha},\bsym{\delta}]$ of bijective type over  an arbitrary well-ordered set $\varX$ of independent variables.
\stopitem
Thus as  left $\ring$-module,  $\RingS$ is free  with basis given by the set $\T$ of standard terms: for the skew polynomial ring $\ring[\varX;\bsym{\alpha},\bsym{\delta}]$, $\T=\T(\varX)$ consist of all ordered product $x_1\dotsm x_n$ with $n\in\N$ and $x_1\leq x_2\leq \dotsm \leq x_n$ in $\varX$; and in the case of the skew Laurent polynomial ring in commuting variables, $\T=\T(\varX\sp{\pm 1})$ simply consists of all commutative words on the alphabet
$\varX\sp{\pm 1}= \varX \cup \varX\sp{-1}$ with $\varX\sp{-1}=\set{x\sp{-1}   \with x\in \varX}$.
We then address the following problem that arises  naturally.
\startprob
\label{prob:subext-skewpolring}
Let $A$ be a polynomial-like subextension of $\RingS$;  essentially one may restrict attention to the case that $A$ contains $\ring$ and is generated as left $\ring$-module by a subset of standard terms. 
Which conditions on the shape of $A$ will ensure that: $\udim(A)=\udim(\RingS)=\udim(\ring)$, while $A$ is left Goldie provided so is   $\ring$? 
\stopprob  
A specialization of Problem~\ref{prob:subext-skewpolring} to the commuting variables context follows.  
\startprob \label{prob:subext-skewLaurentpolring}
Let $A$ be any subextension of the skew Laurent polynomial ring $ \ring[\varX\sp{\pm 1}; \bsym{\alpha}]$, generated over $\ring$  by a family of standard terms from $\T(\varX\sp{\pm1})$. Does it hold that  
$\udim(A)=\udim(\ring)$? Is it true that $A$ is ((semi)prime) left Goldie precisely when so is $\ring$? 
\stopprob 

The feeling  about  Problem~\ref{prob:subext-skewLaurentpolring} is that the answer should be positive, at least when the ambient ring is simply a polynomial ring $\ring[\varX]$.  However even in  the latter case, it appears that one may not reach the expected positive answer by relying on the classical approach~\citep{Shock1972}. The  issue  becomes highly complex and very challenging in the non-commuting variables context of  Problem~\ref{prob:subext-skewpolring}.

\subsection*{Content and main results of the paper}

In section~\ref{sect:gen-ringext},  we review the invariance of the non-necessarily finite uniform dimension for arbitrary ring extensions. Indeed in the original paper by Shock and subsequent papers on the topic,  it was tacitly or explicitly assumed that the ground ring has finite uniform dimension.  For a non-zero left module $M$ over some ring, we say that $M$ has enough uniform submodules to mean that every non-zero submodule in $M$ contains a uniform submodule.   A  key step  towards the invariance of the uniform dimension is provided by 
Proposition~\ref{prop:attained-udim}: 
$M$ has enough uniform submodules   if and only if  $M$ contains an essential submodule which is a direct sum $\oplus\sb{i\in I} M\sb{s}$  of uniform submodules for some indexing set $I$;  in this case,   $\udim(M)$ equals the cardinal number  $\abs(I)$.
This key characterization allows  to provide a general and quick strategy  guaranteeing the invariance of  the  uniform dimension  under general ring extension (Proposition~\ref{prop:invariance-udim.ring-ext}). 
Specifically, we   introduce and   study a concept of  
\emph{nicely essential subring} $A$ of a given ring $\ringS$ (Definition~\ref{def:nicely-essext});   it results that $A$  and $\ringS$ share many properties in common including the invariance of the uniform dimension, the preservation of Goldie conditions and (semi)primeness; see Propositions~\ref{prop:udim-essringext}~and~\ref{prop:udim-essringext-(semi)prime}.  

Section~\ref{sect:subext-skew-Laurent-polring} focusses on the commuting variables setting since   this case is enlightening for the more complex non-commuting setting.   
Relying only on the general strategy investigated in  section~\ref{sect:gen-ringext},  we  solve  Problem~\ref{prob:subext-skewLaurentpolring} completely     in Theorems~\ref{thm:inv-udim-subext-skewLaurentPolring}~and~\ref{thm:Goldie-subext-skewLaurentPolring}. 

In Section~\ref{sect:essspsubext-skewpolring},
we address in full generality Problem~\ref{prob:subext-skewpolring} and formally introduce and investigate   \emph{\esspecial subextensions} of $\RingS$.   We then lift to this general setting the main classical results~\citep{Shock1972, Matczuk1995}, while allowing infinite uniform dimension. More precisely,   Theorems~\ref{thm:invaraince-uniform-univariate-case}~and~\ref{thm:inv-udim-subext-skewpolring} establish the invariance of the uniform dimension and the computation of left singular ideals for \esspecial subextensions $A$.
The main strategy is proving the lifting of uniform or essential left ideals  from the ground ring $\ring$ to $A$, generalizing classical techniques. The final theorem~\ref{thm:Goldie-subext-skewpolring} is an application to  (semi)prime left Goldie rings.  

Section~\ref{sect:esssp-subext-2variate-skewpolring} 
provides 
a family of \esspecial subextensions in the two-variate case, 
highlighting the fact that familiar subextensions in $\RingS$ generally fail to be skew polynomial rings in any meaningful way, while still having the same uniform dimension as $\RingS$.

Though our results are sharp for an ambient skew polynomial ring of bijective type, 
an obvious perspective 
would be to further extend the class of \tqt{specialized} subextensions of $\RingS$ for which the invariance of the uniform dimension holds.


\section{Invariance of the uniform dimension under general ring extensions} \label{sect:gen-ringext}

We consider general ring extensions $\ring \subset \ringS$ and provide  some  strategy allowing, on one hand   to quickly
compare the uniform dimension of both rings $\ring$  and $\ringS$, and on the other hand  to derive Goldie conditions for subrings of left Goldie rings. 

While $\N$ denotes the set of natural numbers including $0$, for $m,n\in\Z$ we write $\irg[m,n]$ for the range of all integers $k$ with $m\leq k\leq n$.  For every subset $U$ of in a left $\ringS$-module $M$,   the \emph{left annihilator} of $U$  in $\ringS$ is the left ideal 
\mthbox{\ann[\ringS](U)=\set{a\in \ringS \with a U=0}}; it is an ideal of $\ringS$  whenever $U$ is a  submodule in $M$. We usually omit the subscript in
$\ann[\ringS](U)$ when $\ringS$ is the base ring  or when it is implied by the context. The \emph{singular submodule}  of  $M$ is defined as  \mthbox{\sing(M)=\set{x\in M \with \ann(x) \text{ contains an essential  left ideal in } \ringS}}. 
We write \mthbox{\sing(\ringS)=\sing({\sb{\ringS}\ringS})}  for the left singular  ideal of the ring $\ringS$. The notation for  singular submodule is being  borrowed from \cite[\S2.2.4]{McConnellRobson1987}.

\startrem
The \ACC  and the \DCC  on left  annihilators  are  inherited by  subrings.   
Thus, a subring of a left Goldie ring is left Goldie if and only it has finite uniform dimension.
\stoprem
%

\subsection{Module with enough uniform submodules} 
\startdef  
We say  that  a non-zero module $M$ over some ring   \emph{has enough uniform submodules} provided every non-zero submodule in $M$ contains a uniform submodule. 
Likewise, a ring \emph{$\ringS$ has enough uniform left ideals} if every non-zero left ideal in $\ringS$ contains a uniform left ideal of $\ringS$. 
\stopdef 

Our first key result  characterizes  those modules whose uniform dimension is attained (in a strong sense), generalizing the finite dimensional case.  
\startprop\label{prop:attained-udim}
Let $M$  be any non-zero left module over some ring. Then   $M$ has enough uniform submodules   if and only if  $M$ contains an essential submodule which is a direct sum $\oplus\sb{i\in I} M\sb{s}$  of uniform submodules for some indexing set $I$.  In this case,  $\udim(M)$ coincides with the cardinal number  $\abs(I)$.
\stopprop 

\startpr
Start with the assumption that  $M$ has enough  uniform submodules.  Thus the family of uniform submodules of $M$ is non-empty and, by an obvious use of Zorn's lemma,  we may consider a maximal independent set $\mc{U}=\set{M\sb{i}\with i\in I}$  of  uniform submodules in $M$. Then the submodule \mthbox{\oplus\sb{i\in I} M\sb{i}}   must be  essential in $M$. 
Indeed assuming the contrary,   there will exist a non-zero submodule $N\subset  M$   with $(\oplus\sb{i\in I} M\sb{i})\cap N=0$, while by assumption  $N$ contains  some uniform submodule 
$ N'$. This yields a larger independent set $\mc{U}\cup \set{N'}$, contradicting the maximality of $\mc{U}$.  Hence as expected,   
$M$ contains an essential submodule which is a direct sum of uniform submodules.

Conversely, assume that  $M$ contains a direct sum $M\sb{I}=\oplus\sb{i\in I} M\sb{s}$  of uniform submodules such that $M\sb{I}$ is essential in $M$.    
Our gold is to prove  that every non-zero submodule in $M$ contains a uniform submodule and  $\udim(M)=\abs(I)$.
In the case  of a finite cardinal number $n\geq 1$, it is  well-known~\cite[Lemma 5.16,~Proposition 5.20]{Goodearl.Warfield2004} that: a module has uniform dimension $n$ precisely when it contains an essential submodule which is a direct sum of $n$ uniform submodules.  
Thus  we may continue  with the case that $\abs(I)$  is  infinite.
Let $N\subset M$ be any non-zero submodule; for every finite subset $\Lambda\subset I$  we  also let $M\sb{\Lambda}=\oplus\sb{i\in \Lambda}M\sb{i}$.  Since $M\sb{I}$ is essential in $M$, the module $N\cap M\sb{I}=N\cap (\oplus\sb{i\in I}M\sb{i})$ is non-zero and consequently, for some finite non-empty subset $\Lambda\subset I$ the module $N\cap M\sb{\Lambda}$ is non-zero as well.  But then (by the final cardinal case), $\udim(M\sb{\Lambda})=\abs(\Lambda)$ is finite and $\udim(N\cap M\sb{\Lambda}) \leq \abs(\Lambda)$ is also finite and non-zero. Thus, $N\cap M\sb{\Lambda}$ contains a uniform submodule.

Now assuming that $M$  contains  a direct sum   $\oplus\sb{s\in I'} M_s'$  of non-zero submodules  with $\abs(I') \geq \abs(I)$, 
our next aim is to  prove that $\abs(I')=\abs(I)$. Write $\aleph\sb{0}=\abs(\N) $ for the countable cardinal; for each finite number $n\in\N$,  
set $\mc{E}\sb{n} =\set{\Lambda \subset I \with  \abs(\Lambda)=n}$.
For every $s\in I'$, since  $M\sb{s}'$ is a non-zero submodule of $M$ while $M\sb{I}$ is an essential submodule of $M$, the module $M\sb{s}'\cap M\sb{I}$   must contain a non-zero  cyclic  submodule $N\sb{s}$  and there is some finite subset $\Lambda[s]\subset I$  with \mthbox{N\sb{s}\subset M\sb{\Lambda[s]}}. 
Next  fix  $1\leq n\in\N$. For every $n$-element subset $\Lambda\subset I$,   form the  subset:  
\mthbox{I\sb{\Lambda}' =\set{s\in I' \with N\sb{s}\subset M\sb{\Lambda}}}, thus the finite direct sum $M\sb{\Lambda}=\oplus\sb{i\in \Lambda} M\sb{i}$ of uniform submodules contains the direct sum $N\sb{{I\sb{\Lambda}'}}=\oplus\sb{s\in I\sb{\Lambda}'} N\sb{s}$ of non-zero cyclic submodules. 
But since $M\sb{\Lambda}$ cannot  contain a direct sum of more than $n$ non-zero submodules, it holds that $\abs(I\sb{\Lambda}')\leq n$.  
Since moreover each element $s\in I'$ is  by definition a member of $I\sb{\Lambda[s]}'$, we also have that 
$I'=\cup\sb{1\leq n\in\N}\cup\sb{\Lambda\in\mc{E}\sb{n}}
I\sb{\Lambda}'$. 
But for  each $1\leq n\in\N$, the set $I$ clearly embeds into each $\mc{E}\sb{n}$ and the latter embeds   into the $n$-fold Cartesian product $I^n$, and since $I$  is an infinite cardinal it follows that  $\abs(I)\leq \abs(\mc{E}\sb{n})\leq \abs(I^n)=\abs(I)$. Hence for each $1\leq n\in\N$  and every
$\Lambda\in \mc{E}\sb{n}$, we have shown that $\abs(I\sb{\Lambda}')\leq n$  and $ \abs(\mc{E}\sb{n})=\abs(I)$, whence the following:
\begin{equation*}
	\abs(I')=\abs({\cup\sb{1\leq n\in\N}\cup\sb{\Lambda\in\mc{E}\sb{n}}
		I\sb{\Lambda}'}) \leq \aleph\sb0\abs(I)=\abs(I). 
\end{equation*}
This shows that $\abs(I')=\abs(I)$, completing the proof of the proposition.
\stoppr 


\subsection{General strategies ensuring the invariance of the uniform dimension and related concepts}
The classical strategy successful for investigating the invariance of the finite uniform dimension for skew polynomial rings is summarized for arbitrary ring extensions by the following proposition, allowing infinite uniform dimension. 

\startprop\label{prop:invariance-udim.ring-ext} Let $A$  be any subring of an $\ring$-ring  such that $A\cap AL$ is essential as left  $A$-submodule in $AL$ for every non-zero left ideal $L$ in $\ring$, and for all  left ideals $L,L'$ with 
$L\cap L'=0$,  $AL \cap AL'=0$. 
Then $\udim(\ring)\leq \udim(A)$. Assume additionally that if $L$ is any  uniform  or  essential  left ideal in $\ring$, then so is the  left ideal $A\cap AL$ in  $A$.
Then, if $\ring$ has enough uniform left ideals
then so does  $A$ and  
$\udim(A)=\udim(\ring)$.
\stopprop

\startpr  To see that $\udim(\ring)\leq \udim(A)$, given any direct sum  $\oplus\sb{i\in I} L\sb{i}$ of non-zero left ideals in $\ring$, we quickly check  that one gets a direct sum  $\oplus\sb{i\in I} (A\cap AL\sb{i})$ of non-zero left ideals in $A$. By assumption,  each $A\cap AL\sb{i}$ is already a non-zero left ideal in $A$ for all $i\in I$. 
And if $i_1,\ldots, i_n$ are  pairwise distinct indices from $I$ for some integer $n\geq 2$, then for every  $s\in\irg[n]$, since $L\sb{i\sb s} \cap \sum\sb{k\neq s} L\sb{i\sb{k}} =0$,  the first assumption from  the proposition grants that $AL\sb{i\sb s} \cap \sum\sb{k\neq s} AL\sb{i\sb{k}}=0$. Thus the sum  $\sum\sb{i\in I} (A\cap AL\sb{i})$ is still a direct sum of non-zero left ideals in $A$. 

Next continue  by assuming that $\ring$ has enough uniform left ideals and,   if $L$  is a uniform (or an essential) left ideal in $\ring$  then so is the  left ideal $A\cap AL$ in $A$.   
Thus by virtue of Proposition~\ref{prop:attained-udim},  the uniform dimension of  $\ring$ is attained in the stronger sense that $\ring$ contains an essential left ideal arising as direct sum  $L=\oplus\sb{i\in I} L\sb{i}$ of uniform left ideals, with $\abs(I)=\udim(\ring)$.  But then, the left ideal $A\cap AL$  is essential in $A$ and it contains the direct sum   $\wt{L}=\oplus\sb{i\in I} (A \cap AL\sb{i})$ of uniform left ideals of $A$.  Let us show that $\wt{L}$ is also essential as left ideal in $A$.
Given any non-zero left ideal $J$ in $A$, since $ J\cap AL=J\cap A \cap AL\neq 0$, there is a finite non-empty subset $\Lambda\subset I$  such that $J\cap \oplus\sb{i\in \Lambda} AL\sb{i} \neq 0$. But  since by assumption each $A\cap AL\sb{i}$ is an essential submodule of the left $A$-module $AL\sb{i}$ for all $i\in\Lambda$, it follows that
$ \oplus\sb{i\in \Lambda}(A\cap AL\sb{i})$ is still an essential submodule of the left $A$-module $\oplus\sb{i\in \Lambda} AL\sb{i}$,
whence $J\cap \oplus\sb{i\in \Lambda}(A\cap AL\sb{i}) \neq 0$ as well. Thus we have shown that  $\wt{L}$ is an essential left ideal in $A$, arising as  $I$-indexed direct sum of uniform left ideals. Now another invocation of Proposition~\ref{prop:attained-udim} yields that
$\udim(A)=\abs(I)=\udim(\ring)$ and  $A$ has enough uniform left ideals. 
\stoppr

It is well-known that the uniform dimension is
preserved under left  localization with respect to a left Ore set of regular elements, \cite[Lemma 2.2.12]{McConnellRobson1987}.
A slightly more general situation is captured  by  the following  concept of nicely essential ring extension, as we will  show.

\startdef \label{def:nicely-essext} Let $L$ be a submodule of a left $\ring$-module $M$. We say that $L$ is \emph{nicely essential in $M$} or $M$ is a \emph{nicely essential extension}\footnote{An alternative terminology  could have been  \tqt{\emph{strongly essential submodule}}, but the latter already appears in few  places in the literature (for different purposes and) to mean a  submodule $L$  of an  $\ring$-module $M$  such that every  $I$-indexed  direct product of copies of $L$ is essential in the corresponding  $I$-indexed  direct product of copies of  $M$ for  all set $I$. {\slshape For a submodule $L$ of an $\ring$-module $M$, while the condition for $L$  to be \tqt{strongly essential}  in $M$ does not implies the property for  $L$  to be a nicely essential submodule in $M$, the latter is evidently  stronger than requiring that every direct sum of copies of $L$ is essential in the corresponding direct sum of copies of $M$}.}  
of $L$ provided for every finite  set   
$E\subset M$ of non-zero elements, there exists \mthbox{a\in  \ring}   with  
\mthbox{0\neq a x\in L} for all $x\in E$. 
If additionally it holds for every $x\in M$ that $cx\in L$ for some left regular element $c$ of $\ring$, then we call $L$ a 
\emph{strongly nicely essential} submodule of $M$, (lacking a more shorter terminology). 

By a  (\emph{strongly}) \emph{ nicely essential subring} of  a ring $\ringS$ we mean any subring  $A$  that is a (strongly) nicely essential  left $A$-submodule in  $\ringS$. In this case, we also call $\ringS$ a (\emph{strongly}) \emph{nicely essential  ring extension} of $A$.
\stopdef



\startprop \label{prop:udim-essringext} The following statements hold for a nicely essential subring  $A$ of a ring $\ringS$.    
\startitem
\item[\pmsf(a)] If $L,L'$ are left ideals in $A$  with $L\cap L'=0$  then $\ringS L \cap \ringS L'=0$;   \mthbox{\udim(\ringS) = \udim(A)}. The left singular ideal  $\sing(A)$  is a nicely essential left $A$-submodule in $\sing(\ringS)$  with $\sing(A)=A \cap \sing(\ringS)$. 
\item[\pmsf(b)] Let $L\subset \ringS$ be any non-zero left $A$-submodule.  Then the left ideal $A\cap L$  is uniform (or resp., essential) in $A$ if and only if so is the left ideal $\ringS L$ in $\ringS$.    
\stopitem
\stopprop

\startpr
Starting with the first claim, since $A$ is at least essential as left $A$-submodule in
$\ringS$, it is already granted that $\udim(\ringS)\leq \udim(A)$. Next assuming that $L,L'$ are two left ideals in $A$ with $L \cap L'=0$, 
we will show that $\ringS L \cap \ringS L'=0$.   Let $v\in \ringS L \cap \ringS L'$; then  $v =\sum\sbsp[i=1,n] v\sb{i}u\sb{i}=\sum\sbsp[i=1,n'] v\sb{i}'u\sb{i}'$ for some $n,n'\in\N$,  and \mthbox{v\sb{i},v\sb{j}'\in \ringS,\, u\sb{i} \in L, u\sb{j}' \in  L'} for $1\leq i\leq n $ and $1\leq j\leq n' $.  
The main assumption on  $\ringS$  applied to the finite set $E=\set{v , v\sb{1},\ldots, v\sb{n}, v\sb{1}',\ldots, v\sb{n'}'}$  yields some $a\in A$ such that for every $(i,j)\in\irg[n] \times {\irg[n']}$ we have:  
$a v\sb{i}, av\sb{j}'\in A,\,  a v  =\sum\sbsp[i=1,n] (av\sb{i})u\sb{i}=
\sum\sbsp[i=1,n'] (av\sb{i}')u\sb{i}' \in L\cap L'$ and 
$a v $ is non-zero whenever so is $ v$. 
But  since $ L\cap L'=0$,  we must have $av =0$ and $v =0$ as well. Hence,   
$ \ringS L \cap \ringS L'=0$.  This also grants (just as in the  first claim of Proposition~\ref{prop:invariance-udim.ring-ext}) that
$\udim(A)\leq \udim(\ringS)$;  whence the equality,
$\udim(\ringS) = \udim(A)$.

Turning to the claim about left singular ideals, 
the assumption that $A$ is a nicely essential left $A$-submodule of $\ringS$  already  grants  for every left ideal $J$ of $S$   that $A \cap J$ is still  a nicely essential left $A$-submodule of $J$. 
Also, $\sing(A)\subset A \cap \sing(\ringS)$ since $A$ is already an essential $A$-submodule in $\ringS$. For the converse inclusion, let $x\in A \cap \sing(\ringS)$ and $0\neq y\in A$. The definition of the left singular ideal $\sing(\ringS)$ yields some $b\in\ringS$  with $by\neq 0$ and  $by\cdot  x=0$, so that by the assumption on $\ringS$, there must also exist some $a\in A$ with $0\neq ab, ab\cdot y\in A$ while $(ab\cdot y)\cdot x=0$ and $x\in \sing(A)$. Whence, $A \cap \sing(\ringS)\subset \sing(A)$.

Next for statement~\pmsf(b), assuming that   $A\cap L$  is a uniform left  ideal  in $A$, let us show that   $\ringS L$ is uniform as left ideal in $\ringS$. So let 
$f=\sum\sbsp[i=1,m]f_iu_i$ and $g=\sum\sbsp[j=1,n]g_jv_j$  be  non-zero elements in $\ringS L$ with $0\neq m,n\in\N$,  
$f_i,g_j\in \ringS$  and $u_i,v_j\in L$  for $1\leq i\leq m$, $1\leq j\leq n$.
Applying the main assumption  to the finite subset $\set{f,g,f\sb i, g\sb j\, \with 1\leq i\leq m,\, 1\leq j\leq n}$  yields  some $v\in A$  with $ vf,vg,vf\sb{i},vg\sb{j}\in A$ (for all $(i,j)\in\irg[m]\times \irg[n]$) and $vf,vg\neq 0$.  In particular, the non-zero elements $vf$  and $vg$  already live in the uniform left ideal $A\cap L$ of $A$ and  there are two elements $a,b\in A $  with $0\neq avf=bvg $; whence  $\ringS L$ is uniform as  
left ideal of $\ringS$. Conversely, assuming that $\ringS L$ is uniform as left ideal in $\ringS$,  we quickly check  that   $A\cap L$ is uniform as left in $A$.  Since $A$ is already essential as left $A$-submodule in $\ringS$, $A\cap L$ must be non-zero. Given any non-zero elements $x,y \in A\cap L$, by assumption  there exist $c,d\in \ringS$  with $0\neq cx=dy$ and another application of the main assumption yields some 
$v\in A$  with $vc,vd \in A$ and $0\neq (v c)x  =(vd)y$, which shows that  $A\cap L$ is a uniform left ideal in $A$. 
Next, to complete  the proof of   statement~\pmsf(b),  it remains to address the claim  about essential left ideals.  Since  $A$ is already essential as left $A$-submodule in $\ringS$, 
if $A\cap L$ is essential as left ideal in $A$, then it is still essential as left $A$-submodule in $\ringS$ and so the left ideal  $\ringS L$ is essential  in $\ringS$. Conversely, supposing that $\ringS L$ is essential as left ideal of $S$, let us check that $A\cap L$ is essential as left ideal of  $A$. 
Given any non-zero $x\in A$, it holds for some $v\in \ringS$ that $0\neq vx=\sum\sbsp[i=1,m]v_iu_i\in \ringS L$ for some positive  $m\in\N$  and some
$(v_i,u_i)\in\ringS \times L$  with $v_iu_i\neq 0$ for $1\leq i\leq m$. The main assumption on $\ringS$ applied to the subset $\set{vx, v_i 
	\with 1\leq i\leq m}$ yields some $a\in A$ with
$0\neq avx, av_i\in A$ for all $i\in\irg[m]$, so that, $0\neq avx=
\sum\sbsp[i=1,m](av_i)u_i\in A\cap  L$ and $A\cap L$ is essential as left ideal in $A$ as expected. This completes the proof of the proposition.
\stoppr

We now characterize (semi)primeness under nicely essential extension. Here, 
\tqt{(semi)prime} is a short form for \tqt{prime (or respectively, semiprime)}.
\startprop 
\label{prop:udim-essringext-(semi)prime} 
Let $A$  be a nicely essential subring of a ring $\ringS$.  Then, if $A$ is (semi)prime then so is $\ringS$; if additionally $A$ is left Goldie then so is $\ringS$.  Conversely, if $\ringS$ is left Goldie, then so is $A$. If  $\ringS$ is (semi)prime left Goldie, then so is every  strongly nicely essential subring in $\ringS$.
\stopprop

\startpr 
Starting with the first statement  of the proposition and assuming that $A$ is prime, we will show that so is $\ringS$. Let $I,J$ be ideals in $\ringS$  with $IJ=0$ and $J\neq 0$. Since $A$ is essential as left $A$-submodule  in $\ringS$,  $J\cap A$ is a non-zero ideal of $A$ while $(I\cap A)(J\cap A)=0$, which  yields that $I\cap A=0$ and $I=0$. Hence, $\ringS$ is prime. The corresponding statement for semiprimeness follows by letting $I=J$.    Now, since a prime ring is semiprime while (by Goldie's theorem \cite[\S2.3.6]{McConnellRobson1987})  a semiprime  left Goldie ring is precisely a semiprime ring with finite uniform dimension and with zero left singular ideal, the second part of  the first statement of the proposition follows thanks to \mbox{Proposition~\ref{prop:udim-essringext}\pmsf(a)}

For the converse, suppose that $\ringS$ is left Goldie. Then since the \ACC on left annihilators is inherited by subrings,  $A$ is left Goldie by virtue of Proposition~\ref{prop:udim-essringext}\pmsf(a).  Now  assume the additional property that 
$A$ is a strongly nicely essential subring in $\ringS$.
Assuming that $\ringS$ is semiprime and given any ideal  $N$    in $A$ with $N\sp{2}=0$, let us show that $N$ is the zero ideal. 
The right annihilator ideal $I=\rann[A](N)$ is known to be essential as left ideal in $A$. Indeed for each non-zero left ideal $L$  in $A$, it  holds for some $k\in\set{0,1}$  that  $N\sp{k}L\neq 0$ and $N\sp{k+1}L=0$, hence 
$0\neq N\sp{k}L \subset L\cap \rann(N)$.   
It follows by \mbox{Proposition~\ref{prop:udim-essringext}\pmsf(b)}  that 
$\ringS I$ is also an essential left ideal of   $\ringS$. But then, because $\ringS$ is   semiprime left non-singular    with finite uniform dimension, we know 
(thanks to \cite[Proposition~2.3.5]{McConnellRobson1987}) that
$\ringS I$ contains some regular element $c$. Write $c=u_1a_1+\dotsm+u_na_n$
with $u_1,\ldots,u_n\in \ringS$  and  $a_1,\ldots, a_n\in I$.
By the additional assumption on $A$, there are left regular elements  $c_1,\ldots,c_n$ of $A$ with
$c_n u_n,c_{n-1} c_n u_{n-1}, \ldots, c_1\dotsm c_n  u_1\in A$. 
In particular,  the product $c'=c_1c_2\dotsm c_n c$ is necessarily a left regular element of $A$ that now lives in $I$; thus $Nc'=0$ which yields $N=0$ as expected.  Finally, assuming that $\ringS$ is prime, we wish to show that $A$ is prime as well.  We know (by Goldie's theorem \cite[\S2.3.6]{McConnellRobson1987}) that the set $\regC(\ringS)$ of all regular elements of $\ringS$ is a left Ore set and the associated left Goldie quotient $Q$ of $\ringS$ is a simple ring.  The following arguments generalize  those proving that $\ringS$ is prime from the assumption that $Q$ is simple.  
Call an element $c\in \ringS$, \emph{left $A$-regular} provided the left annihilator $\ann[A](c) =\set{a\in A \with ac=0}$ vanishes. Let $I,J$ be any non-zero ideals in $A$. Since $Q$ is simple, we have $QJQ=Q$ and 
$1=\sum\sbsp[i=1,n] d\sbsp[i,-1] u\sb{i}a\sb{i} s\sbsp[i,-1] v\sb{i}$ for some positive $n\in\N$  and some $u_i,v_i\in \ringS$,  
$d_i,s_i\in \regC(\ringS)$ and $a\sb{i}\in J$ for $1\leq i\leq n$. Letting $d\in\regC(\ringS)$ be a common left denominator for the $d\sb{i}$'s, we get that $d=d\times 1=\sum\sbsp[i=1,n]   u\sb{i}' a\sb{i} s\sbsp[i,-1] v\sb{i}$ for some $u\sb{i}'\in \ringS$, $1\leq i\leq n$. Now applying the additional assumption for $A$, there are left regular elements $c_1,\ldots, c_n$ of $A$  with
$c_n u_n',c_{n-1} c_n u_{n-1}', \ldots, c_1\dotsm c_n  u_1'\in A$.
But then, since moreover $d$ is a regular element of the larger ring $\ringS$,  the product $c=c_1 c_2\dotsm c_n d$  is necessarily a left $A$-regular element  which now lives in $JQ$. It follows in particular that $IJQ\supset Ic\neq 0$. Hence $IJ\neq 0$ and $A$ is prime. This completes the proof of the proposition.      
\stoppr 

\subsection{Easy examples illustrating the invariance of the uniform dimension or its failure}   

\startexpl \label{expl1}
Fix $1\leq q\in\N$ and two $q$-roots of unity  $\varepsilonup,\varepsilonup'$  in the center of $\ring$.  Consider the skew polynomial ring $\RingS=\ring[x,y][z;\alpha]$  where $\alpha$ is the
$\ring$-ring automorphism with $\alpha(x)=\varepsilonup y$ and $\alpha(y)=\varepsilonup' x$; clearly $\alpha[,2q]$ is the identity map.   

The  subextension $A\subset \RingS$ generated by $x,z$ is free as left $\ring$-module and its left $\ring$-basis is the following  subset of terms:
\startcent $\sT=\set{x\sp{m}y\sp{n}z\sp{p} \with 
	m,n,p\in \N \text{ and, if } n\neq 0 \text{ then } p\neq 0}$.
\stopcent 
The ring $A$ cannot be presented as skew polynomial ring (even if one allows an infinite set of variables). 
But interestingly, the term $z\sp{2q} \in \sT$ is a central element of $\RingS$ with  $z\sp{2q} \RingS \subset A$. Thus  $A$ is a nicely essential subring in $\RingS$, in particular Proposition~\ref{prop:udim-essringext} grants that
$\udim(A)=\udim(\RingS)$.  
%
\stopexpl

\startexpl\label{expl2} {\slshape A three-variate (semi)prime left Goldie graded skew polynomial ring of bijective type over $\ring$ admitting a monoid $\ring$-subring that is isomorphic to a free non-commutative polynomial ring.}

Modifying the previous example, consider the graded Ore extension   $\RingS=\ring[x,y][z;\alpha]$  where this time $\alpha$ is the
$\ring$-ring automorphism with $\alpha(x)= y$ and $\alpha(y)=x+y$.
Thus we know that $\RingS$ is (semi)prime left Goldie precisely when so is the coefficient  $\ring$.  Letting $(a_n)_n$  denote the Fibonacci series with $a_0=0,a_1=1$ and $a\sb{n+1}=a\sb{n-1}+a\sb{n}$ for all $n\geq 1$, we have: 
\begin{equation*}
	\alpha[,n](x)=a\sb{n-1}x+a\sb{n}y
	\text{ and } \alpha[,n](y)=\alpha[,n+1](x)=a\sb{n}x+a\sb{n+1}y
	\text{ for } 1\leq n\in\N.
\end{equation*} 
It easily follows that $\alpha$ will have finite order precisely when $\ring$ has a non-zero characteristic.  Assume that  $\ring$  has characteristic $0$, and  let $u=xz,\, v=x\sp{2}z$.  

Then the terms $u$  and $v$ lie in the monoid $\regC(\RingS)$ of all regular elements of $\RingS$ since the variables $x,y,z$ are regular in $\RingS$.  We claim that the submonoid  $\Gamma=\tup(u,v)\subset \regC(\RingS)$ generated by $u$  and $v$  is  the free  monoid on two non-commuting variables.   Indeed given any $k$-length word  $w=u_1u_2\dotsm u_k$    with $u_i\in\set{u,v}$ for $1\leq i\leq k$,  we have:  $\deg[z](u_i)=1$ and 
$\deg[z](w)=k$;  thus two words in the alphabet $\set{u,v}$, with different lengths, cannot be equal as elements of $\Gamma$.  
Next, there are exactly four  distinct elements of $\Gamma$ having length $2$: 
\begin{equation*}
	u^2=xyz\sp{2},\, uv=xy\sp{2}z\sp{2}, vu=x\sp{2}yz\sp{2} \text{ and } v^2=x\sp{2}y\sp{2}z\sp{2}.
\end{equation*}
The $8$ three-length words in terms of $u$  and $v$ are also pairwise distinct  in $\Gamma$:
\begin{align*}
	u^3 &=xy(x+y)z\sp{3}, & v^3 &=x\sp{2}y\sp{2}(x+y)\sp{2}z\sp{3};\\
	u\sp{2}v &=xy(x+y)\sp{2}z\sp{3}, & v^2u &=x\sp{2}y\sp{2}(x+y)z\sp{2}\\
	uvu &=xy\sp{2}(x+y)z\sp{3}, &  
	vuv &=x\sp{2}y(x+y)\sp{2}z\sp{3},\\
	uv^2 &=xy\sp{2}(x+y)\sp{2}z\sp{2}, &
	vu\sp{2}&= x\sp{2}y(x+y)z\sp{3}.
\end{align*}
Let  $\Gamma[\alpha]$ be the submonoid of $\regC(\RingS)$ generated by the elements \mthbox{\alpha[,n](x),\, 1\leq n\in \N}.  
For every $k$-length word  $w=u_1u_2\dotsm u_k$ with $k\geq 3$ and $u_1,\ldots, u_k\in\set{u,v}$, 
\startcent  there exists  some $\tau\in\Gamma[\alpha]$  with $w=x\tau z\sp{k}$ if $u_1=u$ or  $w=x\sp{2}\tau z\sp{k}$ if $u_1=v$.  
\stopcent 
Since moreover, $\gcd(x,\tau)=1$  for all $\tau\in\Gamma[\alpha]$,
these last observations together with the already computed $2$-length words in the alphabet $\set{u,v}$  show that for all $k,l\in\N$ and terms   $u_1,\ldots,u_k,v_1,\ldots,v_l\in\set{u,v}$, the products  
$u_1\dotsm u_k$ and $v_1\dotsm v_l$ are equal in $\RingS$ precisely when $k=l$  and $u_i=v_i$ for all $1\leq i\leq k$. Hence the submonoid $\Gamma= \tup(u,v)\subset \regC(\RingS)$ is  the free monoid on two non-commuting variables.  Consequently,  the monoid subring $\ring[u,v]=\ring\Gamma\subset \RingS$ is a two-variate free non-commutative polynomial ring, and it has infinite left or right uniform  dimension (even when the ambient ring is (semi)prime left Goldie). 
\stopexpl  

The subring $\ring[u,v]$ in the previous example is a very easy example of a monoid $\ring$-subextension of  $\RingS$ generated by a subset of standard terms that obviously fails to be \emph{essentially special} in $\RingS$  (as introduced by Definition~\ref{def:special-subext}).  This shows that the main results of  section~\ref{sect:udim-esssp-subext-skewpolring} are also sharp. 

For (semi)prime left Goldie skew polynomial rings of injective type over $\ring$, it is even more easy to form monoid $\ring$-subextensions with infinite uniform dimension. 

\startexpl\label{expl3} {\slshape A family of two-variate (semi)prime left Goldie skew polynomial rings  of injective  type over $\ring$ containing a two-variate free non-commutative polynomial ring.} 

Fix an integer $e\geq 2$ and consider the monoid ring $\RingS=\ring\Gamma[e]$  where $\Gamma[e]=\tup(x,y\with yx=x\sp{e}y)$ is the quotient of the free monoid on two generators $x,y$ modulo the relation that $yx=x\sp{e}y$. Alternatively,  $\RingS=\ring[x,y;\alpha]=\ring[x][y;\alpha]$ is the graded skew polynomial ring  where \Seq{\alpha: \ring[x] \to \ring[x]} is the $\ring$-ring injective map with $\alpha(x)=x\sp{e}$.
Thus when the coefficient ring $\ring$ is semiprime left Goldie,  we know that $\RingS$ is semiprime left Goldie and $\udim(\RingS)=\udim(\ring )$~\cite[Theorems 3.4, 3.8]{LeroyMatczuk1995}. 
Next, we claim that the terms $xy,x\sp{2}y$ freely generate a non-commutative submonoid $\tup(xy,x\sp{2}y) \subset \Gamma[e]$.
Given any $k$-length word $w= u_1\dotsm u_k$ with $u_1\ldots u_k \in\set{ xy, x\sp{2}y}$, one  observes that $\deg[y](u_i)=1$ and 
$\deg[y](w)=k$, so that, two words in the alphabet $\set{xy,x\sp{2}y}$  with different lengths  cannot be equal as elements of $\Gamma[e]$.  Next, $\deg[x](w) \equiv \deg[x](u_1) \mod e$ (while $1=\deg[x](xy)\neq \deg[x](x\sp{2}y) =2$).
These two observations  show that for all $k,l\in\N$ and terms  $u_1,\ldots,u_k,\, v_1,\ldots,v_l\in\set{xy,x\sp{2}y}$, the products  
$u_1\dotsm u_k$ and $v_1\dotsm v_l$ are equal in $\Gamma[e]$ precisely when $k=l$  and $u_i=v_i$ for all $1\leq i\leq k$; which yields our claim. 
Consequently, the monoid subring $\ring[xy,x\sp{2}y]\subset \ring\Gamma[e]$  is a two-variate free non-commutative polynomial ring over $\ring$;  
and $\ring[xy,x\sp{2}y]$  has infinite left or right uniform  dimension. 
\stopexpl



\section{Subextensions in skew Laurent polynomial rings}
\label{sect:subext-skew-Laurent-polring}

We fix an ambient skew Laurent polynomial ring $\ring[\varX\sp{\pm1};\bsym{\alpha}]$, always assumed of bijective type; here   $\varX$ is an  arbitrary non-empty set of freely commuting variables over $\ring$, 
$\varX\sp{\pm1}=\varX\cup \varX\sp{-1}$ with $\varX\sp{-1}=\set{x\sp{-1} \with x\in \varX}$.
\startalignl For each $x\in \varX$, the associated conjugation map 
\Seq{\alpha[x]: \ring\to\ring} is a ring automorphism with $xa=\alpha(a)x$  and $x\sp{-1}a=\alpha[x,-1](a)x\sp{-1}$  for all $a\in \ring$.
\stopalignl
For every $\lambda=x\sbsp[1,k_1]\dotsm x\sbsp[n,k_n]$  with $n\in\N$,  $k_1,\ldots,k_n\in\N$ and $x\sb1,\ldots x\sb{n}\in\varX$, the associated conjugation map is the composite automorphism
\begin{equation*}
	\alpha[\lambda]=\alpha[{x_{1}},k_1]\cdot\dotsm  \alpha[{x_{n}},k_n] :\ring \to[/isymb] \ring.
\end{equation*}
The skew Laurent polynomial ring   $\ring[\varX\sp{\pm1};\bsym{\alpha}]$  may be alternatively termed as the skew monoid ring extension of $\ring$ by the free abelian
group $\mcT(\varX\sp{\pm1})$ over $\varX$:  
\begin{equation*}
	\mcT(\varX\sp{\pm1})
	=\set{x\sbsp[1,k_1]\dotsm x\sbsp[n,k_n] \with n\in \N,\,  k_1,\ldots,k_n\in\Z
		\text{ and } x\sb1,\ldots, x\sb{n} \in \varX}.
\end{equation*}
Thus, a subextension of 
$\ring[\varX\sp{\pm1};\bsym{\alpha}]$ (generated over $\ring$ by a family of terms from $\T(\varX\sp{\pm1})$)  is precisely a skew monoid ring extension  of $\ring$  by a submonoid $\sT$ in a free abelian group.  
Every skew polynomial ring of bijective type and  in commuting variables embeds into a skew Laurent polynomial ring; see Proposition~\ref{prop:localizedmonoidring} below. 

For the investigation of left singular ideals and semiprimeness of subextensions in $\ring[\varX\sp{\pm1};\bsym{\alpha}]$,   it shall be necessary to consider a  well-ordering \tqt{\mth{\leq}} on  $\varX$, which is always possible (at least by axiom of choice). Then the latter  is extended to a linear ordering on 
$\varX\sp{\pm1}$  by setting $y\sp{-1}<x\sp{-1}$ for all 
all variables $x,y\in \varX$ with $x<y$. 
The  induced  reverse lexicographic ordering on $\mcT(\varX\sp{\pm1})$ is defined as it follows. For every natural number $n\geq 1$ and every $n$-subset or $n$-tuple $\bsym{x}$ 
formed by variables $x\sb{1}<\dots<x\sb{n}$ in $(\varX,\leq)$,   
for all $n$-tuples $\bsym{k}=(k_1,\ldots,k_n), \bsym{k}'=(k_1',\ldots,k_n')\in\Z\sp{n}$, 
setting $\bsym{x}\sp{\bsym{k}}=x\sbsp[1,k_1]\dotsm x\sbsp[n,k_n]$
we have: 
\begin{equation*}
	\bsym{x}\sp{\bsym{k}} < \bsym{x}\sp{\bsym{k}'} \iff \, \exists   j\in\irg[1,n],\,    k_j<k_j' \text{  and }  k_i=k_i' \text{ for } j<i\leq n.
\end{equation*} 
We  get a corresponding \emph{leading term function} 
\Seq{\lt:\ring[\varX\sp{\pm1}; \bsym{\alpha}] \to  \mcT(\varX\sp{\pm1})} where for each 
$f=a_1 \tau[1]+\dotsm+a_p\tau[p]$ with $1\leq p\in \N$,
$ 0\neq a_1,\ldots, a_p\in \ring$ and $\tau[1]<\dotsm<\tau[p]$ in  $\mcT(\varX\sp{\pm1})$, we set $\lt(f)=\tau[p]$ and the associated leading coefficient
is $\lc(f)=a_p$.
For convenience, we usually set: $\lc(0)=0=\lt(0)$ and $0<\tau$ for all $\tau\in \mcT(\varX\sp{\pm1})$.   

\startrem \label{rem:term-fct-skwelaurentpolring}  
We get  a linearly ordered free abelian group 
\mthbox{(\mcT(\varX\sp{\pm1}),\leq)}:  the linear-ordering on $\mcT(\varX\sp{\pm1})$  is compatible with the multiplication in the sense that,  
\begin{equation*}
	\text{for all } 	\lambda,\tau,\tau'\in \mcT(\varX\sp{\pm1}),
	\text{ if } \tau<\tau' \text{ then }\lambda\tau <\lambda\tau'.
\end{equation*} 
Thus, the  leading term function \Seq{\lt:\ring[\varX\sp{\pm1}; \bsym{\alpha}] \to  \mcT(\varX\sp{\pm1})} is a pseudo valuation on $\ring[\varX\sp{\pm1}; \bsym{\alpha}]$; more precisely  
for  all \mthbox{0\neq f,g\in \ring[\sT,\bsym{\alpha}]}, 
it holds that     
\mthbox{fg=\lc(f)\alpha[\lt(f)](\lc(g))\lt(f)\lt(g) +h} for some polynomial $h$   with $\lt(h)<\lt(f)\lt(g)$.
\stoprem 

%
%
%
%

By  a \emph{left normal element} in a monoid $\Gamma$  we mean any $\upsilon \in \Gamma$   with $\upsilon\Gamma\subset \Gamma \upsilon$.
We include the following proposition which should be well-known.
\startprop \label{prop:localizedmonoidring}
Let $\ringS=\ring[\Gamma;\bsym{\alpha}]$  be a skew monoid ring of injective type.
\startitem
\item[\pmsf(a)] If $\Gamma$ is left or right cancellative, then each  $\lambda\in \Gamma$ is left (or resp., right)  regular as element of $\ringS$. 
\item[\pmsf(b)] Assume that $\Gamma$ is cancellative. 
Then  the subset $\msf{Z}$  of left normal terms in $\Gamma$ is a left Ore set  of  regular elements for $\ringS$.  If each conjugation map $\alpha[\upsilon]$ for $\upsilon\in\msf{Z}$ is bijective,   then the left ring of fractions $\msf{Z}\sp{-1} \ringS$ is a skew monoid ring $\ring[G;\bsym{\alpha}]$  with $G=\msf{Z}\sp{-1}\Gamma=\set{\upsilon[,-1]\lambda \with (\upsilon, \lambda)\in \msf{Z} \times \Gamma}$  and $\alpha[\upsilon\sp{-1}\lambda]=\alpha[\upsilon,-1]\circ \alpha[\lambda]$ for all $(\upsilon,\lambda) \in \msf{Z} \times  \Gamma$.   
\stopitem 
\stopprop 

\startpr
Part~\pmsf(a) of the proposition easily follows by a direct verification. Turning to part~\pmsf(b), assume that $\Gamma$ is cancellative. Then it is immediate  by \pmsf(a) that every element of  $\Gamma$ is regular as element of $\ringS=\ring[\Gamma;\bsym{\alpha}]$.  Let $\msf{Z}$ consists of all left normal elements in $\Gamma$ and observe that $\msf{Z}$ is already closed under multiplication. Clearly, for  every $\upsilon \in\msf{Z}$ and $f\in \ringS$, it holds for some $f'\in \ringS $ that $\upsilon f=f'\upsilon \in  \ringS   \upsilon \cap \msf{Z} f$. We conclude that    $\msf{Z}$ is a left Ore set  of regular elements  for  $\ringS $; in  particular one can form the left ring of fractions $\msf{Z}\sp{-1}\ringS$.   For the last claim in \pmsf(b), we now suppose for every $\upsilon\in \msf{Z}$ that each conjugation map \Seq{\alpha[\upsilon] : \ring \to \ring} is bijective.
Given each $\upsilon \in \msf{Z}$,  for every term $\lambda\in \Gamma$ there is a unique term
$\lambda[\upsilon] $  with $\upsilon \lambda =\lambda[\upsilon] \upsilon$; and for every scalar $a\in \ring$ we have $\upsilon a=\alpha[\upsilon](a)\upsilon$.
In particular, it holds in $\msf{Z}\sp{-1}\ringS$ that:
\begin{equation*}
	\lambda \upsilon[,-1] =\upsilon[,-1]\lambda[\upsilon] 
	\text{ and }   \upsilon[,-1] a =\alpha[\upsilon,-1](a)\upsilon[,-1]   
	\text{ for all } a\in \ring.
\end{equation*}
We obtain in $\msf{Z}\sp{-1}\ringS$ a submonoid  $G=\msf{Z}\sp{-1}\Gamma=\set{\upsilon[,-1]\lambda \with (\upsilon, \lambda)\in \msf{Z} \times \Gamma}$ which contains $\Gamma$ as submonoid; 
and $G$ is a left $\ring$-generating set for $\msf{Z}\sp{-1}\ringS$.
To complete the proof, it remains to see that $G$ freely generates the left $\ring$-module 
$\msf{Z}\sp{-1}\ringS$. Thus  let   $f=\sum\sbsp[i=1,n] a\sb{i} \upsilon[i,-1]\lambda[i]=0$  with 
$0\neq n\in\N$, $a\in\ring$ and  $(\upsilon[i],\lambda[i]) \in \msf{Z} \times \Gamma$ such that the $\upsilon[i,-1]\lambda[i]$'s  are pairwise distinct elements of $\msf{Z}\sp{-1}\ringS$. We must check that the $a\sb{i}$'s are all zero.   Since $\msf{Z}$ consists of left normal elements of $\Gamma$, the product-term 
$\upsilon=\upsilon[1]\dotsm  \upsilon[n]$ readily serves as  common left denominator for the  $\upsilon[i,-1]\lambda[i]$'s and we have:
$\upsilon\cdot \upsilon[i,-1]\lambda[i] = \lambda[i]'$ for some pairwise distinct $\lambda[i]'\in \Gamma$, $1\leq i\leq n$.  It follows that:
$0=\upsilon f=\sum\sbsp[i=1,n] \alpha[\upsilon](a\sb{i})  \lambda[i]'$, showing that 
$\alpha[\upsilon](a\sb{i})=0$ and $a_i=0$ for all $i$.
Hence $G$ is a left $\ring$-basis for $\msf{Z}\sp{-1}\ringS$ and the latter arises as skew monoid ring $\ring[G;\bsym{\alpha}]$. 
\stoppr 

\startthm\label{thm:inv-udim-subext-skewLaurentPolring}
Let $A=\ring[\sT;\bsym{\alpha}]$  be any subextension in a skew Laurent polynomial ring of bijective type over $\ring$. 
Then $\sT$  is a denominator set for $A$,  the associated ring of fractions is  a skew Laurent polynomial ring of bijective type $\RingS=\ring[\varX\sp{\pm1};\bsym{\alpha}]$, and 
$\udim(A)=\udim(\RingS)=\udim(\ring[\varX;\bsym{\alpha}])$, $\sing(A)=\sing(\ring)\sT$. 
If $\ring$ has enough uniform left ideals, then so does  $A$  and 
$\udim(A)=\udim(\ring)$.
\stopthm

\startpr In view of the structure theorem for free modules over principal ideal 
domains~\cite[Thm IV. 6.1]{Hungerford1974},   every subgroup in a free abelian group is still a free abelian group.
The definition of $A=\ring[\sT;\bsym{\alpha}]$ yields that $\sT$ is a submonoid in some free abelian group. We can therefore form a free abelian group $G$ containing $\sT$ as submonoid such that $G=\sT\sT\sp{-1}=\set{\lambda\tau[,-1]\, \with \lambda,\tau \in \sT}$, and  we may let $\varX$  be a basis of $G$.  
One may view $\varX$ just as a set of feely commuting  variables over $\ring$, and write $G=\mcT(\varX\sp{\pm1})$. Naturally setting $\alpha[\lambda\tau\sp{-1}] =\alpha[\lambda]\alpha[\tau,-1]$  for  all  $\lambda,\tau\in \sT$,   one can now  form the graded skew polynomial ring $\ringS=\ring[\varX;\bsym{\alpha}]$, of bijective type  and with commuting variables; recall a natural two-sided $\ring$-basis of $\ringS$ is given by the set  $\T(\varX)$ (of standard) terms.  Likewise we get the  skew Laurent polynomial ring of bijective type
$\RingS=\ring[\varX\sp{\pm1};\bsym{\alpha}]$.
Obviously as clarified by Proposition~\ref{prop:localizedmonoidring}, $\RingS$ is a ring of   
fractions for $\ringS=\ring[\varX;\bsym{\alpha}]$  with respect to $\varX$, as well as for the ring $A=\ring[\sT;\bsym{\alpha}]$   with respect to $\sT$ (since $G=\mcT(\varX\sp{\pm1})=\sT\sT\sp{-1}$).
In particular, both rings $\ringS$ and $A$ are nicely essential subextensions in $\RingS$;
more precisely,
\startalignl
for finite subsets $Y\subset \T(\varX)$  and  $Z\subset \sT$, there are terms
$\lambda\in \sT$ and $\tau\in \T(\varX)$  with
$\lambda Y\subset \sT$ and $\tau Z\subset \T(\varX)$.
\stopalignl
Then  Proposition~\ref{prop:udim-essringext}  grants that  
\mthbox{{\udim(A)=\udim(\RingS)=\udim (\ringS)}},  \mthbox{{\sing(A)=A\cap \sing(\RingS)}} and
\mthbox{{\sing(\ringS)=\ringS \cap \sing(\RingS)}}.

Next, we will show that $\sing(\ringS)=\sing(\ring)\cdot \T(\varX)$ and derive a similar formula   for  $A$.  The proof is similar to the classical case of polynomial rings  
in an arbitrary set of commuting variables. 
For every finite subset $\varX' \subset \varX$, write $\ringS\sb{\varX'}=\ring[\varX' ;\bsym{\alpha}]$ for the skew polynomial  subextension  of $\ringS$ generated over $\ring$ by the variables from $\varX'$.
The classical result~\cite[Theorem~2.7]{Matczuk1995}  already grants that:
if $L$ is a uniform or an essential left ideal in $\ring$, then so is the left ideal
$\ringS\sb{\varX'}L$ in $\ringS\sb{\varX'}$.
Since $\ringS$ is clearly a union of the skew polynomial subrings 
$\ringS\sb{\varX'}$ for $\varX'$ running over the family of all finite subsets in $\varX$, 
it still holds that:
\startalignl 
if $L$ is a uniform or an essential left ideal in $\ring$, then so is the left ideal $\ringS L$ in $\ringS$. Thanks to Proposition~\ref{prop:udim-essringext},  the same holds for $\RingS=\ring[\varX\sp{\pm1};\bsym{\alpha}]$ and  $A$. 
\stopalignl 
It follows immediately that $\sing(\ring) \T(\varX)\subset \sing(\ringS)$ and
$\sing(\ring)\sT\subset \sing(A)$.  Proceeding by contradiction and assuming that 
$\sing(\ringS)$ is not contained in $\sing(\ring) \cdot \T(\varX)$, one can form  a polynomial $h\in \sing(\ringS) $   with $h\notin \sing(\ring) \cdot \T(\varX)$. 
By Remark~\ref{rem:term-fct-skwelaurentpolring}, the free commutative monoid 
$(\T(\varX),\cdot )$ enjoys a  linear  ordering  such that:  
all $0\neq f,g\in  \ringS$, we have      
$fg=\lc(f)\alpha[\lt(f)](\lc(g))\lt(f)\lt(g) +w$ for some polynomial $w$   with $\lt(w)<\lt(f)\lt(g)$.  We may therefore let  $\lambda=\lt(h)\in \T(\varX)$   and $a=\lc(h)$.
Since $a\notin \sing(\ring)$, choose $b\in \ring $ non-zero such that 
$\ring b\cap \ann[\ring](a)=0$. We claim that $\ringS b \cap \ann[\ringS](h)=0$, which will produce a contradiction to fact that $h$ lies in $\sing(\ringS)$.  
Let $f=a_1\tau[1]+\dotsm+a_p\tau[p]$ with $1\leq p\in \N$,
$ 0\neq a_1,\ldots, a_p\in \ring$ and $\tau[1]<\dotsm<\tau[p]$ in  $\T(\varX)$ such that 
$fb h=0$. Then we must have 
$a_p\alpha[{\tau[p]}](ba)=0$, that is,  
$\alpha[{\tau[p]},-1](a_p)ba=0$, forcing that $\alpha[{\tau[p]},-1](a_p)b=0$.
Thus  $a_p\tau[p]b=0$ and $ (f-a_p\tau[p])b=0$, hence inductively, on gets that
$a_i\tau[i]b=0$ for all  $i\in\irg[p]$ and $fb=0$.  This yields our claim that
$\ringS b \cap \ann[\ringS](h)=0$   and produces the needed  contradiction to the fact that $h$ lives in $\sing(\ringS)$. 
Hence, $\sing(\ringS) \subset \sing(\ring)\T(\varX)$ and one concludes that  $\sing(\ringS)=\sing(\ring)\T(\varX)$.  
But then, for all $f=a\sb{1}\tau[1]+\dotsm +a_p\tau[p] \in \sing(A)$ with $(a_i,\tau[i])\in \ring \times \sT$, there exists some term $\lambda\in \T(\varX)$  with
$\lambda\tau[i]\in \T(\varX)$ for all $i\in\irg[p]$. Hence having in view the fact that $\sing(A)=A\cap \sing(\RingS)$ while $\sing(\ringS)=\ringS \cap \sing(\RingS)$, we get that  
$\lambda f=\sum\sbsp[i=1,p]\alpha[\lambda](a\sb{i})\lambda\tau[i] \in \ringS \cap \sing(\RingS)=\sing(\ringS) =\sing(\ring) \T(\varX)$, so that for all $i\in\irg[p]$, $\alpha[\lambda](a_i)\in \sing(\ring)$ and $a_i\in\sing(\ring)$ because $\alpha[\lambda]$ is an automorphism of $\ring$.  Thus the equality $\sing(A)=\sing(\ring) \sT $ holds as well.

Finally with the additional assumption that $\ring$ has enough uniform left ideals
we will show that the same holds for $A$ and $\udim(A)=\udim(\ring)$.  
Since $\sT$ is a two-sided $\ring$-basis for $A$, every direct sum $\oplus\sb{i\in I}L\sb{i}$ of left ideals in $\ring$ extends to a direct sum $\oplus\sb{i\in I}\sT L\sb{i}$ of left ideals in $A$. Moreover it is already granted (as shown more above) that
if $L$ is a uniform or an essential left ideal in $\ring$, then so is the left ideal
$A L$ in $A$. Hence, Proposition~\ref{prop:invariance-udim.ring-ext} grants
that $\udim(A)=\udim(\ring)$ and 
$A$ has enough uniform left ideals, 
completing the proof of the theorem.
\stoppr


\startprop \label{prop:semiprime-skewmonoidring}
Let  $A=\ring[\sT;\bsym{\alpha}]$ be any subextension of a  skew Laurent polynomial ring  of injective type over $\ring$. Suppose that $\ring$ has the \ACC on left annihilators.   
Then if $\ring$ is (semi)prime, then so is   $A$.
\stopprop

\startpr  Thanks to Remark~\ref{rem:term-fct-skwelaurentpolring}, $\sT$ is endowed  with a linear-ordering and the leading term function
\Seq{\lt: \ring[\sT;\bsym{\alpha}] \to \sT} is a pseudo valuation 
on $\ring[\sT;\bsym{\alpha}]$. More precisely,     
for  all $0\neq f,g\in \ring[\sT;\bsym{\alpha}]$ with $\lambda=\lt(f)$ and $\tau=\lt(g)$, we have:     
$fg=\lc(f)\alpha[\lambda](\lc(g))\lambda\tau +w$ for some $w\in \ring[\sT;\bsym{\alpha}]$   with $\lt(w)<\lt(f)\lt(g)$.
With these  facts, the proof of the  proposition continues by repeating exactly the same arguments as in the univariate case by \cite[Proposition~3.6]{LeroyMatczuk1995}.
\stoppr

\startthm \label{thm:Goldie-subext-skewLaurentPolring} Let $A=\ring[\sT;\bsym{\alpha}]$  be any subextension in a skew Laurent polynomial ring of bijective type over  a (semi)prime ring $\ring$. Then, the coefficient ring $\ring$ is  left Goldie  precisely when so is  $A$, in which case,  
$A$  is also (semi)prime and 
$\udim(A)=\udim(\ring)$.
\stopthm 

\startpr  Recall that in a ring with finite uniform  dimension, every non-zero left ideal always contains a uniform left ideal.  By Goldie's theorem~\cite[2.3.6]{McConnellRobson1987}, a semiprime left Goldie ring is a semiprime ring with finite uniform dimension and  with a zero left singular ideal. Thus if $\ring$ is a (semi)prime left Goldie ring, then Theorem~\ref{thm:inv-udim-subext-skewLaurentPolring} together with Proposition~\ref{prop:semiprime-skewmonoidring} show that
the ring $A=\ring[\sT;\bsym{\alpha}]$  is also (semi)prime left Goldie
and $\udim (\ring) =\udim(A)$. 
Conversely, assume that  $A$  is  left Goldie. Then since already every direct sum $\oplus\sb{i\in I}L\sb{i}$ of left ideals in $\ring$ yields a direct sum
$\oplus\sb{i\in I}\sT L\sb{i}$ of left ideals in $A$, we get 
$\udim(\ring)\leq \udim(A)$. But as subring, $\ring$ already inherit from $A$ the \ACC on left annihilators. Hence
$\ring$ is left Goldie and, consequently, $\udim(\ring)=\udim(A)$ by virtue of 
Theorem~\ref{thm:inv-udim-subext-skewLaurentPolring}. 
\stoppr 

\section{\Esspecial subextensions in skew polynomial rings}
\label{sect:essspsubext-skewpolring}
For the sequel, $\RingS$ is an ambient skew polynomial ring 
of bijective type. 

\subsection{Skew polynomial rings and nice polynomials} 
\label{sect:skewpolring-nicepol} 

There is a need to review  skew polynomial rings in an arbitrary well-ordered set  
of variables.
Let \mthbox{\varX=\set{x\sb{s} \with s\in \ordn}} be indexed  by some non-necessarily finite  ordinal $\ordn$;  the well-ordering on $\varX$ being induced by that of the ordinal $\ordn$. 
The reader  is reminded that the class of all ordinals enjoys a well-ordering \tqt{\mth{\leq}} with:
\startcent  $\kappa <\eta$  if and only if $\kappa\in \eta$ for all ordinals $\kappa,\eta$. 
\stopcent
We have the following set  of \emph{standard terms},  as well as corresponding truncated subsets  for each ordinal $\kappa\leq \ordn$:
\begin{equation} \label{equ:stdTerms}
	\begin{aligned}
		& \T=\T(\varX) =\set{x\sb{s\sb{1}}\dotsm x\sb{s\sb{n}}  \with n\in\N,\,   0\leq s\sb{1}\leq \dotsm \leq s\sb{n} <\ordn}; \\
		&  \T[\kappa]=\T(\varX[\kappa]) \text{ with } \varX[\kappa]=\set{x\sb{s} \with s\in\kappa}\subset \varX;
		\text{ thus } \T[\ordn]=\T.
	\end{aligned}
\end{equation}
The skew polynomial ring $\RingS=\ring[\varX;\bsym{\alpha},\bsym{\delta}]$  and its  truncated subextensions $\RingS[\kappa]$ (for ordinals $\kappa \leq \ordn$)
are inductively   described as it follows.
\startitem
\item \mbox{$\RingS[\kappa]=\ring\T[\kappa]$ is the  free left $\ring$-module  with basis $\T[\kappa]$; $\RingS=\RingS[\ordn]$.}   As rings and  for a finite ordinal \mthbox{\kappa=\set{0,1,\ldots,k}},    \mthbox{\RingS[\kappa]=\ring[x\sb{0};\alpha[0],\delta[0],\ldots,x\sb{k};\alpha[k],\delta[k]]} 
is the usual iterated Ore extension of $\ring$, and \mthbox{\RingS[0]=\ring}.
\item  Let $\kappa<\ordn$ be an arbitrary ordinal and assume  we have defined the multiplication of the  truncated skew polynomial subextension $\RingS[\kappa]=\ring\T[\kappa]=\ring[\varX[\kappa]; \bsym{\alpha},\bsym{\delta}]$. Then $\RingS[\kappa+1]=\RingS[\kappa][x\sb{\kappa};\alpha[\kappa],\delta[\kappa]]$ is  the univariate skew polynomial ring extension of $\RingS[\kappa]$, where   the conjugation map \Seq{\alpha[\kappa]:\RingS[\kappa] \to\RingS[\kappa]} is an  endomorphism of $\RingS[\kappa]$ while
\Seq{\delta[\kappa] ; \RingS[\kappa] \to \RingS[\kappa]} is an $\alpha[\kappa]$-derivation.
\item We then notice for a limit ordinal $\ordn$ that  $\RingS$  and  any  subring $A\subset \RingS$   arise as increasing unions of  truncated  subrings:
\mthbox{\RingS=\cup\sb{\kappa\in\ordn} \RingS[\kappa]} and 
\mthbox{A=\cup\sb{\kappa\in\ordn} A\sb{\kappa}} with $A\sb{\kappa}=A\cap \RingS[\kappa]$ for $\kappa\in\ordn$.  
\stopitem  

\startnote
For every $\kappa\in\ordn$,  the associated  conjugation map and derivation map are suitably extended to maps \Seq{\alpha[\kappa], \delta[\kappa]: \RingS \to \RingS} defined over all of $\RingS$, with $\alpha[\kappa]$ the ring endomorphism of $\RingS$  and $\delta[\kappa]$ the $\alpha[\kappa]$-derivation such that:
\startcent
$\alpha[\kappa](\tau)=\tau$ and $\delta[\kappa](\tau)=0$ for 
every $\tau\in\T(x\sb{\eta} \with \kappa\leq \eta \in \ordn)$.
\stopcent 
For all  $n\in \N$ and every standard term $\lambda=x\sb{s\sb{1}}\dotsm x\sb{s\sb{n}}$
with $s\sb{1}\leq \dotsm \leq s\sb{n} <\ordn$, we shall equally use the alternative notation:
\startcent
$\alpha[\lambda]=\alpha[s\sb{1}]\circ\dotsm\circ \alpha[s\sb{n}]$.
\stopcent
\stopnote
Our main assumption  says that: 
\emph{the skew polynomial ring $\RingS$ is of bijective type}, that is, 
each conjugation map  $\alpha[\kappa] $ for  \mthbox{\kappa\in\ordn}  is an automorphism. In particular, each set   $\T(x\sb{\kappa})=\set{x\sbsp[\kappa,m] \with m\in\N}$ is  both a left and a right   $\RingS[\kappa]$-basis for $\RingS[\kappa+1]$. 



\startprop \label{prop:ext-dsum-left.subextskewpol} Let $A$  be any subring  in $\RingS$.  Then every direct sum $\oplus\sb{i\in I}L\sb{i}$  of left ideals in $\ring$ yields  a direct sum $\oplus\sb{i\in I}(A\cap AL\sb{i})$  
of left ideals in $A$. 
\stopprop 

\startpr We first prove the proposition when $A=\RingS$. Clearly proceeding by (transfinite) induction on the indexing ordinal $\ordn$, the base case is empty while the limit ordinal case is immediate (since in this case  $\RingS$ arises as increasing union of truncated skew polynomial rings). Now turn to the case that $\ordn=\kappa+1$ is a successor ordinal and let $\oplus\sb{i\in I}L\sb{i}$ be a direct sum of left ideals in $\ring$. The induction hypothesis grants that  the sum $\sum\sb{i\in I} \RingS[\kappa] L\sb{i}$ is direct while
$\RingS \RingS[\kappa] L\sb{i}=\RingS  L\sb{i}$ for all $i\in I$.  
But then, the set of univariate terms $\T(x\sb{\kappa})=\set{x\sbsp[\kappa,n] \with n\in\N}$ is also a right $\RingS[\kappa]$-basis of the univariate skew polynomial ring
$\RingS=\RingS[\kappa][x\sb{\kappa};\alpha[\kappa],\delta[\kappa]]$ and it follows that the sum  $\sum\sb{i\in I} \RingS L\sb{i}$ is still a direct sum  of left ideals in $\RingS$. Hence the proposition  is proved for $A=\RingS$. For a general subring $A\subset \RingS$,   since  each $A\cap AL\sb{i}$ is a left $A$-submodule in $\RingS L\sb{i}$ for all $i\in I$,  the sum $\sum\sb{i\in I} (A\cap AL\sb{i})$ is   direct as well. 
\stoppr

We  shall consider on standard terms the \emph{reverse lexicographic ordering} 
\tqt{\mth{\rlexeq}}  induced by the well-ordering on $\varX$  and  inductively defined as it follows. For all $\kappa\in\ordn$ and $x=x\sb{\kappa}\in\varX, \lambda,\tau\in\T[\kappa]$ and  $k,l\in \N$, we have:
\begin{equation*}
	\lambda x\sp{k} \rlex \tau x\sp{l} \iff k<l \text{ or }k=l \text{ and }
	\lambda \rlex \tau.  
\end{equation*}
We continue to denote this ordering simply by \tqt{\mth{\leq}}.  
As confirmed by the next technical lemma,  we get in this way 
a well-ordering on $\T$.  Hence there is an induced \emph{leading term function} \Seq{\lt: \RingS \to \T \cup \set{0}}:   for each non-zero   polynomial $f=a\sb1\tau[1]+\dotsm+a\sb m\tau[m]$  with $0\neq m\in\N$, $0\neq a\sb1,\ldots, a\sb{m}\in \ring$ and  terms $\tau[1]>\dotsm>\tau[m]$ in $\T$, we have $\lt(f)=\tau[1]$; the leading coefficient of $f$ is  
$\lc(f)=a\sb{1}$. For convenience, one sets $\lt(0)=0=\lc(0)$ and $0<\lambda$ for every $\lambda\in\T$.

\startlem \label{lem:stdTerm-well-ordering}
The set $\T$ of standard terms is well-ordered by the reverse lexicographic ordering induced by the well-ordering on $\varX$.
\stoplem

\def\Lv{\msf{Lvar}}
\startpr  Assume the non-trivial case, that is, $\varX$ is non-empty.  We reason by (transfinite) induction on well-ordered sets (or if one prefers, on ordinals).  The desired result is obvious if $\varX$ has only one variable, (that is,  if  $\ordn$ is the ordinal \tqt{\mth{1}}).
Suppose $\varX$ is not a singleton and the result holds for all proper well-ordered subset  in $\varX$.    Let $\Sigma$  be a non empty subset of $\T=\T(\varX)$, whose set of \tqt{leading variables}  is $\Lv(\Sigma)=\set{x\in\varX \with \Sigma \cap \T(\varX[\leq x])x \neq \vnoth}$, where  for  each ordinal $s \in \ordn$ we have  the truncated subsets
$\varX[<x\sb{s}]=\varX[s]=\set{z\in \varX \with z< x\sb{s}}$ and  
$\varX[\leq x\sb{s}]=\varX[s+1]=\set{z\in \varX \with z\leq x\sb{s}}$.   
By the well-ordering principle of $\varX$ and of $\N$,  we may let $y$ be the smallest variable in $\Lv(\Sigma)$ and $m$ be the smallest positive integer with  $\Sigma \cap \T(\varX[< y])y\sp{m} \neq \vnoth$.
Then, form the non-empty subset $\Sigma'=\set{\lambda\in \T(\varX[< y]) \with \lambda y\sp{m} \in \Sigma}$.  By the induction hypothesis,  $\Sigma'$ contains a smallest term, say, $\lambda$.
The inductive definition of the reverse lexicographic ordering on $\T$ readily shows that the term $\lambda y\sp{m}$ is the smallest term in $\Sigma$. 
\stoppr   

In the process of achieving the  main  gold of this section, 
the most difficult and important  step occurs when the indexing ordinal is a successor ordinal.  
So from now henceforth,  we  shall let  \mthbox{\extringS=\ringS[x;\alpha,\delta]} 
be a univariate skew polynomial ring  of bijective type over  an $\ring$-ring  $\ringS$, 
where $\ringS$ is  meant to be a specialized subextension of  a multivariate skew polynomial ring of bijective type. As  a basic assumption,   
\startalignl
$\ringS$ is required to be free  as left $\ring$-module with basis denoted as $\T[\ringS]$. 
\stopalignl 
Thus, the set \mthbox{\T[\ringS]\T(x)=\set{\lambda x\sp{n} \with \lambda\in\T[\ringS] \text{ and } n\in\N}} is a left $\ring$-basis  for $\extringS$.  We regard elements of $\T[\ringS]\T(x)$ as \tqt{\emph{terms}} and view $\extringS$ as a \tqt{polynomial-like} ring.

We now continue with  a review of some key strategies useful for investigating  the uniform dimension of polynomial-like rings.   
For our purpose,
it shall be necessary  to suitably review  the concept of \tqt{nice polynomials} and gather  some crucial technical facts about them.

\startdef \label{def:nice-annihilation} 
By an \emph{\mth{\ring}-nice element} of $\ringS$ is meant any non-zero $f\in \ringS$ whose coefficients in $\ring$ share a common left $\ring$-annihilator. Thus  $f=\sum\sbsp[i=1,m] a\sb{i}\lambda\sb{i}$ with $0\neq a\sb{i} \in \ring$,
$\lambda[i] \in \T[\ringS]$  and the $\lambda[i]$'s  pairwise distinct, and $\ann[\ring](a\sb{i})=\ann[\ring](a\sb{j})$ for all $i,j\in\irg[m]$.
We also say that the set $\coef(f)=\set{a\sb1,\ldots, a\sb{m}}$ is  \emph{\mth{\ring}-nice}. 
\stopdef  

We start with the following  first key lemma which for the purpose of the present framework includes and amplifies~\cite[
Lemma~2.1 and Proposition~2.2 ]{Shock1972}. 

\startlem \label{lem:eql-ann} 
Let $L$ be any essential left ideal of $R$ and $E\subset R$  any finite subset with $E\setminus\set{0}$ non-empty.
Then for a suitable $c\in R$,  $c E\setminus\set{0}$ is a    non-empty $\ring$-nice subset of $L$. If moreover  $E$ contains some element $a$  with $\ann[\ring](a)\subset \ann[\ring](E)$, 
then   $c\in R$ may be chosen such that  $\ann[\ring](c E)=\ann[\ring](c a)$.
\stoplem 


\startpr 
For  the  main part of the lemma   and starting with the case that $E$  contains some element $a$ with $\ann[\ring](a)\subset \ann[\ring](E)=\cap\sb{x\in E}\ann[\ring](x)$, we will show that there is some $c\in R$ such that $c E \neq \set{0}$, $c E\setminus\set{0}$ is $\ring$-nice   and $\ann[\ring](c E)=\ann[\ring](ca)$.
Proceeding  by induction on the number of  non-zero elements in $E$, 
since the base case is obvious, 
we may further suppose that there exists some $0\neq b\in E$  with    $\ann[\ring](b)\neq \ann[\ring](a)$. Thus pick some $c_1\in \ann[\ring](b)\setminus \ann[\ring](a)$; then   
$c_1 E$ has fewer non-zero elements than $E$,   $0\neq c_1a\in c_1E$   and
$\ann[\ring](c_1a)\subset \ann[\ring](c_1E)$. Hence by the way of induction, 
there is some $c_2\in R$ such that  $c_2\cdot (c_1a)\neq 0$  and $\ann[\ring](c_2c_1a)$ is the  left annihilator of every non-zero element  of $c_2c_1E$, and we are done.
Next in the case that  $E$  does not contain any element  whose left annihilator in $R$ is contained in the left annihilator of $E$ in $R$, one simply repeats the same arguments as before (but without any reference to $a$)  to  show  for a suitable $c\in R$
that $c E\setminus \set{0}$ is a non-empty   $R$-nice subset.

To complete the proof of the main part of the lemma, letting $E=\set{a_1,\ldots, a_n} \subset \ring\setminus\set{0}$ be $\ring$-nice,  we have to check that there is   some $c\in R$ with   $0\notin c E \subset L$ (since  $c E$  is still obviously  $R$-nice).   
Since by assumption  $L$ is an essential left ideal in $R$ while the $a\sb{i}$'s share a common left annihilator in $\ring$, inductively there are scalars $c_1,\ldots, c_n$  with
$0\neq c_1a_1\in L, 0\neq c_2c_1 a_2\in L,\,\ldots, 0\neq (c_n\dotsm c_1)a_n\in L$. Thus  letting 
$c=c_n\dotsm c_1$ yields the desired result. 
\stoppr


We include the following characterization of left annihilators for $\ring$-nice   multivariate polynomials. 
\startprop \label{prop:ann-skewpol-bijcase} Let $0\neq g\in \extringS=\ringS[x;\alpha,\delta]$ and  $p=\lc[\ringS](g)$;  suppose that  $p$ is $\ring$-nice and $\ann[\ring](p) =\ann[\ring](g)$, which is implied whenever $g$ is $\ring$-nice.   Then    
$\ann[\extringS](g) = \extringS \ann[\ring](g)$ provided   $\ann[\ringS](p)=\ringS \ann[\ring](p)$. 
In particular it holds  for every $\ring$-nice polynomial $g\in\RingS$  that $\ann[\RingS](g)=\RingS \ann[\ring](g)$.
\stopprop 

\startpr  The first claim of the proposition is proved by imitating the same arguments  for the univariate analogue~\cite[Lemma 1.1]{Matczuk1995}.   The proof of the second claim   proceeds by (transfinite) induction on the indexing ordinal $\ordn$, which essentially reduces to the case of a successor ordinal. For  clarity let us provide the needed evidence.  Let $g\in\RingS$ be $\ring$-nice.  The base of induction is trivial.
If $\ordn$ is a limit ordinal while the result holds for all ordinal $\kappa<\ordn$,   then $\RingS$  is an increasing union of truncated skew polynomial subrings $\RingS[\kappa]$ for $\kappa<\ordn$  and the induction hypothesis applies yields that $\ann[\RingS](g)=\cup\sb{\kappa<\ordn} \ann[\RingS[\kappa]](g)=\cup\sb{\kappa<\ordn} \RingS[\kappa]\ann[\ring](g)=\RingS\ann[\ring](g)$.
Now the  case that $\ordn=\kappa+1=\kappa\cup\set{\kappa}$ is a successor ordinal has been handled by the first claim, completing the proof of the proposition. 
\stoppr 

\subsection{\Esspecial subextensions and nice polynomials with prescribed degrees} 
\label{sect:esssp-subext-skewpolring-nicepol}

It's worth recalling from \cite{Matczuk1995} that  the natural requirement ensuring the invariance of the uniform dimension for a univariate skew polynomial ring $B=\ring[t;\sigma,\vtheta]$ of bijective type  amounts to the condition that  each $ B a$ contains an   $\ring$-nice polynomial of degree one for every non-zero scalar  $a\in\ring$. Let us examine this requirement for $\ring$-subextensions  $A\subset \extringS=\ringS[x,\alpha,\delta]$.  By assumption  $\alpha$ is bijective; but only the injectivity of $\alpha$ will be  enough here.  We call $ n\in \N$   an \emph{attainable  $x$-degree  in $A$} if there exists $f\in A$ with $\deg[x](f)=n$. We   set:
\begin{equation} \label{eq:attainable-xdegree}
	\begin{aligned}
		\Deg[x](A) &=\set{n\in\N \with \text{\mth{n} is an attainable \mth{x}-degree  in \mth{A}}}
		\text{ and, } \\
		\dd[A] &=\gcd(\Deg[x](A)\setminus\set{0})
		\text{ if } A\not\subset \ringS;  \text{ otherwise, } 
		\dd[A](x)=0.
	\end{aligned} 
\end{equation}
Clearly, $0\in \Deg[x](A) \subset \dd[A] \N$. 
If it happens that  $\dd[A]$ is not already an attainable $x$-degree in $A$, then 
there will be some $0\neq n,e\in \N$ with $n,n+e\in \Deg[x](A)$ but $e\notin \Deg[x](A)$,   and   given an  $\ring$-nice polynomial $f\in A$ with $n=\deg[x](f)$,  $Af$  cannot contain  any  polynomial with $x$-degree $n+e$. 
Such  a weird situation will definitely prevent one from adapting to $A$   classical strategies necessary for studying the  uniform dimension of skew polynomial rings. To resolve the last issue, we embed $A$ into another subextension $\eextA\subset \extringS$ where the jumps in attainable $x$-degrees disappear while $A$ and  $\eextA$ are still close enough to guarantee the invariance of the uniform dimension.

So we arrive at the following class  of specialized subextensions in skew polynomial rings, where we recall for any ring $A$  that
\begin{equation*}
	\regC(A)=\set{u\in A: \text{ \mth{u} is a regular element}}.
\end{equation*}

\startdef \label{def:special-subext}   
\startitem
\item[\pmsf(a)] 
By a \emph{\special subextension}  of the univariate skew polynomial ring  $\extringS=\ringS[x;\alpha,\delta]$ we mean   any $\ring$-subring 
$A$ of  $\extringS$ with the following special shape:
\startalignl[/vskip=0]
\mthbox{A = \sum\limits\sb{m\in\N}\hskip-0.4em (A\cap \ringS x\sp{m})}, \mthbox{\Deg[x](A)= \dd[A]\N}; to each  \mthbox{0\neq m \in \Deg[x](A)} are  associated some \mthbox{\splterm[m] =\upsilon[m]x\sp{m}\in  A}  with 
\mthbox{\upsilon[m]\in  \regC(\ringS)}, some \mth{\ring}-subring
\mthbox{\isubext(A)[m] \subset A\cap \ringS} and automorphism 
\mth{\alpha[(m)]} of   \mth{\isubext(A)[m]} (induced by \mth{\alpha})   such that: 
for all \mthbox{p\in \isubext(A)[m]} and
\mthbox{q\in \ringS} with \mthbox{qx\sp{m}\in A}, it holds that 
\mthbox{\deg[x]\bigl(\splterm[m]p -\alpha[(m)](p) \splterm[m]\bigr)<m}  and  \mthbox{uqx\sp{m}\in \isubext(A)[m]\splterm[m]}  for   some \mthbox{u\in \isubext(A)[m]\cap \regC(\ringS)}.
\stopalignl[/vskip=0] 
Notice  that for  $m=0$, the above requirements trivially hold  with
$\isubext(A)[0]=A\cap \ringS$,  $\splterm[0]=1=x\sp{0}$  and $\alpha[(0)]$ given by the identity map.  
%
\item[\pmsf(b)]
A  (strongly) nicely essential  subring of a \special  subextension of $\extringS$  is briefly referred  to as an \emph{\esspecial subextension} (or a \emph{strongly \esspecial subextension}) in $\extringS$.
\item[\pmsf(c)] 
An (or a \emph{strongly})  \emph{\esspecial subextension} of the multivariate skew polynomial ring $\RingS$  is inductively defined as it follows.
\startitem
\item If the indexing ordinal $\ordn$ is the zero ordinal, then $A$  is  a (strongly) nicely essential subring  of $\ring$. For a  limit ordinal $\ordn$,   $A=\cup\sb{\kappa \in\Lambda}A\sb{\kappa}$ arises as increasing union of rings with
$\Lambda\subset \ordn$ such that each $A\sb{\kappa}$ for $\kappa\in \Lambda$ is an (or a strongly) \esspecial subextension of $\RingS[\kappa]$.
\item If \mthbox{\ordn=\kappa+1=\kappa\cup\set{\kappa}} is a successor ordinal, then  $A$ is a (strongly) nicely essential subring  of  some prescribed \special subextension $\eextA$  of     \mthbox{\RingS=\RingS[\kappa][x\sb{\kappa};\alpha[\kappa],\delta[\kappa]]} regarded as  a univariate skew polynomial ring over $\RingS[\kappa]$,   while for all  $m\in \Deg[{x}\sb{\kappa}](\eextA)$ the associated $\ring$-subring $\isubext(\eextA)[m]$ arises again as an (or a strongly) \esspecial subextension of $\RingS[\kappa]$.  
\stopitem
\stopitem  
\stopdef 

\startrem \label{rem:spext-degree}
Keeping  notations as in part~\pmsf(a) of the above definition and letting $A$ be a  \special subextension in $\extringS$, for each polynomial $f\in A$ with $\deg[x](f)=m$, the special shape of $A$ ensures for some regular element $u\in\isubext(A)[m]$ that 
$uf=\splterm[m] p+f'$ with  $0\neq p\in \isubext(A)[m]$,  $f'\in A$ and $\deg[x](f')<m$.
\stoprem

\startexpl  The following observations are in order.
\startitem 
\item[\pmsf(i)] For a univariate skew polynomial ring  $\ring[x;\sigma,\vtheta]$ of bijective type,   
every $\ring$-subring $A$ that  is generated as  left  $\ring$-module 
by some subset  $\T(x\sp{d})\subset \T(x)$, with $d\in \N$, is evidently a \special subextension of $\ring[x;\sigma,\vtheta]$.
\item[\pmsf(ii)]  The ambient skew polynomial ring $\RingS$ is obviously a \special subextension of itself. If $A$ is any (strongly)  \esspecial subextension in $\RingS$, then so is every (strongly) nicely essential subring of $A$. 
\item[\pmsf(iii)] Easy but non-trivial examples
are obtained as it follows:  let $A\subset \RingS$, free as left $\ring$-module with basis $\sT\subset \T=\T(\varX)$, have  the following simple  structure: $A\sb{\kappa+1}=A\sb{\kappa}\cdot(\sT\cap\T(x\sb{\kappa}))$ for all $\kappa\in \ordn$.
Then $A$ is obviously a \special subextension  in $\RingS$, yet this ring needs not be a skew polynomial ring in the usual sense.
\stopitem

More involved examples of special subextensions  are the object of
section~\ref{sect:esssp-subext-2variate-skewpolring}. 
\stopexpl

We then include the following second key lemma, a more stronger and generalized version  of the ordinary  univariate analogue~\cite[Lemma~2.1]{Matczuk1995}.
\startlem \label{lem:nicepol} For the skew polynomial ring \mthbox{\extringS=\ringS[x;\alpha,\delta]} of bijective type,
assume that $\ann[\ringS](p)=\ringS\ann[\ring](p)$ for each $\ring$-nice  $p\in\ringS$.  Let   $A$  be any \special subextension of $\extringS$, not necessarily requiring that $\dd[A]$ be a member of $\Deg[x](A)$.
Then  statements~\pmsf(a)~and~\pmsf(b) below  are equivalent.
\startitem
\item[\pmsf(a)] 
$Aa$ contains an $\ring$-nice polynomial of degree $\dd[A]$ for  all 
\mthbox{0\neq a\in \ring}. In this  case, \mthbox{\Deg[x](A)= \dd[A]\N}.
\item[\pmsf(b)]
\mbox{For  each $\ring$-nice polynomial
	$f\in A$ and  all  $n\in\Deg[x](A)$ with    
	$n\geq \deg[x](f)$}, $ A f $ contains an $\ring$-nice polynomial 
$gf$ with \mthbox{\deg[x](gf)\hskip-0.15em=\hskip-0.15em n\hskip-0.15em=\hskip-0.15em\deg[x](g)\hskip-0.15em+\hskip-0.15em\deg[x](f)}.
\item[\pmsf(c)] Assume  that $\sing(\ring)\subset \sing(\ringS)$ and $\delta(\sing(\ringS))=0$.  Then  statement~\pmsf(a) holds  
whenever $\dd[A]$ is an attainable $x$-degree in $A$ and
$\ring \cap \sing(\isubext(A)[{{\dd[A]}}]) \subset \sing(\ring)$.  
\stopitem
\stoplem  

\startpr Consider the non-trivial  case that  the set $\Deg[x](A)$ of attainable $x$-degrees in $A$ does not reduce  to $\set{0}$;  so the natural number 
\mthbox{\dd=\dd[A]=\gcd(\Deg[x](A)\setminus\set{0})} is non-zero 
and \mthbox{\Deg[x](A) \subset \dd\N}.  For each $0\neq m\in\Deg[x](A)$,  the special shape of $A$ (Definition~\ref{def:special-subext}\pmsf(a))  grants a  special element
\mthbox{\splterm[m]=\upsilon[m]x\sp{m}\in A \cap \regC(\extringS)} 
with $\upsilon[m]\in \regC(\ringS)$, together with an induced automorphism $\alpha[(m)]$ of the associated $\ring$-subring $\isubext(A)[m]$ such that
\startcent 
$\deg[x](\splterm[m] p-\alpha[(m)](p)\splterm[m])<m$ for all 
$p\in \isubext(A)[m]$.
\stopcent 
We quickly see that $\Deg[x](A)$ is a submonoid of the additive monoid $(\N,+)$. Indeed for $0\neq m,n\in\Deg[x](A)$, we have:
$\splterm[m] \splterm[n]=\upsilon[m]\alpha[,m](\upsilon[n])x\sp{m+n} + \Delta$ with $\upsilon[m]\alpha[,m](\upsilon[n]) \neq 0$, $\Delta\in A$ and $\deg[x](\Delta)<m+n$, showing that  $m+n\in \Deg[x](A)$. 

Starting with the assumption that \pmsf(b) holds, let us prove \pmsf(a). If  $\dd$ is an attainable $x$-degree in $A$,  then \pmsf(a) becomes  a special case of   \pmsf(b).  So  we only need to show that $\dd$ is an attainable $x$-degree in $A$, which will also yield the equality $\Deg[x](A)= \dd\N$. Let $0\neq m\in \N$ be smallest such that $\dd m$ is an   attainable $x$-degree in $A$.  
Thanks to Lemma~\ref{lem:eql-ann}, we may choose an $\ring$-nice polynomial 
$w \in \ring \splterm[\dd m] \subset A$  with $\deg[x](w)=\dd m$. Our aim is to show that $m=1$. Assuming the contrary, by the definition of $\dd$  there must exist some   $\dd n \in \Deg[x](A)$ with $ n=ms+r$ for
some $1\leq s, r\in \N$ with $r<m$.  But by assumption, we can form an $\ring$-nice polynomial $g=px\sp{\dd ms}+g'\in Aw$  with $0\neq  p\in \ringS$, $g'\in A $ and $\deg[x](g')<\dd ms$;  another use of the assumption yields that we can choose a polynomial $h\in A$ with minimal degree, says $\dd k$, such that $\deg[x](hg)=\dd n=\dd ms+\dd r$.   Remark~\ref{rem:spext-degree} yields some  $u\in \isubext(A)[\dd k] \cap \regC(\ringS) \subset \regC(\extringS)$  
such that $uh= \splterm[\dd k]  q+h'$    with $0\neq q\in \isubext(A)[\dd k]$, $h'\in A$ and   $\deg[x](h')<\dd k$; in particular we still have $\deg[x](uh)=\deg[x](h)$ and $\deg[x](uhg)=\deg[x](hg)=\dd n$.  But then,
$\splterm[\dd k]=\upsilon[\dd k]x\sp{\dd k} \in A$ with 
$\upsilon[\dd k]\in \regC(\ringS)$     and 
$uhg=\upsilon[\dd k]\alpha[,\dd k](qp) x\sp{\dd\cdot (ms+k)} +\Delta$   with  $\Delta\in \extringS$ and \mthbox{\deg[x](\Delta)<\dd \cdot (ms+k)}.
Since $g$ is  $\ring$-nice, so is its leading polynomial-coefficient $p=\lc[\ringS](g)$ and $\ann[\ring](g)=\ann[\ring](p)$; 
the assumption about left annihilators in $\ringS$ yields that
$\ann[\ringS](p)=\ringS\ann[\ring](p) \subset \ann[\ringS](g)$. Hence $qp$ cannot be the zero polynomial because otherwise, $q g$ will be $0$ and $uhg=h'g$ while $\deg[x](h'g)=\deg[x](uhg)=\dd n$ and the polynomial $h'\in A$  satisfies $\deg[x](h')<\deg(h)$,  contradicting the minimality of $\deg[x](h)$.  
Since moreover $\upsilon[\dd k]$ is a regular element in $\ringS$, it follows 
that $\upsilon[\dd k]\alpha[,\dd k](qp)\neq 0$ as well, so that,
$\dd\cdot (ms+k) =\deg[x](uhg)= \dd n=\dd ms +\dd r$, showing that $k=r$. Thus the integer $\dd r=\deg[x](h)$ lies in $\Deg[x](A)$ with $1\leq r<m$, a contradiction to the minimality of $m$. Hence,  $m=1$ and
$\dd\in \Deg[x](A)$ as  expected. 

Next assuming  \pmsf(a), we will show that \pmsf(b) holds. Since the assumption already grants that  $\Deg[x](A)= \dd \N$, given  any $\ring$-nice polynomial $f\in A$ with $ \deg[x](f)=m\in\N$,   
it suffices to show that  $  A f $ contains an $\ring$-nice polynomial
$gf$ with $g\in A$,  $\deg[x](gf)=m+\dd$ and $\deg[x](g)=\dd$. Letting  $p=\lc[{\ringS}](f) \in \ringS$, the polynomial $f'=f-px\sp{m} \in \extringS$ satisfies    $\deg[x](f')<m$.  Since $f$ is $\ring$-nice, 
so is $p$ and, picking some  coefficient $a\in\coef[\ring](p)$, 
it holds that  
\begin{equation*}
	\ann[\ring](p)=\ann[\ring](a)=\ann[\ring](f) \text{ and } \ann[\extringS](f)=\extringS\ann[\ring](a)=\ann[\extringS](a),
\end{equation*}
where the two rightmost equalities are granted by   Proposition~\ref{prop:ann-skewpol-bijcase}.  
In particular, we obtain an isomorphism of left $\extringS$-modules
\Seq{\extringS a \to[/isymb] \extringS f} mapping each $ha$ to $hf$ for all $h\in\extringS$, and $\ann[\extringS](hf)=\ann[\extringS](ha)$ for all $h\in\extringS$. Now by \pmsf(a), the family $\mc{E}=\set{h\in A \with  \deg[x](ha)=\dd  \text{ and   \mth{ha}  is \mth{\ring}-nice}}$
is non-empty. 

Picking some $h\in\mc{E}$ with minimal $x$-degree $e$, we will show that \mthbox{e=\dd}. 
One may write
$ha=qx\sp{\dd} +w$  with $q\in\ringS$, $w\in A$ and $\deg[x](w)<\dd$; since $ha$ is $\ring$-nice, so is its leading polynomial-coefficient  $q=\lc[\ringS](ha)$ and $\ann[\ring](ha)=\ann[\ring](q)$, so that, 
$\ann[\extringS](ha)=\extringS \ann[\ring](q)=\ann[\extringS](q)$ (still by Proposition~\ref{prop:ann-skewpol-bijcase}).  
Remark~\ref{rem:spext-degree}  yields some   $u\in \isubext(A)[e]$, regular in $\ringS$ and  such that
$uh=\splterm[e] q'+h'$  with $0\neq q'\in \isubext(A)[e]$, 
$\deg[x](uh)=e=\deg[x](\splterm[e])$, $h'\in A$ and $\deg[x](h')<e$.
In particular, $\deg[x](uha)=\deg[x](ha)=\dd $ while $uha=\splterm[e] q'a+h'a$. We claim that $q'a\neq 0$. If this is not the case, then  since the polynomial $h'a=uha=uqx\sp{\dd}+uw$ has $x$-degree $\dd$ with $\deg(uw)<\dd$ while $\ann[\ring](uq)=\ann[\ring](uha)$, the first key Lemma~\ref{lem:eql-ann} shows for some suitable $c\in\ring$  that $ch'a$ is $\ring$-nice and still has $x$-degree $\dd$; but then, this will produce a contradiction to the minimality of $e=\deg(h)$  because the polynomial $ch'\in A$ has $x$-degree less than $e$.
Hence, $q'a\neq 0$ as claimed. Now, $\splterm[e]=\upsilon[e]x\sp{e}$ with $\upsilon[e]\in \ringS$ regular, thus  $uha=\splterm[e] q'a+h'a=\upsilon[e]\alpha[,e](q'a)x\sp{e}+h''$
with $0\neq \upsilon[e]\alpha[,e](q'a) \in\ringS$, 
$h''\in A$ and $\deg[x](h'')<e$.  So,
$e=\deg[x](uha)=\deg[x](ha)=\dd$. 

Then  we may let $h=vx\sp{\dd}+v'$  for some  $(v,v')\in \ringS\times A$  with  $\deg[x](v')<\dd$;  thus, 
$ha=v\alpha[,\dd](a) x\sp{\dd}+ h'$ with 
$0\neq v\alpha[,\dd](a)\in \ringS$, 
$h' \in A$   and  $\deg[x](h')< \dd$. 
Since   $ha$  is $\ring$-nice,  so is $v\alpha[,\dd](a)$  and $\ann[\ring](ha)=\ann[\ring](v\alpha[,\dd](a))$.
It follows that  $hf$ is a non-zero polynomial with $\ann[\ring](hf)=\ann[\ring](ha)=\ann[\ring](v\alpha[,\dd](a))$. 
But  we have:
\begin{equation*}
	hf= 
	v\alpha[,\dd](p) x\sp{m+\dd} +f''   \text{ with } f''\in \extringS  \text{ and }  \deg[x](f'') < m+\dd.
\end{equation*}  Next, we wish to prove that $\ann[\ring](v\alpha[,\dd](p)) \subset \ann[\ring](hf)$. Arbitrary let $b\in\ring $ with $b v\alpha[,\dd](p)=0$. 
Then, $\alpha[,-\dd](bv)p=0$  while for the coefficient  $a\in\coef[\ring](p)$ it is already granted that $\ann[\extringS](p)=\extringS\ann[\ring](a)$. 
It follows that  $\alpha[,-\dd](bv)a=0$ and
$b\in \ann[\ring](v\alpha[,\dd](a))=\ann[\ring](hf)$.  Hence  $\ann[\ring](v\alpha[,\dd](p)) \subset \ann[\ring](hf)$ as expected, and since moreover $hf\neq 0$ we must have $v\alpha[,\dd](p)\neq 0$, so that, $hf$ has $x$-degree $m+\dd$. 
Now,  invoking the first key Lemma~\ref{lem:eql-ann}  yields a suitable scalar $c\in\ring$ such that  $chf$ is $\ring$-nice while
$cv\alpha[,\dd](p) $ is still non-zero. Thus $\deg[x](chf)=m+\dd$ and $\deg[x](ch)=\dd$, completing  the proof  that  \pmsf(a) implies \pmsf(b).

Continuing with part~\pmsf(c) of the lemma, suppose that:  
$\sing(\ring)\subset \sing(\ringS)$,  $\delta(\sing(\ring))=0$,  $\dd$ is an 
attainable $x$-degree in $A$ 
while $\ring \cap \sing(\isubext(A)[{{\dd}}]) \subset \sing(\ring)$. 
We have to prove for every non-zero scalar $a\in \ring$  that $Aa$ contains an $\ring$-nice polynomial  with $x$-degree $\dd$.   Setting $B= \isubext(A)[\dd]$,  
we are granted with some special elements   $\upsilonup[\dd] \in\regC(\ringS)$ and  $\splterm=\splterm[\dd]=\upsilonup[\dd] x\sp{\dd} \in A \cap \regC(\RingS)$, together with an induced automorphism $\alpha[(\dd)]$ of the associated $\ring$-subring  $B\subset \ringS$ such that $\splterm a=\alpha[(\dd)](a) \splterm +\Delta$ with $0\neq \alpha[(\dd)](a)\in B$, $\Delta\in A$ and $\deg[x](\Delta)<\dd$.  First consider the case that $a \notin \sing(\ring)$. Then since by assumption $\ring \cap \sing(B) \subset \sing(\ring)$ it follows that $a\notin \sing(B)$ and 
$\alpha[(\dd)](a)\notin \sing(B)$.
Hence one may let $0\neq w\in B$ 
with $B w\cap \ann[B](\alpha[(\dd)](a)) =0$. In particular, since $\splterm$ is regular as element of $\extringS$, we have: 
\startcent 
$w\alpha[(\dd)](a)\splterm \neq 0$ and $\ann[\ring](w\alpha[(\dd)](a)\splterm) =\ann[\ring](w\alpha[(\dd)](a)) \subset \ann[\ring](w) \subset  \ann[\ring](w\splterm a)$.
\stopcent
Invoking Lemma~\ref{lem:eql-ann} yields some suitable  $c\in \ring$ such that $cw \splterm a$ is  $\ring$-nice with $cw\alpha[(\dd)](a) \splterm\neq 0$; hence
$\deg[x](cw \splterm a)=\dd$.  We then turn to the case that $a$ lies in 
$\sing(\ring)$. For every $n\in \N$, the Ore-rule grants that:
\startcent 
$x\sp{n}a=\alpha[,n](a)x\sp{n} + \delta[x\sp{n}](a)$  with
$\delta[x\sp{n}](a)=\sum\sbsp[k=1,n]x\sp{k-1}\delta(\alpha[,n-k](a))x\sp{n-k}$.
\stopcent 
But  by assumption, $\sing(\ring)\subset \sing(\ringS)$ while 
the derivation map $\delta$ vanishes on $\sing(\ringS)$ and   each  automorphism \Seq{\alpha[,k]: \ringS  \to\ringS} for   $k\in\N$  
sends $\sing(\ringS)$ to $\sing(\ringS)$. Consequently, $\delta[x\sp{n}](a)=0$ and  $x\sp{n} a=\alpha[,n](a)x\sp{n}$ for all $n\in\N$. 
In particular,  $\splterm a=\upsilonup[\dd] x\sp{\dd} a=\upsilonup[\dd]\alpha[,\dd](a) x\sp{\dd}$  with $0\neq\upsilonup[\dd]\alpha[,\dd](a) \in \ringS$ (because $\upsilonup[\dd]$ is a regular element of $\ringS$ and $\alpha[,\dd](a)\neq 0$). Thus we only need to invoke Lemma~\ref{lem:eql-ann} to choose  a suitable $c\in\ring$  such that  
the polynomial $c\splterm a=c \upsilonup[\dd]\alpha[,\dd](a)x\sp{\dd}$ is   $\ring$-nice  with degree $\dd$. This completes the proof of \pmsf(c), as well as the proof of the lemma. 
\stoppr

The next technical lemma is one of the most trickiest result,  specific to the general setting of subextensions in skew polynomial rings,  dealing  with left ideals in   \special subextensions that are induced from  the coefficient ring. This result complements the \tqt{nice polynomial} strategy provided by the above lemma, and all together they are crucial for implementing the  proof of the main  theorems more below. Continue to keep the notations introduced by Definition~\ref{def:special-subext}.

\startlem \label{lem:ext-leftIdeal.subextskewpol} 
Let  $A$ be a \special subextension of $\extringS=\ringS[x;\alpha,\delta]$, 
$L$ be any left ideal in $\ring$.  Then
it holds for every  $0\neq f\in AL$  with $m=\deg[x](f)$ that 
$f\in \sum\sbsp[k=0,m](A\cap \ringS x\sp{k})L$ and there exists some  $u\in \isubext(A)[m]\cap\regC(\ringS)$ such that   
$ \lc[\ringS](uf)x\sp{m} \in \alpha[(m)] (\isubext(A)[m] L)\splterm[m]$, with $u=1$ if $m=0$. 
\stoplem 

\startpr  Definition~\ref{def:special-subext} already includes the requirement that  \mthbox{A=\sum\sb{k\in \N} A\cap \ringS x\sp{k}}. 
Let  $0\neq f\in AL$, while $\deg[x](f)$  is not set. 
We consider the   set $\mc{E}$ of all pairs $(m,u)$ defined as it follows:
\startitem
\item $m\in\N$; if $f$ is a constant polynomial (with respect to variable $x$)  then   $u=1$; otherwise, $u$  is any  element in  $\isubext(A)[m]$ that is regular as element of $\ringS$;
\item  the polynomial $uf$ expresses as sum of elements of the form $px\sp{k}a$ with $px\sp{k}\in A\cap \ringS x\sp{k}$, $a\in L$  and 
$k\leq m$.
\stopitem 
Since $0\neq f\in AL$,  the set $\mc{E}$ is clearly non-empty. We  therefore choose a pair $(m,u)\in\mc{E}$ with $m$ minimal; in the course of the proof we will show that $m=\deg[x](f)$.  The minimality of $m$ already yields that there is a non-zero  polynomial  $h=\sum\sbsp[i=1,n]p\sb{i}x\sp{m}a\sb{i}$ for some positive $n\in\N$ 
with $p\sb{i}x\sp{m}a\sb{i}\neq 0$, $p\sb{i}\in \ringS, p\sb{i}x\sp{m}\in A, a_i \in  L$   for all $i$, and such that the polynomial $uf-h$ is either zero (if $m=0$) or is  a 
sum of elements of the form $px\sp{k}a$ with $px\sp{k}\in A\cap \ringS x\sp{k}$, $a\in L$  and 
$k < m$. Notice since $u$ is a regular element of $\ringS$ that $\deg[x](uf)=\deg[x](f)$. 

If  $m=0$, then $uf=h\in (A\cap \ringS) L $, in particular $f$ was already constant and by the definition of $\mc{E}$ we already have   $u=1$; thus $f$ lies in $(A\cap \ringS) L $  as expected.  So, we continue with  the case that $m\geq 1$. 
In view of  the special shape  of $A$ (Definition~\ref{def:special-subext}\pmsf(a)), the special element  \mthbox{\lambda=\splterm[m] \in A\cap \ringS x\sp{m}} is regular as element of $\extringS$   and we are granted with some element \mthbox{u'=u_nu\sb{2}\dotsm u_1} with each   $u\sb{i}$ lying in the associated $\ring$-subring $\isubext(A)[m]\subset A\cap \ringS$ and regular   as an element of $\ringS$,  such that 
\startcent
$u' p\sb{i}x\sp{m}= u\sb{i}' \lambda$ for some $u\sb{i}' \in \isubext(A)[m]$,
$1\leq i\leq n$.
\stopcent
Readily,  $(m,u'u) \in \mc{E}$ while the polynomial \mthbox{u'uf-u'h=u'\cdot (uf-h)} expresses  as 
sum of elements of the form $px\sp{k}a$ with $px\sp{k}\in A\cap \ringS x\sp{k}$, $a\in L$  and  $k\leq m-1$.   Still by the special shape of $A$,  $\isubext(A)[m]$ comes endowed with an induced automorphism $\alpha[(m)]$   such that:
\begin{equation*}
	\deg[x](\lambda p-\alpha[(m)](p)\lambda) \leq m-1
	\text{ for all } p\in \isubext(A)[m]. 
\end{equation*}   With this observation, 
it holds  that $u\sb{i}'\lambda =\lambda\alpha[(m),-1](u\sb{i}') +\Delta[i]$ for some  $\Delta[i]\in A$ with $\deg[x](\Delta[i])\leq m-1$ for $1\leq i\leq n$, so that,
\begin{align*}
	u'h &=\sum\sbsp[i=1,n]\lambda  \cdot 
	(\alpha[(m), -1](u\sb{i}') +\Delta[i] )a\sb{i} = \lambda 	\alpha[(m) ,-1 ](v) + \Delta 	\text{ with } \\ 
	v &= \sum\sbsp[i=1,n]  u\sb{i}' \alpha[(m)](a\sb{i})   
	\text{ and } \Delta=	\sum\sbsp[i=1,n]\Delta[i]a\sb{i} \in A L, 
	\text{ with } \Delta[i] \in A,\,  \deg[x](\Delta[i])\leq m-1. 
\end{align*}
Thus the minimal property of the pair $(m,u)\in \mc{E}$ shows that the polynomial $v$ just  defined above must be non-zero. 
But then,  we also get that 
$u'h=v\lambda +\Delta'$ for some $\Delta'\in A$  with $\deg[x](\Delta')\leq m-1$,   
and  the polynomial $u'uf\in AL$ now  satisfies: 
$\deg[x](u'uf)=\deg[x](v\lambda)=m$ and $\deg[x](f)=m$  because   the element $u'u \in  \isubext(A)[m]\subset A\cap \ringS$ is regular in $\ringS$.  It also follows in particular that  $u'uf$  lives in $ \sum\sbsp[k=0,m](A \cap\ringS x\sp{k})L$ as expected.  Moreover,    
$\lc[{\ringS} ](u'uf)x\sp{m}=v\lambda=\sum\sbsp[i=1,n]  u\sb{i}' \alpha[(m)](a\sb{i})\lambda$ lives in $\isubext(A)[m]\alpha[(m)](L) \lambda=\alpha[(m)](\isubext(A)[m]L) \lambda$ as claimed by the last part of the lemma.  
\stoppr

\subsection{The uniform dimension of \esspecial subextensions}
\label{sect:udim-esssp-subext-skewpolring}
It helps recalling as set before that   $\extringS=\ringS[x;\alpha,\delta]$ is a univariate skew polynomial ring of bijective type over  an $\ring$-ring $\ringS$ that is free as left $\ring$-module with left $\ring$-basis denoted as $\T[\ringS]$.  For our main theorems, a specific choice for $\ringS$  will be    \anesspecial $\ring$-subextension in a multivariate skew polynomial ring of bijective type.  So, from now henceforth we  require the additional condition:
\begin{equation*}
	\begin{aligned}
		& \text{every \mth{\ring}-nice element \mth{p\in \ringS} has an \mth{\ring}-nice annhilator in \mth{\ringS}:  \mth{\ann[\ringS](p)=\ringS \ann[\ring](p)}}; \\
		& \text{and for all ideals
			\mth{L,L'} in \mth{\ring}, if \mth{L\cap L'=0} then  
			\mth{\ringS L\cap \ringS L'=0},} 
	\end{aligned}
\end{equation*} 
which, in view of Proposition~\ref{prop:ann-skewpol-bijcase}~and~\ref{prop:ext-dsum-left.subextskewpol}, holds for $\RingS$, as  well as for all \esspecial subextensions in $\RingS$  as we shall prove more latter below. 

We are now in the position to state the first main theorem of this section; 
where for  \anesspecial subextension $A$  
of $\extringS$, we write as usual  $\eextA$ for a prescribed \special subextension of $\extringS$ containing $A$ as nicely essential subring.

%

\startthm \label{thm:invaraince-uniform-univariate-case} 
Let   $A$ be  \anesspecial subextension of $\extringS$ such that $\eextA a$ contains an $\ring$-nice polynomial of $x$-degree $\dd[\eextA]$ for  all $0\neq a\in \ring$.       
Then, every uniform  or essential  left ideal $L$  of $\ring$ lifts to a 
uniform (or resp., an essential) left ideal $A\cap AL$ in $A$ provided the same holds for $\isubext(\eextA)[m]$ for all $m\in\Deg[x](\eextA)$.  When this is the case, 
if $\ring$ has enough uniform left ideals, then  so does $A$ and $\udim(A) = \udim(\ring)$.
\stopthm


\startpr  By virtue of Proposition~\ref{prop:udim-essringext} it holds that: $\udim(A)=\udim(\eextA)$, and for every left ideal $L$ in $\ring$, the left ideal $A\cap AL$ is uniform (or essential) in $A$ if and only  if so is the left ideal $\eextA L$ in $\eextA$. Hence we only need to prove the theorem when $A=\eextA$ is a \special subextension in $\extringS$.
Assuming the latter,  for the rest of the proof we let $m \in \Deg[x](A)$.
The now assumed special shape of $A$ grants   
some special  element  $\lambda=\splterm[m] \in A\cap \ringS x\sp{m}$ regular in $\extringS$,  an induced  automorphism $\alpha[(m)]$ of  the associated $\ring$-subring $\isubext(A)[m] \subset A\cap \ringS$, such that:
\begin{equation} \label{eq:pr-lem:extuniform-ess-ideal} 
	\tag{\mth{\star}}
	\begin{aligned}
		& \text{for all } p\in \isubext(A)[m], \deg[x](\lambda p-\alpha[(m)](p)\lambda) < m; 
		\text{and  for all } q\in\ringS \text{ with }   \\
		& qx\sp{m} \in A, \text{it holds  for some }
		u\in \isubext(A)[m]   \text{  regular in \mth{\ringS} that } 
		uqx\sp{m} \in \isubext(A)[m]\lambda. 
	\end{aligned}
\end{equation}
Of course (in Definition~\ref{def:special-subext}), we  set: $\lambda[0]=1=x\sp{0}, \isubext(A)[0]=A\cap \ringS$   and $\alpha[(0)] =\id$.   Also remember that the special shape of $A$ already includes the requirement that 
$A=\sum\sb{k\in \N} A\cap \ringS x\sp{k}$.  

\sectitem{Let $L$ be any uniform left ideal of $\ring$ that lifts to a uniform  left ideal  in $\isubext(A)[m]$  for all $m\in\Deg[x](A)$. Proving that $AL$ is still uniform as left ideal in $A$} 
Given non-zero polynomials $f,g\in A L$,  set  
$k=\deg[x](f)$, $m=\deg[x](g)$. We may further suppose that $k\leq m$ and,  by multiplying $f$ and $g$ by suitable scalars from the ground ring  $\ring$ we may assume thanks to the first key Lemma~\ref{lem:eql-ann} that $f$ and $g$ are $\ring$-nice.  
We will proceed  by ordinary induction on the natural number $k+m$ to show that $Af \cap Ag\neq 0$.

If $k+m=0$  then $f,g\in AL\cap \ringS$ while the technical Lemma~\ref{lem:ext-leftIdeal.subextskewpol} shows that $AL\cap \ringS=(A\cap \ringS) L=\isubext(A)[0]L$ while by assumption it holds that \mthbox{(A\cap \ringS) f\cap (A\cap \ringS) g \neq 0} and we are done. So turning to the induction step, suppose that  $k+m \geq 1$ and that for all non-zero polynomials 
$f',g' \in A L$  with $\deg[x](f')+\deg[x](g')<k+m$, it holds that
$A f'\cap A g'\neq 0$.  Lemma~\ref{lem:nicepol}  grants that $A f$  contains an $\ring$-nice polynomial $\wt{f}$ with $\deg[x](\wt{f})=m$. Then considering the leading polynomial-coefficients $p=\lc[{\ringS}](\wt{f})$  and $q=\lc[{\ringS}](g)$, we get $\ring$-nice polynomial 
in $\ringS$  such that  
$\wt{f} - px\sp{m}, g-qx\sp{m} \in A$ are polynomials  with $x$-degree less than $m$ while $\ann[\ring](p)=\ann[\ring](\wt{f})$ and $\ann[\ring](q)=\ann[\ring](g)$.  In particular, Proposition~\ref{prop:ann-skewpol-bijcase} ensures that
$\ann[\extringS](\wt{f})=\extringS\ann[\ring](p)=\ann[\extringS](p)$ 
and $\ann[\extringS](g)=\extringS\ann[\ring](q)=\ann[\extringS](q)$.
Next having in view~\eqref{eq:pr-lem:extuniform-ess-ideal} above,   the  technical Lemma~\ref{lem:ext-leftIdeal.subextskewpol} shows that  
$upx\sp{m}=p'\lambda$  and   $u'qx\sp{m}=q'\lambda$ for some    elements $u,u'\in \isubext(A)[m] \cap \regC(\ringS)$ and   non-zero polynomials $p',q'\in \alpha[(m)](\isubext(A)[m]L)$.
In particular (since $x\sp{m}, \lambda$ are regular in $\extringS$)   it holds that:    
\begin{align*} 
	\ann[\extringS](up)\hskip-0.2em=\hskip-0.2em\ann[\extringS](p')  \text{ and } \ann[\extringS](u'q)\hskip-0.2em=\hskip-0.2em\ann[\extringS](q');\; 
	\ann[\extringS](u\wt{f}) \hskip-0.2em=\hskip-0.2em \ann[\extringS](p')
	\text{ and }   \ann[\extringS](u'g) \hskip-0.2em=\hskip-0.2em \ann[\extringS](q').
\end{align*}  
But since  by  assumption   $\isubext(A)[m]L$ is a uniform left ideal in $\isubext(A)[m]$,  there are elements 
$p'',q''\in \isubext(A)[m]$   with 
\begin{equation*}
	0\neq p''p' =q''q',  \text{ so that } p''u\wt{f} \neq 0 \text{ and } q''u'g\neq 0 \text{ as well}.
\end{equation*}
The polynomial $h=q''u'g-p''u\wt{f}$ lies in $AL$  with
$\deg[x](h)<m$. If $h=0$ then we are done since in  this case,
$0\neq q''u'g=p''u\wt{f} \in A f \cap Ag$. Thus suppose that $h\neq 0$. Then the  induction hypothesis  yields some $f',g'\in A$  with
\begin{equation*}
	0\neq f'f=g'h=g'q''u'g-g'p''u\wt{f} \text{ and }
	g'q''u'g =f'f+g'p''u \wt{f}.
\end{equation*} 
Since $\ann[\extringS](u\wt{f})= \ann[\extringS](p')$  and  
$\ann[\extringS](u'g)= \ann[\extringS](q')$  while $p''p'=q''q'$, 
it follows that 
$\ann[\extringS](g'q''u'g)= \ann[\extringS](g'q''q')=
\ann[\extringS](g'p''p')=\ann[\extringS](g'p''u\wt{f})$.  Thus since also  $g'q''u'g-g'p''u \wt{f}\neq 0$, it follows that $g'q''u'g\neq 0$,
so that the polynomial $ g'q''u'g=f'f+g'p''u \wt{f}$ is a non-zero element in $Af\cap Ag$, completing the proof that $AL$ is a uniform left ideal in $A$.

\sectitem{Let $K$ be an essential left ideal of $\ring$ that lifts to an essential  left ideal  in $\isubext(A)[m]$  for all $m\in\Deg[x](A)$. Proving   that $AK$ is essential as left ideal in $A$}
Letting $0\neq f\in A$ and $m=\deg[x](f)$, we will show by   induction on  $m$ that  $ AK \cap Af\neq 0$. 
In view  of the first key Lemma~\ref{lem:eql-ann}, left multiplying $f$ by a suitable scalar  we may assume that 
$f$ is $\ring$-nice.  If  $m=0$, then $f\in A\cap \ringS $ and the assumption  already  grants  that
$(A\cap \ringS) K\cap (A\cap \ringS) f\neq 0$, so that,   
$A K\cap A f \neq 0$ as well.
Next, suppose that  $m\geq 1$ and the desired result holds in $x$-degree less than $m$.  
In view of \eqref{eq:pr-lem:extuniform-ess-ideal}, recall that  
the  special element $\lambda \in A \cap\regC(\extringS)$, 
the associated $\ring$-subring $\isubext(A)[m]\subset A\cap \ringS$ 
and the  induced automorphism  $\alpha[(m)]$ of
$\isubext(A)[m]$  satisfy:  
for all $p'\in \isubext(A)[m]$,  the polynomial $\lambda p'-\alpha[(m)](p')\lambda \in A$ has $x$-degree less than $m$; and   for some   
$u\in \isubext(A)[m]\cap\regC(\ringS)$ and  some non-zero polynomial 
$p\in \isubext(A)[m]$ it holds  that the polynomial $f'=uf-p\lambda \in A$ has $x$-degree less than $m$.   
The assumption on $K$ ensures that $\isubext(A)[m] K$ is an essential left ideal in $\isubext(A)[m]$;  
thus the left ideal $\isubext(A)[m] \alpha[(m)](K)=\alpha[(m)](\isubext(A)[m] K)$ is essential in $\isubext(A)[m]$  and it holds for some $v\in \isubext(A)[m]$ that  $0\neq vp \in \alpha[(m)](\isubext(A)[m] K)$.  Since by definition
$\lambda$ is a regular element of $\extringS$ it follows that 
$ vp\lambda\neq 0$ and the polynomial $vuf=vp\lambda +vf'$  still has  $x$-degree   given by $m$. 
Notice having in view Lemma~\ref{lem:eql-ann} that,  if $\ann[\ring](vp\lambda)\nsubset\ann[\ring](vuf)$  then left multiplying  $vuf$  by a suitable scalar $c$ we may form an $\ring$-nice polynomial $cvuf\in A$ with $\deg[x](cvuf)<m$  and then  apply the   induction hypothesis to produce a non-zero element in $AK\cap Af$.  
So we may continue our discussion by assuming that 
$\ann[\ring](vp\lambda)\subset \ann[\ring](vuf)$. Then still  by Lemma~\ref{lem:eql-ann}, we may further assume that both $vp\lambda$ and $vuf$ are  $\ring$-nice polynomials with $\ann[\ring](vuf)=\ann[\ring] (vp\lambda)$. While the non-zero polynomial $\alpha[(m),-1](vp)\in \isubext(A)[m]$ needs not be $\ring$-nice, by  left multiplying  the latter by a suitable scalar from $\ring$ to form an $\ring$-nice element in $\isubext(A)[m]\subset\ringS$  and then   invoking   the second key Lemma~\ref{lem:nicepol}, we can produce some $\ring$-nice polynomial \mthbox{g=(qx\sp{m}+g') \alpha[(m),-1](vp)}  with $q\in \ringS, qx\sp{m}, g'\in A$ such that $\deg[x](g')<m$ and  $\deg[x](g)=m$. But  we already  know that
$wq x\sp{m}=q'\lambda$ for  some $w\in \isubext(A)[m]\cap \regC(\ringS)$ and some non-zero polynomial $q'\in \isubext(A)[m]$.  Thus,   
\startalignl
$wg=(q'\lambda +wg')\alpha[(m),-1](vp)=q'vp\lambda +\Delta$   with $\Delta\in A$ and $\deg[x](\Delta)<m$; \\
and letting $ p'=\lc[{{\ringS}}] (g)$ it holds  that  $q'vp\lambda=\lc[{{\ringS}}] (wg)x\sp{m}=wp'x\sp{m}$.
\stopalignl 
Since $g$ is $\ring$-nice,  so is its leading  polynomial-coefficient 
$ p'=\lc[{{\ringS}}] (g)$  and  
$\ann[\ring](g)=\ann[\ring](p')=\ann[\ring](p'x\sp{m})$.
Now  since $vuf=vp\lambda +vf'$  with $\deg[x](vf')<m$, the polynomial $ q'vuf-wg \in A$  has $x$-degree less than $m$.
If  $q'vuf-wg=0$, then $0\neq  q'vuf= wg\in AK$  and we are done.
Next  assuming that $q'vuf-wg$ is non-zero,  the  induction hypothesis grants some $w'\in A $  with
\startcent 
$0\neq w' \cdot (q'vuf-wg)\in A K$.
\stopcent  
But then,  the polynomials $vuf,vp\lambda,g$ and $p'x\sp{m}=\lc[{{\ringS}}](g)x\sp{m}$ are $\ring$-nice with
\startcent 
$\ann[\ring](vuf)= 
\ann[\ring](vp\lambda)$,  $\ann[\ring](g) = 
\ann[\ring](p'x\sp{m})$  and   $q'vp\lambda =wp'x\sp{m}$. 
\stopcent 
Thus   Proposition~\ref{prop:ann-skewpol-bijcase} shows  that
$\ann[\extringS](vuf)=\ann[\extringS](vp\lambda)$ while 
$\ann[\extringS](g)=\ann[\extringS](p'x\sp{m})$, and it follows that:\,  
$\ann[\extringS](q'vuf)=\ann[\extringS](q'vp\lambda)=
\ann[\extringS](wp'x\sp{m})=  \ann[\extringS](wg)$. 
%
Whence, the fact  that $w' \cdot (q'vuf-wg)$ is non-zero forces $w'q'vuf$ to be non-zero and consequently, $0\neq w'q'vuf =w'\cdot (q'vuf-w g)+w'wg \in AK$.  This completes the proof that $AK$ is an essential left ideal in $A$.  

For the last claim of the theorem, the additional assumption yields that 
every uniform or essential left ideal in $\ring$ lifts to a uniform (or resp., an essential) left ideal in $A$.  Next given any left ideals $L,L'$ in $\ring$  with $L\cap L' =0$, let us quickly verify that $AL\cap AL'=0$. By our setting for the $\ring$-ring $\ringS$ we know  that  $\ringS L\cap \ringS L'=0$. Hence for the univariate skew polynomial ring of bijective $\extringS$, Proposition~\ref{prop:ext-dsum-left.subextskewpol} grants that $\extringS L\cap \extringS L'=0$,  thus as  $\ring$-subring in $\extringS$, it holds as well that $A L\cap A L'=0$.
Now,   Proposition~\ref{prop:invariance-udim.ring-ext} applies showing that $\udim(A)=\udim(\ring)$,  completing the proof of the theorem.
\stoppr 

\startcor \label{cor:invaraince-uniform-univariate-case} 
Let   $A$ be \anesspecial subextension  in a univariate skew polynomial ring of bijective type  $\ring[x;\alpha,\delta]$, such that $\eextA a$ contains an $\ring$-nice polynomial of $x$-degree $\dd[\eextA]$ for  all $0\neq a\in \ring$.       
Then, every uniform  or essential  left ideal $L$  of $\ring$ lifts to a 
uniform (or resp., an essential) left ideal $A\cap AL$ in $A$; 
$\sing(A)=A\cap \sing(\ring) \T(x)$;       
if $\ring$ has enough uniform left ideals, then  so does $A$ and $\udim(A) = \udim(\ring)$.
\stopcor

\startpr Except for the claim about the left singular ideal, the other claims of the corollary are granted by Theorem~\ref{thm:invaraince-uniform-univariate-case}
specialized to the case that $\ringS$ coincides with the ground ring $\ring$.
It follows in particular   that $\sing(\ring) \subset \sing(AR)$ with $AR$ considered as  left $A$-module. Now the claim that
$\sing(A)=A\cap \sing(\ring) \T(x)$ is the object of a more general result given by Proposition~\ref{prop:singual-ideal-subext-skewpolyring} below.
\stoppr

Recall by Lemma~\ref{lem:stdTerm-well-ordering} that the reverse lexicographic ordering yields a well-ordering on  $\T=\T(\varX)$; correspondingly for each $f\in \RingS$,   
we let $\lt(f) \in \T$ and  $\lc(f) \in \ring$  be respectively  the leading  term and the leading coefficient of $f$.


\startprop  \label{prop:singual-ideal-subext-skewpolyring}
Let $A\subset \RingS$ be any  subring   such that $A\cap AL\neq 0$ for all non-zero left ideal $L$ in $\ring$. Then for all  $ f\in \RingS \setminus \sing(\ring) \T$, there exists a non-zero left ideal $L\subset \ring$  with  $\RingS L\cap \ann[\RingS](f)=0$. Consequently, $\sing(A)\subset \sing(\ring) \T$,   and $\sing(A) = A\cap \sing(\ring)\T$ provided $\sing(\ring)\subset \sing(A\ring)$,  with $A \ring$ viewed as left $A$-submodule in $\RingS$. 
\stopprop

\startpr  
Let $0\neq f\in \RingS$, $a=\lc(f)$.
We start by assuming that $a\notin \sing(\ring)$. Then it holds for some non-zero left ideal  $L$ of $\ring$  that $L\cap \ann[\ring](a)=0$. 
We shall reason  by transfinite induction on the indexing ordinal $\ordn$  to show that $\RingS L\cap \ann[\RingS](f)=0$.  Since the base case is obvious, we turn to the induction step and suppose that the desired result holds for all ordinal less than $\ordn$.  Now $\ordn$ is either a limit ordinal or a successor ordinal. The case that $\ordn$ is a limit ordinal is also straightforward. Indeed for this case, $\RingS=\cup\sb{\kappa<\ordn} \RingS[\kappa]$ is an increasing union of truncated skew polynomial rings of bijective type over $\ring$ while for some fixed $\kappa[0]$, $f \in \RingS[\kappa]$ for $\kappa[0] \leq \kappa <\ordn$.  But, 
$\ann[\RingS](f)=\cup\sb{\kappa_0\leq \kappa<\ordn} \ann[{\RingS[\kappa]}](f)$ and
$\RingS L=\cup\sb{\kappa_0\leq \kappa<\ordn} \RingS[\kappa] L$, while the induction hypothesis grants that $\RingS[\kappa]L \cap \ann[{\RingS[\kappa]}](f)=0$ for $\kappa_0\leq \kappa <\ordn$, which shows that $\RingS L \cap \ann[\RingS](f)=0$.
Next consider the case  of a successor ordinal 
\mthbox{\ordn=\kappa+1=\kappa\cup\set{\kappa}}. Then  setting $x=x\sb{\kappa},\alpha=\alpha[\kappa]$ and $\delta=\delta[\kappa]$, we get that  $\RingS=\RingS[\kappa][x;\alpha,\delta]$ is a univariate skew polynomial ring of bijective type over the truncated skew polynomial ring $\RingS[\kappa]$ of bijective type over $\ring$. Letting $m=\deg[x](f)$, one writes: $f=px\sp{m} +f'$  with $p\in\RingS[\kappa]$, $f'\in \RingS$  and $\deg[x](f')<m$. By the definition of the reverse lexicographic ordering, $\lc(p)=\lc(f)=a$. 
So the induction hypothesis grants that $\RingS[\kappa]L\cap \ann[{\RingS[\kappa]}](p)=0$. Since the conjugation map  $\alpha$ is bijective over $\RingS[\kappa]$, the set of univariate terms $\T(x)=\set{x\sp{n} \with   n\in \N}$ also serves as right $\ring$-basis for 
$\RingS$. So every element $g\in \RingS L$ expresses  as 
$g=x\sp{n_1} w_1+\dotsm +x\sp{n_k}w_k$  with
$k\in\N$,  $n_1>\dotsm >n_k$ and $w_1,\ldots,w_k\in \RingS[\kappa]L$. Suppose that $gf=0$. Then
since $gf=\alpha[,n_1](w_1p)x\sp{n_1+m}+g'$ for some $g'\in\RingS$  with $\deg[x](g')<n_1+m$, one must have $w_1p=0$, so that, $w_1=0$ by the transfinite induction hypothesis. Thus inductively, one finds that $w_1,w_2,\ldots, w_k$ all vanish and $g=0$.   Hence we have proved that $\RingS L\cap \ann[\RingS](f)=0$ under the assumption that $\lc(f) \notin \sing(\ring)$.    

Now for a general polynomial  $f\in \RingS$ that does not live in $\sing(\ring)\T$, we may write $f=g+f'$  with $g\in \sing(\ring) \T$  and
$\lc(f')\notin \sing(\ring)$. So let $L'$ be a non-zero left ideal of $\ring$  with
$L'\cap \ann[\ring](\lc(f'))=0$. From above,  it holds that $\RingS L' \cap \ann[\RingS](f')=0$. Since in any ring the intersection of finitely many essential left ideals  must be an essential left ideal, we get that $\ann[\ring](g)$ is an essential left ideal of $\ring$. Hence,  setting  $L=\ann[\ring](g) \cap L'$ yields a non-zero left ideal in $\ring$  with $\RingS L \cap \ann[\RingS] (f)=0$ as expected. 

Turning to the second part of the proposition, let  $A \subset \RingS$ be  any subring  with $A\cap AL\neq 0$ for every non-zero left ideal $L$ in $\ring$. Given  any polynomial  $f\in A$   with $f \notin \sing(\ring) \T$,  the above paragraph grants  some non-zero left ideal $L$ of $\ring$ with   $\RingS L\cap \ann[\RingS](f)=0$, so that,  $(A\cap A L)\cap \ann[A](f)=0$ while $A\cap AL$ is a non-zero left ideal in $A$; hence, $f\notin \sing(A)$, which shows   that  $\sing(A) \subset \sing(\ring)\T$. Finally given the additional assumption that $\sing(\ring)$ is contained in the singular ideal of the left $A$-submodule $A \ring $ of $\RingS$, we have to prove that 
$ A\cap \sing(\ring)\T \subset \sing(A)$. Let 
$g=a\sb1\tau[1]+\dotsm +a\sb{m} \tau[m]\in A$ with $0\neq m\in\N$ and $(a\sb{i},\tau[i]) \in \sing(\ring) \times \T$ for $1\leq i\leq m$. 
Since  $a\sb{i}$ lies in   $\sing(\ring) \subset \sing(A \ring )$,  each $\ann[A](a\sb{i})$ is an essential left ideal in $A$ and so, 
$\cap\sbsp[i=1,m] \ann[A](a\sb{i})$ is still an essential left ideal in $A$ that is contained in $\ann[A](g)$. Hence, $g\in \sing(A)$. This completes the proof of the proposition. 
%
\stoppr

We  now continue with a multivariate version of Theorem~\ref{thm:invaraince-uniform-univariate-case}.

\startthm \label{thm:inv-udim-subext-skewpolring} Suppose that each derivation map associated with $\RingS$ vanishes on $\sing(\ring)\T$. Let $A$ be  \anesspecial  subextension in $\RingS$. Then, every uniform or essential left ideal $L$ in  $\ring$ lifts to a uniform (or resp., an essential)  left ideal $A\cap AL$ in $A$; $\sing(A)=A\cap \sing(\ring) \T$. 
And if  $\ring$ contains enough  uniform left ideals,   then so does  $A$ and $\udim(A) = \udim(\ring)$.  
\stopthm 
It's worth pointing our that, even when one specializes to the classical setting of the ambient skew polynomial ring $\RingS$ (as special subextension of itself), 
the above theorem is still more stronger than its classical analogue \cite[theorem~2.4]{Matczuk1995} in the following aspects:
\startitem
\item $\udim(\ring)$ needs not be finite and the set $\varX$  of variables needs not be finite;
\item the  assumption  that $\delta[\kappa](\sing(\ring)\T)=0$ for all $\kappa\in\ordn$ is implied  by the classical more stronger condition that $\sing(\ring)=0$; moreover, we obtain that  $\sing(\RingS) =\sing(\ring) \T$, generalizing the classical fact that $\sing(\RingS)=0$  if  $\sing(\ring)=0$ while  $\varX$ is finite.
\stopitem   

\startpr  Using induction on $\ordn$, we start by  proving the following first claim.

\sectitem{Every non-zero left ideal $L$ of $\ring$ lifts to a
	nicely essential  submodule $A\cap AL$ of the non-zero left $A$-module $AL$}
If  $\ordn$ is the zero ordinal, then $\RingS$ reduces to the ground ring $\ring$ and $A$ becomes a nicely essential subring in $\ring$,  so that,  $A\cap AL$ is obviously  a nicely essential left $A$-submodule in  $AL$. Suppose that $\ordn$ is non-zero
while the desired result holds for all ordinal $\kappa<\ordn$.  
For a limit ordinal $\ordn$,  (by the definition of  \anesspecial subextension of $\RingS$)  it holds for some subset $\Lambda\subset \ordn$ that $A=\cup\sb{\kappa \in\Lambda} A\sb{\kappa}$ is an increasing union of truncated subrings where each $A\sb{\kappa}$ for $\kappa\in\Lambda$ is   \anesspecial subextension in $\RingS[\kappa]$. The induction  hypothesis applies yielding that each $A\sb{\kappa}\cap A\sb{\kappa} L$ for $\kappa\in\Lambda$ is a nicely  essential left $A\sb{\kappa}$-submodule in $A\sb{\kappa} L$.  Since every finite subset in $AL$ must entirely live  in some $A\sb{\kappa} L$ with $\kappa\in\Lambda $, it follows that
$A\cap AL$ is still a nicely essential left $A$-submodule in $AL$. Next turning to the successor ordinal case,  $\ordn=\kappa+1=\kappa\cup \set{\kappa}$ for some ordinal $\kappa$.  Here $A=A\sb{\kappa+1}$ is by  Definition~\ref{def:special-subext} a nicely essential subring of some $\ring$-subring $\eextA$.  Thus $A$ is a nicely essential left $A$-submodule in $\eextA$, and so is $A\cap AL$ in the non-zero left $A$-submodule  $AL$ of $\eextA$. This completes the proof of our first claim. 

As an immediate consequence,
Proposition~\ref{prop:ext-dsum-left.subextskewpol} yields that every direct sum $\oplus\sb{i\in I}L\sb{i}$ of non-zero left ideals of the ground ring  $\ring$  lifts to  a direct sum
$\oplus\sb{i\in I}(A\cap A L\sb{i})$  of non-zero left ideals in $A$.
Thus in view of applying Proposition~\ref{prop:invariance-udim.ring-ext} (in its full strength),  the last part of the theorem follows by the claim about uniform and essential  left ideals, which we  shall now prove. 

\sectitem{Given a uniform left ideal $L\subset\ring$ and an essential left ideal $K\subset \ring$,  $A\cap AL$ is uniform while $A\cap AK$ is essential as left ideals in $A$. Moreover,  
	$\sing(A)=A\cap \sing(\ring)\T$}

We shall reason by transfinite induction  on the  ordinal $\ordn$.  
If $\ordn$ is the zero ordinal, then $\RingS$  reduces to the ground ring $\ring$ and $A$ is a nicely essential subring in $\ring$, so that, the desired result is granted by Proposition~\ref{prop:udim-essringext}.
So turning to the induction step, we suppose that the ordinal  $\ordn$   is non-zero  while the desired result holds   for every ordinal  smaller than $\ordn$.    

\sectitem{The limit ordinal case}  As  already observed a couple of times before, for the lifting of uniform or essential ideals  we shall see that this limit ordinal case  essentially  reduces to an  application of the induction hypothesis; but for the  characterisation of  the left singular ideals, an indirect trick is necessary.   Definition~\ref{def:special-subext}\pmsf(c) yields some subset $\Lambda\subset\ordn$  with \mthbox{A=\cup\sb{\kappa\in\Lambda} A\sb{\kappa}} arising as an increasing union of subrings such that each $A\sb{\kappa}$ for $\kappa\in\Lambda$ is \anesspecial subextension  of the truncated   skew polynomial ring $\RingS[\kappa]$.   Thus the  induction hypothesis yields for all $\kappa\in\Lambda$ that $A\sb{\kappa}\cap A\sb{\kappa}L$  and $A\sb{\kappa}\cap A\sb{\kappa}K$ are respectively uniform and essential as left  deals in $A\sb{\kappa}$.  
Given  any non-zero polynomials $f,g\in A\cap AL$  and $h\in A$, they depend on  finitely many  variables and it holds  for some $\kappa \in\Lambda $  that  $f,g \in A\sb{\kappa}\cap A\sb{\kappa} L$ while $ h \in A\sb{\kappa} $;  it follows that 
$A f\cap A g \supset A\sb{\kappa}f\cap A\sb{\kappa}g\neq 0$ and   $(A\cap A K)\cap A h \supset (A\sb{\kappa}\cap A\sb{\kappa} K)\cap A\sb{\kappa} h \neq 0$. 
In particular, for  every $a\in \sing(\ring)$, $A\cap A\ann[\ring](a)$ is an essential left ideal of $A$ contained in  $\ann[A](a)$, showing that  $\sing(\ring)\subset \sing(A\ring)$. Since moreover  every non-zero left ideal $I$ of $\ring$ lifts to a non-zero left ideal $A\cap AI$ in $A$,  Proposition~\ref{prop:singual-ideal-subext-skewpolyring}  applies showing that
$\sing(A)=A\cap \sing(\ring)\T$.   
So the desired result is proved in the limit ordinal case.

\sectitem{The  successor ordinal case:    $\ordn=\kappa+1=\kappa\cup\set{\kappa}$}  
As we shall now see,  this  very involved step has actually been handled by Theorem~\ref{thm:invaraince-uniform-univariate-case}.  Throughout,  we set:
\startcent  $x=x\sb{\kappa}$,
$\alpha=\alpha[\kappa]$, $\delta=\delta[\kappa]$ and $\ringS=\RingS[\kappa]$
so that, $\RingS=\ringS[x;\alpha,\delta]$.
\stopcent
In view of applying  main Theorem~\ref{thm:invaraince-uniform-univariate-case}, 
the ring   $\ringS=\RingS[\kappa]$ is a multivariate skew polynomial ring of bijective type over the coefficient ring  $\ring$, 
having the following  properties, where \pmsf(i,ii) are expected by the setting of 
Theorem~\ref{thm:invaraince-uniform-univariate-case}.
\startitem   
\item[\pmsf(i)] $\ringS$ is already free as left $\ring$-module with $\ring$-basis given by the truncated subset of standard terms $\T[\kappa]=\T(x\sb{\eta} \with  \eta \in \kappa)$; and by virtue of Proposition~\ref{prop:ann-skewpol-bijcase}, 
every $\ring$-nice polynomial $p\in \ringS$ enjoys a nice annihilator in $\ringS$: $\ann[\ringS](p)=\ringS\ann[\ring](p)$. 
\item[\pmsf(ii)] Every direct sum of left ideals in $\ring$ lifts to a direct sum of left ideals in $\ringS$, (Proposition~\ref{prop:ext-dsum-left.subextskewpol}).
\item[\pmsf(iii)] By the transfinite induction hypothesis, $\sing(\ringS)=\sing(\ring)\T[\kappa]$. By the main assumption of the theorem, $\delta(\sing(\ringS))=\delta(\sing(\ring)\T[\kappa]) =0$.
\stopitem
Now, by Definition~\ref{def:special-subext}\pmsf(c) we know that $A$ is a nicely essential subring of some prescribed \special subextension $\eextA$ in  $\RingS =\ringS[x;\alpha,\delta]$. Since the induced subring $\isubext(\eextA)[\dd[\eextA]] \subset \ringS$ is at least an $\ring$-subring,  Proposition~\ref{prop:singual-ideal-subext-skewpolyring} grants that 
$\sing(\isubext(\eextA)[\dd[\eextA]]) \subset \sing(\ring) \T[\kappa]$, so that, 
$\ring \cap \sing(\isubext(\eextA)[\dd[\eextA]]) \subset \sing(\ring)$. This fact together with  property~\pmsf(iii) above show that
Lemma~\ref{lem:nicepol}\pmsf(c) applies, yielding the crucial fact that:  
\startalignl   
$\eextA$ already satisfies the main condition of Theorem~\ref{thm:invaraince-uniform-univariate-case} about $\ring$-nice polynomial of prescribed degree: $\eextA a$ contains an $\ring$-nice polynomial of $x$-degree $\dd[\eextA]$ for all  non-zero scalar $a\in \ring$.
\stopalignl 
Moreover and still by Definition~\ref{def:special-subext}\pmsf(c),   for each  attainable $x$-degree  \mthbox{m\in\Deg[x](\eextA)}  the induced subring $\isubext(\eextA)[m]$ is  still \anesspecial  subextension of the truncated  skew polynomial ring  $\ringS=\RingS[\kappa]$.    
Thus the induction hypothesis guaranties that $\isubext(\eextA)[m] L$ is uniform while $\isubext(\eextA)[m] K$ is essential as left ideals in $\isubext(\eextA)[m]$.
We have therefore  shown that, for the   \esspecial  subextension $A$ of the univariate skew polynomial ring $\RingS=\ringS[x;\alpha,\delta]$,   all the hypotheses of Theorem~\ref{thm:invaraince-uniform-univariate-case} are satisfied, granting that the left ideals $A\cap A L$ and $A\cap AK$ are respectively uniform and essential in $A$.   In particular, $\sing(\ring)  \subset \sing(A\ring)$ while it is already granted by  the first claim of the present proof that  $A\cap AI\neq 0$ for all non-zero left ideal $I$ of $\ring$. Hence,  Proposition~\ref{prop:singual-ideal-subext-skewpolyring} applies yielding that  $\sing(A)=A\cap \sing(\ring)\T$.
This completes the proof that every   uniform  or essential  left ideal $L$ in $\ring$ lifts to a uniform (or resp., an essential)  left ideal $A\cap AL$ in $A$ and $\sing(A)=A\cap \sing(\ring)\T$.  Hence, the proof of the theorem is complete.
%
\stoppr

We shall now investigate the (semi)primeness of \esspecial subextensions in $\RingS$;  as usual,
\tqt{(semi)prime} is  a short form for \tqt{prime (or respectively, semiprime)}.  
Allowing an ambient univariate skew polynomial ring
$\extringS=\ringS[x;\alpha,\delta]$ of injective  type,  for an $\ring$-subring $A \subset \extringS$ we say that:
\startalignl[/vskip=0]
$A$ is a \emph{\special subextension of injective type} in   $\extringS$ 
provided \mth{A} has the special shape  of Definition~\ref{def:special-subext}\pmsf(a) in which, for  
\mth{0\neq m \in \Deg[x](A)}, 
the induced endomorphism \mth{\alpha[(m)]} of the associated  $\ring$-subring   \mth{\isubext(A)[m]} is only required to be injective.  
\stopalignl[/vskip=0] 

\startprop \label{prop:semiprime-skewpoly} Let $\extringS=\ringS[x;\sigma,\delta]$ be a univariate skew polynomial of injective type, and $A \subset \extringS$ any   \special subextension of injective type. 
Assume that $A\cap \ringS$ has the \ACC on left annihilators. Then $A$ is (semi)prime  provided so is each induced subring $\isubext(A)[m]$ for every $m\in \Deg[x](A)$. 
\stopprop 

\startpr  It helps recording by \cite[Proposition~3.6]{LeroyMatczuk1995} that:  if  $\ringS$ is (semi)prime, then so is $\extringS$. The arguments proving this result shall be adapted and generalized to the present setting. Let $I$ be an ideal in $A$ and $m\in \Deg[x](A)$;  for  every $n\in \N$ we define:
\begin{equation*}
	I\sb{m,n}=\set{p\in \isubext(A)[m] \with \exists w\in \extringS,\, 
		p\splterm[m,n]+w\in I,\, \deg[x](w)<nm}.
\end{equation*}
\sectitem[\rmNb(i)]{If $I$ contains  a non-zero polynomial of degree $m$, then $I\sb{m,1}\neq 0$ and $I\sb{m+d,1}\neq 0$ for every $d \in \Deg[x](A)$} 
Indeed,  suppose that   $I$   contains  a non-zero polynomial $f$ with $\deg[x](f)=m$. We  have $f=px\sp{m}+f'$ for some $m\in \N$, $0\neq p\in \ringS$ and $f'\in A$  with $\deg[x](f')<m$. The special shape of $A$ grants some $u\in \isubext(A)[m]$  regular as element of $\ringS$ with 
$upx\sp{m}=p'\splterm[m] \neq 0$, which shows that $0\neq p'\in I\sb{m,1}$. Next let $d$ be any  attainable $x$-degree in $A$, thus we are granted with  the special element $\splterm[d]=\upsilonup[d]x\sp{d}\in A\cap \regC(\extringS)$ with
$\upsilonup[d]\in\regC(\ringS) \subset \regC(\extringS)$. 
We get that 
\begin{equation*}
	\splterm[d]uf=\upsilonup[d]x\sp{d}p'\upsilonup[m]x\sp{m} +
	\splterm[d] uf' = \upsilonup[d]\alpha[,d](p'\upsilonup[m])  x\sp{m+d}+ \Delta \in I, 
\end{equation*}  
where  the element $\upsilonup[d]\alpha[,d](p'\upsilonup[m]) \in\ringS$ is non-zero (because so  is $\alpha[,d](p'\upsilonup[m])$ while $\upsilonup[d]$ is a regular element of $\ringS$) and  $\Delta$ is some polynomial belonging to $A$ and satisfying  $\deg[x](\Delta)<m+d$. This shows that $I$ contains a polynomial of $x$-degree $m+d$ and so, as before, $I\sb{m+d,1}\neq 0$.

\sectitem[\rmNb(ii)]{Fix $m\in\N$. Then,   $(I\sb{m,n})\sb{n\in \N}$ is an increasing  sequence of left ideals in $\isubext(A)[m]$;  for some $N\sb{m}$ we have  $\rann[{A\cap \ringS}](I\sb{m,n})=\rann[{A\cap \ringS}](I\sb{m,2n})$  for all $n\geq N\sb{m}$; and for all $k,n\in\N$, $\alpha[(m),k](I\sb{m,n})\subset I\sb{m, n+k}$}  Clearly, each $I\sb{m,n}$ is a left ideal in $\isubext(A)[m]$; and right multiplying elements of $I$  by the special element $\splterm[m]\in A\cap \regC(\extringS)$ shows that $I\sb{m,n}\subset I\sb{m,n+1}$  for all $n\in\N$. 
Thus since the ring $A\cap \ringS$ satisfies the \ACC on left annihilators,
it satisfies the \DCC on right annihilators  
and it holds for some $N\sb{m}\in \N$ that $\rann[{A\cap \ringS}](I\sb{m,n})=\rann[{A\cap \ringS}](I\sb{m,2n})$  for all $n\geq N\sb{m}$.
Next,  the  special shape of $A$  guaranties for all $v\in \isubext(A)[m]$ that: 
\begin{equation*}
	\begin{aligned} & 
		\deg[x](\splterm[m] v-\alpha[(m)](v) \splterm[m])<m , \text{ which implies for all  \mth{ n\in \N } that  }  \\ &
		\deg[x](\splterm[m,n] v-\alpha[(m),n](v) \splterm[m,n])<nm.
	\end{aligned}
\end{equation*}
It follows for all $k,n\in\N$  that $\alpha[(m),k](I\sb{m,n})\subset I\sb{m, n+k}$.

\sectitem[\rmNb(iii)]{Now proving that $A$ is (semi)prime if so is each induced subring $\isubext(A)[m]$ for every attainable $x$-degree $m$} 
With sections~\rmNb(i,ii) of the present proof,  the following remaining  arguments
(essentially) repeat those proving the ordinary
univariate analogue of~\cite[Proposition~3.6]{LeroyMatczuk1995}.
Let $I,J$ be any ideals of  $A$ with $IJ=0$.  

Given  each attainable $x$-degree $m$ in $A$, we first check that  the natural number $N\sb{m}$ (afforded by section~\rmNb(ii)) satisfies: $I\sb{m,n}J\sb{m,n}=0$ for all $n\geq N\sb{m}$. 
Indeed letting  $n\geq N\sb{m}$,   it holds that  
$\rann[{A\cap \ringS}](I\sb{m,n})=\rann[{A\cap \ringS}](I\sb{m,2n})$. The relation $IJ=0$ yields in particular that $I\sb{m,n}\alpha[(m),n](J\sb{m,n})=0$, so that, $\alpha[(m),n](J\sb{m,n}) \subset \rann[A\cap \ringS](I\sb{m,n})=\rann[A\cap \ringS](I\sb{m,2n})$ and $I\sb{m,2n}\alpha[(m),n](J\sb{m,n})=0$ as well. But   $\alpha[(m),n](I\sb{m,n})\subset I\sb{m,2n}$
and  we get that:
$\alpha[(m),n](I\sb{m,n}J\sb{m,n})=\alpha[(m),n](I\sb{m,n})\alpha[(m),n](J\sb{m,n})=0$, so that, 
$I\sb{m,n}J\sb{m,n}=0$ as expected.  

Now given the assumption that each  subring
$\isubext(A)[m]$ is prime for every attainable $x$-degree $m$ in $A$, and further assuming that  $J\neq 0$, we have to check that $I=0$. Let $m$ be any fixed attainable $x$-degree. Since $J\neq 0$, section~\rmNb(i) of the proof  shows that $J\sb{d,1}\neq 0$ for some $d\in \N$, and since $m$ is already an attainable $x$-degree in $A$,  we equally have still by \rmNb(i) that  $J\sb{m+d,1}\neq 0$, so that $J\sb{m+d,n}\neq 0$ for all $n\in \N$ (because the sequence $(J\sb{m+d,n})\sb{n}$ is increasing). But then,   $I\sb{m+d,n}J\sb{m+d,n}=0$ for all $n\geq N\sb{m+d}$ and  the primeness of $\isubext(A)[m+d]$ grants that $I\sb{m+d,n}=0$ for all $n\geq N\sb{m+d}$.
In particular, $I\sb{m+d,1}=0$ and  still by \rmNb(i)  it follows that  $I\sb{m,1}=0$. The latter being true for all attainable $x$-degree $m$, another application of section~\rmNb(i) yields that $I=0$ as expected, so that  the ring $A$ is prime.  Likewise, starting with the assumption that   each    $\isubext(A)[m]$ is semiprime, 
the previous discussion with $J=I$ shows that
$A$ is semiprime as well.   
\stoppr 

\startthm\label{thm:Goldie-subext-skewpolring} Let $\ring$ be a   (semi)prime left Goldie ring and  $A$ any \esspecial subextension of $\RingS$.  
Then  $A$ is left Goldie with $\udim(A)=\udim(\ring)$,   and $A$ is also (semi)prime provided it is a strongly \esspecial subextension in $\RingS$.  
\stopthm

Before writing the proof, one observes that for  every $\ring$-subring $B$ of $\RingS$ generated as left $\ring$-module by a subset of $\T$, if $B$ is left Goldie then so is $\ring$ because $\udim(\ring)\leq \udim(B)$ by virtue of 
Proposition~\ref{prop:ext-dsum-left.subextskewpol}. 

\startpr
The coefficient ring $\ring$  being by assumption (semi)prime left Goldie, Goldie's theorem  yields that $\ring$ is also left non-singular. It follows   by Proposition~\ref{prop:singual-ideal-subext-skewpolyring}  that  both $\RingS$ and $A$  are left non-singular, as well as every subring $B\subset \RingS$ with $B\cap BL\neq 0$ for all non-zero left ideal $L$ in $\ring$.  Since moreover as a ring with finite uniform dimension,  $\ring$ has enough uniform left ideals,   Theorem~\ref{thm:inv-udim-subext-skewpolring} yields that   $\udim(\RingS)=\udim(\ring) =\udim(A)$.
Hence, to conclude in particular that the ambient skew polynomial ring 
$\RingS$ is (semi)prime left Goldie, it only remains by Goldie's theorem to prove that $\RingS$ is (semi)prime. 

With the facts already granted by the  previous paragraph,  
we shall  continue proving  the theorem  using transfinite induction on the indexing ordinal $\ordn$.
When $\ordn$ is the zero ordinal,   $\RingS$ coincides with the (semi)prime left Goldie coefficient ring $\ring$ while $A$  is a nicely essential subring in $\ring$ or a strongly nicely essential subring in $\ring$ in case  the additional assumption on $A$ is taken into account. Hence  thanks to Proposition~\ref{prop:udim-essringext-(semi)prime}, all the claims of the theorem  hold for the case of the zero ordinal.

We then turn  to the induction step with $\ordn$ non-zero and   assume that  the theorem holds for every ordinal less than $\ordn$.  Starting with the limit ordinal case,  
the ambient skew polynomial ring    
$\RingS=\cup\sb{\kappa\in\ordn}\RingS[\kappa]$ is evidently (semi)prime as increasing union of rings each of which is (semi)prime by the induction hypothesis.  Thus together with the introductory paragraph, we get that $\RingS$ is a (semi)prime left Goldie ring and, since moreover the \ACC on left annihilators is transferred to subrings,  $A$ is left Goldie,  as well as each $A\sb{\kappa}$ for $\kappa\in\ordn$ such that $A\sb{\kappa}$ is \anesspecial subextension in $\RingS[\kappa]$. Next additionally assume that $A$ is a strongly \esspecial subextension of $\RingS$. Then Definition~\ref{def:special-subext}\pmsf(c)  yields some subset $\Lambda\subset\ordn$  such that  $A=\cup\sb{\kappa\in\Lambda} A\sb{\kappa}$ and for every  $\kappa\in\Lambda$, each $A\sb{\kappa}$  is still a  strongly \esspecial subextension of the truncated skew polynomial ring $\RingS[\kappa]$. Thus the induction hypothesis grants that  $A\sb{\kappa}$ is (semi)prime for $\kappa\in\Lambda$, and so the ring $A=\cup\sb{\kappa\in\Lambda} A\sb{\kappa}$  is (semi)prime as union of (semi)prime rings.  This completes the proof of the theorem for the limit ordinal case. 

We now continue with the successor ordinal case:  $\ordn=\kappa+1=\kappa\cup\set{\kappa}$ for some ordinal $\kappa$.  As usual, set $x=x\sb{\kappa}, \alpha=\alpha[\kappa]$ and $\delta=\delta[\kappa]$, so that, 
$\RingS=\RingS[\kappa][x;\alpha,\delta]$ is a univariate skew polynomial ring of bijective type
over $\RingS[\kappa]$.   For all
ordinal $\kappa'\leq \kappa$, the induction hypothesis grants that  each $\RingS[\kappa']$ is (semi)prime and hence (semi)prime left Goldie  as a (semi)prime left non-singular ring with finite  uniform dimension. In particular, $\RingS[\kappa]$  already has  the \ACC on left annihilators, and  Proposition~\ref{prop:semiprime-skewpoly} applies showing that $\RingS$ is (semi)prime. Thus,  $\RingS$ is semiprime left Goldie, and its subring $A$ (which already has finite uniform dimension) is left Goldie.  It remains to see that $A$ is (semi)prime under the additional assumption that $A$ is a strongly \esspecial subextension in $\RingS$.  By definition  the ring $A=A\sb{\kappa+1}$ is a strongly nicely essential subring in some \special  subextension $\eext(A)$ of the univariate skew polynomial $\RingS[\kappa][x;\alpha,\delta]$, and for every attainable $x$-degree $m$ in $\eextA$, it is required that the associated subring   $\isubext(\eext(A))[m]$   be also a strongly \esspecial subextension of the truncated skew polynomial ring $\RingS[\kappa]$. So the induction hypothesis grants that each $\isubext(\eext(A))[m]$ is (semi)prime and invoking again   Proposition~\ref{prop:semiprime-skewpoly} yields that $\eext(A)$  is (semi)prime. But  then,  thanks to Proposition~\ref{prop:udim-essringext}\pmsf(a), $\udim(\eext(A))=\udim(A)$ is already finite while $\eext(A)$ is already left non-singular (as granted by the the introductory paragraph of the present proof). 
So we have reached the conclusion that  $\eext(A)$ is a (semi)prime left Goldie ring,   and it follows  by virtue of Proposition~\ref{prop:udim-essringext-(semi)prime} that $A$  is   (semi)prime, completing the proof of the theorem.     
\stoppr



\section{Some  \esspecial subextensions in a two-variate  skew polynomial ring} \label{sect:esssp-subext-2variate-skewpolring}
%
We consider integer parameters  
\begin{equation*}
	p,e,r\geq 2 \text{ and } c,d\geq 1, \text{ with
		\mth{p}  a prime number and  \mth{d} a divisor of \mth{e-1}}.
\end{equation*} 
The  ambient ring is  a two-variate Ore extension \mthbox{\RingS=\ring\tup(x,y: yx=(x+px\sp{e}))}  with natural left $\ring$-basis  given by the set 
$\T=\T(x,y)=\set{x\sp{n}y\sp{m} \with m,n\in\N}$ of standard   terms in the independent variables $x,y$. 
Write  $\term(f)$ for  the set of all terms appearing in each $f \in \RingS$.   
Roughly our setting will enable us  to investigate for each pair $a,b$ of positive integers  
the uniform dimension of the  subextension $A$  of  $\RingS$ generated over $\ring$  by the term $x\sp{a}y\sp{b}$ and satisfying the condition that $\term(f)\subset A$  for every non-zero polynomial $f$ in $A$. 

Recall for all $ m, n\in\N$ that,   the binomial coefficient $\binom{n}{m}=\dfrac{n!}{m!(n-m)!}$ vanishes   when  $0\leq  n<m$. We  need the following computation of the $p$-adic valuation of  $\binom{p^{r}}{n}p^n$  for $1\leq n\leq p^r$, where $\mpzc{v}_p(n)$ is the greatest natural number $k$  with  $p^{k} \mid n$. 
\startlem\label{lem:binom.primepower}
For each $n\in\irg[1,p\sp{r}]$,  write $\mpzc{s}_p(n)$ for the sum of the base-$p$ digits of $n$.  We have $\mpzc{v}_p\bigl(\binom{p^{r}}{n}p^n\bigr)\geq r+n-\mpzc{v}_p(n!) \geq r+\frac{(p-2)n+\mpzc{s}_p(n)}{p-1} \geq  r+1$.
\stoplem 

\startpr  

Since $p$ is a prime number,   \mbox{$\mpzc{v}_p(uv)=\mpzc{v}_p(u)\mpzc{v}_p(v)$  for all $0\neq u,v\in\N$};    and it  is clear that  
$\mpzc{v}_p\bigl(\binom{p^{r}}{n}p^n\bigr)
\geq r+n-\mpzc{v}_p(n!)$.
But then, letting $\len$ be  the number of base-$p$ digits of $n$, 
the Legendre formula  gives   
\mthbox{\mpzc{v}_p(n!) =\sum\sbsp[s=1,\len] \lfloor n/p^s\rfloor =
	\frac{n-\mpzc{s}_p(n)}{p-1} \leq  \frac{n-1}{p-1}}, which yields the desired result.
\stoppr 

\startprop \label{prop:expl2-ncom} Let $\ring$ has characteristic  $p\sp{r+1}$. Then we get an $\ring$-ring automorphism  \Seq{\alpha: \ring[x] \to \ring[x]}  with $\alpha(x)=x+px\sp{e}$,   having finite order $p\sp{r}$.
Thus the terms   $x\sp{p^r}$  and $ y\sp{p^r}$  are central elements 
in  the Ore extension $\RingS=\ring[x][y;\alpha]$. 
\stopprop 

\startpr 
For $0\neq g\in \ring[x]$,  the valuation $\mpzc{v}_x(g)$ denotes the largest $k\in\N$  with $x^k$ dividing $g$; while   $\mpzc{v}_x(0)=\infty$.
The definition of $\alpha$ readily yields for every $f\in\ring[x]$ that $\alpha(f)=f+pg$  for some $g\in\ring[x]$   with  $\mpzc{v}_x(g)>\mpzc{v}_x(f-f(0))$.  
Thus one  can inductively  form 
elements $u\sb1,u\sb2, \ldots$  in $\ring[x]$    with:
\begin{align*} &
	u_0=x, u_1=x^e ,\; \alpha(u_m)=u_m+pu_{m+1},\,
	\mpzc{v}_x(u\sb{m+1})> \mpzc{v}_x(u\sb{m}) \text{ and } \\   &
	\alpha\sp{m}(x) =\hskip-0.3em\tsum\sbsp[k=0,m] \textstyle{\binom{m}{k}} p\sp{k} u\sb{k}  \text{ for all } m\in \N.
\end{align*}
%
%
%
%
Thanks to Lemma~\ref{lem:binom.primepower},  the natural number $\binom{p^{r}}{k}p^k$ is a multiple of the characteristic   $p^{r+1}$ of $\ring$ for every $k\geq 1$. And it follows from the above relation that $p^r$ is the least positive integer $m$ with   $\alpha\sp{m}(x) =x$, proving that $\alpha$ is  an $\ring$-ring automorphism  with   finite order $p\sp{r}$.
As an  immediate consequence,  
$\alpha(x^{p^r})=(x+px^e)^{p^r}= x^{p^r} \sum\sbsp[k=0,p^r] \binom{p^{r}}{k}p^k x^{k(e-1)}=x^{p^r}$.  Thus  $x^{p^r}, y^{p^r}$ live in the center of $\RingS$. 
\stoppr 

Next for every $m\in\N$, we define  the following bounds $\Sn[m]=\Sn[m](e), \mrk=\mrk[r,c](e)\in \N$ and  $ \bSn[m]=\bSn[m,r,c](e)$: 
\begin{equation} \label{eq:Sn}
	\begin{aligned}
		\Sn[m] &=\tsum\sbsp[s=0,m-1]e\sp{s}=
		\frac{e\sp{m}-1}{e-1}; \; 
		\mrk\text{ satisfies the relation: } 
		c\Sn[\mrk-1] \leq r <c\Sn[\mrk], \\ 
		\bSn[m] &=c\Sn[m] \text{ for } 0\leq m\leq \mrk;\; 
		\bSn[\mrk+s]=c(s+\Sn[\mrk]) + sr(e-1) \text{ for } s\in\N. 
	\end{aligned}
\end{equation}
It is clear  that: 
$\Sn[0]=0,\, \Sn[1]=1\leq r,\, \Sn[m+1]=1+\Sn[m]e$ for all $m\in\N$.

\def\eexT{\wt{\sT}}
We then consider the two following parametrized  
sets of terms (depending on  the integer parameters $c,d\geq 1$ set above wit $d\mid e-1$):
\begin{equation} \label{eq:T[d]}
	\begin{aligned}
		& \eexT=\set{x\sp{n}y\sp{m} \with m,n\in\N \text{ and } n- cm   \in  d\Z}; \\
		& \sT=\set{x\sp{n}y\sp{m} \with m,n\in\N,\, n-c m\in d\N,\,
			n\leq \bSn[m]}.
	\end{aligned} 
\end{equation}

\startprop \label{prop:expl2-ncom+} Let $\ring$ has characteristic  $p\sp{r+1}$.   
Then  the left $\ring$-submodules  $A=\ring\sT$ and $\eextA=\ring\eexT$ of the Ore extension $\RingS=\ring\tup(x,y\with yx=(x+px\sp{e})y)$ 
are  subrings  such that   $A$ is a strongly nicely essential subring in $\eextA$ while the latter is is a \special subextension in $\RingS$.
\stopprop

Before we turn to the deeply involved proof of the proposition, the reader is prompted on the following points, some of which follow from the details of the proof below.
\startitem
\item[\pmsf(i)] If $d=e-1$, then $A$  is the subextension  of $\RingS$ characterized by:  $x\sp{c}y\in A$ and  $\term(f)\subset A$ for all $f\in A$.
\item[\pmsf(ii)] The  ring   $A$ is not a skew monoid ring  over $\ring$ and it may not be viewed as  skew polynomial ring  extension of $\ring$  in any useful way. The same is true for $\eextA$, unless $d$ equals $1$ where in the last case $\eextA$ coincides with $\RingS$.  
\item[\pmsf(iii)] Assume the case that  $d\neq 1$. Then  the \esspecial subextension $A$ of $\RingS$  is not a \special subextension and one may not rely on classical techniques to directly investigate the uniform dimension of $A$. The \special subextension $\eextA$ of $\RingS$ is not an essential subring of $\RingS$.  
\item[\pmsf(iv)] Consequently,   
$\udim(A)=\udim(\eextA)$.
If $\ring$ has enough uniform left ideals, then
Theorem~\ref{thm:invaraince-uniform-univariate-case} (or better, Theorem~\ref{thm:inv-udim-subext-skewpolring}) applies yielding that: 
the same holds for $A$ and $\eextA$ and  $\udim(A)=\udim(\eextA)=\udim(\ring)=\udim(\RingS)$.
If $\ring$ is (semi)prime left Goldie then the same holds for both $A$  and $\eextA$ thanks to Theorem~\ref{thm:Goldie-subext-skewpolring}.    
\stopitem

\startpr Given  $m,n\in\N$ we first need to explicitly compute $\alpha^m(x^n)$.  
To this last end,  let us inductively define for all $m,n,k\in\N$  the following natural numbers:  
\begin{align*}
	& \bcoef(n,1,k) =\bcoef(n,k)   = \binom{n}{k}; \;  
	\bcoef(n,m+1,k)  = \hskip-0.4em
	\tsum\sbsp[l=0,{N\sb{n,m}(k)}] \hskip-1em \bcoef(n,m,l) \bcoef(n+l(e-1),k-l)= \hskip-0.4em \tsum\sbsp[l={\mpzc{l}\sb{k,n}} ,{{N\sb{n,m}(k)}}] 
	\hskip-0.4em \bcoef(n,m,l)   \bcoef(n+l(e-1),k-l), \\
	& \text{ where one  puts } N\sb{n,m}(k)=\min(k,n\Sn[m]) \text{ and } 
	{\mpzc{l}\sb{k,n}=\max(0,\lceil(k-n)/e\rceil)}.
\end{align*}
With these settings  
we easily derive   the next  formula:
\begin{align} \label{eq:power-conjmap}
	\alpha(x^n) =& (x+px^e)^n=\hskip-1em\tsum\sbsp[k=0,{\min(r,{n})}] 
	\hskip-0.75em\bcoef(n,k) p\sp{k} x\sp{n+k(e-1)};\, 
	\alpha[,m](x^n) = \hskip-1em\tsum\sbsp[k=0,{N\sb{n,m}(r)}] 
	\hskip-0.8em\bcoef(n,m,k) p\sp{k} x\sp{n+k(e-1)}.
\end{align}
The  formula is clear for $m=0,1$. Assume that it holds
for a fixed $m\in\N$. Having in view that $\ring$ has characteristic $p\sp{r+1}$ and  
setting 
\mthbox{N'(k\sb{1})=\min(r-k\sb1,{n+k\sb1(e-1)})} for all $k\sb{1} \in \irg[N\sb{n,m}(r)]$, 
we have:
\begin{equation*}
	\alpha[,m+1](x^n)  =  
	\tsum\sbsp[k\sb1=0,{N\sb{n,m}(r)}]
	\tsum\sbsp[k\sb2=0,{N'(k\sb{1})}] \bcoef(n,m,k\sb1)\bcoef({n+k\sb1(e-1)},k\sb2) p\sp{k\sb1+k\sb2} x\sp{n+(k\sb1+k\sb2)(e-1)}. \\
\end{equation*}
Letting \mth{k=k\sb{1}+k\sb{2}}, one has 
\mth{k\sb1+k\sb2\leq r} \text{ and } \mth{k\sb1+k\sb2\leq n\Sn[m]+(n+n\Sn[m](e-1))=n\cdot(1+\Sn[m]e)=n\Sn[m+1]}, so we get:
\begin{align*}
	\alpha[,m+1](x^n) &=  \hskip-0.88em \tsum\sbsp[k=0,{N\sb{n,m+1}(r)}] \hskip-0.92em \bcoef(n,m+1,k)  p\sp{k} x\sp{n+k(e-1)},  \text{ with }  \bcoef(n,m+1,k)   =\hskip-1em
	\tsum\sbsp[l=0,{N\sb{n,m}(k)}] \hskip-0.88em \bcoef(n,m,l)\bcoef(n+l(e-1),k-l) \text{ as expected.}
\end{align*}	

We now continue  by splitting the claim of the proposition in three statements. 

\sectitem[\rmNb(i)]{Both $A$ and $\eextA$ are  subrings in $\RingS$}
Given any terms \mthbox{\lambda=x\sp{n}y\sp{m},\,\tau=x\sp{n'}y\sp{m'} \in  \eexT},  we  will show that the product $\lambda\tau$ (in $\RingS$) belongs to   \mthbox{\eextA}.  By the definition of $\eexT$, it holds that  
\mthbox{n-cm, n'-cm'\in d\Z}.
There is noting to show  if $m=0$; so we  assume that $m\geq 1$. 
We have:  
\begin{equation*} 
	\lambda\tau=x\sp{n}\alpha[,m](x\sp{n'}) =
	\tsum\sbsp[k=0,{\min(r,{n'\Sn[m]})}] \bcoef(n',m,k) p\sp{k} x\sp{n+n'+k(e-1)}y\sp{m+m'}. 
\end{equation*}
Consider the exponents \mthbox{n\sb{k}=n+n'+k(e-1)} for $0\leq k\leq \min(r,{n'\Sn[m]})$.  Each number \mthbox{n_k-c(m+m')=(n-cm)+(n'-cm')+k(e-1)}  lives in $d\Z$  since this is the case for both $n-cm$ and $n'-cm'$  while  $d$ is a divisor of $e-1$. Thus at this step, we conclude that $\lambda\tau$ lies in $\eextA$ and $\eextA$  is a subring of $\RingS$.

Next continuing with the assumption that $\lambda, \tau \in \sT$, 
we  want to show that the product $\lambda\tau$ lies in $A$. 
For this, it  only remains to verify that  in the above expression of $\lambda\tau$
the exponents \mthbox{n\sb{k}=n+n'+k(e-1)} for $0\leq k\leq \min(r,{n'\Sn[m]})$ 
are bounded by $\bSn[m+m']$. Here, recall by \eqref{eq:Sn} that  
\begin{align*} 
	& \bSn[m]=c\Sn[m] \text{ for } 0\leq m\leq \mrk, \text{ and }
	\bSn[\mrk+s]=c(s+\Sn[\mrk])+ sr(e-1) \text{ for all } s\in\N; \\
	& \mrk\geq 2  \text{ is the natural number with } 
	c\Sn[\mrk-1] \leq r <c\Sn[\mrk].
\end{align*}  
By the definition of $\sT$,  the   numbers \mthbox{m,n,m',n'\in\N} 
are constrained  as it follows: 
\startcent 
\mthbox{mc\leq n\leq \bSn[m]},
\mthbox{m'c\leq n'\leq \bSn[m']}.
\stopcent
Thus the  largest value for the $n_k$'s is 
\mthbox{N=\bSn[m]+\bSn[m'] +\min(r, \Sn[m]\bSn[m'])(e-1)}
and we wish to check that 
$N\leq \bSn[m+m']$. By symmetry (of $N$ in terms of $m,m'$), one may further assume that  $1\leq m\leq m'$. 
\startitem
\item[\pmsf(i).] {\slshape  Case that $m+m'\leq \mrk$}.   Here, 
$\bSn[m]=c\Sn[m]$, $\bSn[m']=c\Sn[m']$, $\bSn[m+m]=c\Sn[m+m']$, so that, 
\begin{align*}
	N &\leq c[\Sn[m]+\Sn[m']+\Sn[m]\Sn[m'](e-1)] =
	c[\Sn[m]+\Sn[m'](1+\Sn[m]e-\Sn[m])] \\
	&=c[\Sn[m]+\Sn[m'](\Sn[m+1]-\Sn[m])]=c[\Sn[m]+\Sn[m']e^{m}]=c\Sn[m+m']=\bSn[m+m'].
\end{align*} 
\item[\pmsf(ii).] {\slshape The case   $m+m' > \mrk$}.   
Here, there are three subcases to consider. 
\startitem
\item {\slshape The subcase: $m\geq \mrk$}. 
Then, $m'\geq \mrk$ as well and 
$\Sn[m]\bSn[m'] \geq  \bSn[\mrk]>r$. Thus setting
$m=\mrk+s$  and $m'=\mrk+s'$ with $s,s'\in\N$, we compute:
\begin{align*}
	N &= \bSn[m]+\bSn[m'] + r(e-1) \\ 
	&=c(s+\Sn[\mrk]) +sr(e-1) +c(s'+\Sn[\mrk]) +s'r(e-1) + r(e-1) \\
	&=c(s+s'+2\Sn[\mrk]) +(s+s'+1)r(e-1) \\
	\bSn[m+m'] &=c(\mrk+s+s'+\Sn[\mrk])+ (\mrk+s+s')r(e-1); \\
	\bSn[m+m']-N &=c(\mrk-\Sn[\mrk]) +(\mrk-1)r(e-1) \\
	\intertext{Applying the fact that  \mth{\Sn[\mrk]=1+\Sn[\mrk-1]e} while
		\mth{c\Sn[\mrk-1]\leq r<c\Sn[\mrk]}, we  have:}
	\bSn[m+m']-N &=c(\mrk-1) -c\Sn[\mrk-1]e+(\mrk-1)r(e-1) \\ 
	&\geq c(\mrk-1) +[(\mrk-1)(e-1)-e]r. 
\end{align*}
Since $e,\mrk\geq 2$, the last expression above is clearly positive if $\mrk\geq 3$. And if $\mrk=2$, then by definition we must have
$c=c\Sn[1]\leq r<c\Sn[2]=c(1+e)$ and
$\bSn[m+m']-N=c(1-e) +r(e-1)=(r-c)(e-1) \geq 0$.
Whence the conclusion that  $N\leq \bSn[m+m']$ as desired.
\item {\slshape The subcase: $1\leq m<\mrk\leq m'$}. 
Here again, it is already granted that $\Sn[m]\bSn[m'] \geq  \Sn[\mrk]>r$, while this time, $\bSn[m]=c\Sn[m]<c\Sn[\mrk]$ and
\begin{align*}
	N &= \bSn[m]+\bSn[m'] + r(e-1) \\
	&= c\Sn[m]+c(m'-\mrk+ \Sn[\mrk])+(m'-\mrk+1) r(e-1);\\
	\bSn[m+m'] &=c(m+m'-\mrk+ \Sn[\mrk])+(m+m'-\mrk)r(e-1); \\
	\bSn[m+m']-N &=cm-c\Sn[m] +(m-1) r(e-1) \;\; (\text{and using that } c\Sn[m]\leq r), \\
	&\geq  c[m+((m-1)(e-1)-1)\Sn[m]]\geq 0,
\end{align*}
where the last inequality is clear because   for $m=1$ we have $\Sn[1]=1$. 
\item {\slshape The subcase: \mthbox{1\leq m\leq m'<\mrk}}.
Here, 
\begin{equation*}
	N= c\Sn[m]+c\Sn[m'] + \min(r,\Sn[m]\bSn[m'])(e-1)\leq c\Sn[m]+c\Sn[m'] +r(e-1).
\end{equation*}
\mbox{Let \mth{m+m'=\mrk+s}   for some $s\geq 1$;  we have:} 
\mthbox{\bSn[m+m']=c(s+\Sn[\mrk]) +sr(e-1)}  and
\mthbox{\bSn[m+m']-N  \geq c(s+\Sn[\mrk]-\Sn[m]-\Sn[m'])+(s-1)r(e-1)}.
Since  \mth{\Sn[\mrk]=1+\Sn[\mrk-1]e}  while  \mth{c\Sn[m],c\Sn[m']\leq c\Sn[\mrk-1]\leq r}, it follows that:   
\begin{align*}	 
	\bSn[m\hskip-0.01em+\hskip-0.01em m']\hskip-0.2em-\hskip-0.2em N &\geq c[1\hskip-0.2em+\hskip-0.2em s \hskip-0.2em+\hskip-0.2em (\Sn[\mrk\hskip-0.01em-\hskip-0.01em 1]\hskip-0.2em-\hskip-0.2em\Sn[m]) \hskip-0.2em+\hskip-0.2em
	(\Sn[\mrk\hskip-0.01em-\hskip-0.01em 1]\hskip-0.2em-\hskip-0.2em\Sn[m'])] 
	\hskip-0.2em+\hskip-0.2em c\Sn[\mrk\hskip-0.01em-\hskip-0.01em 1](e\hskip-0.2em-\hskip-0.2em 2)\hskip-0.2em+\hskip-0.2em(s\hskip-0.2em-\hskip-0.2em 1)r(e\hskip-0.2em-\hskip-0.2em 1)\\  
	&\geq c(1+s) >0.
\end{align*}
\stopitem 
\stopitem 
This completes the proof that the product (in $\RingS$) of any two terms from $\sT $ lies in $A$, establishing that
$A$ is a subring in $\RingS$.  Notice in particular that, since  the term  $x\sp{c}y$ lies in $\sT$,
it follows that $A$  already contains the subextension $B$ of $\RingS$  characterized by the condition that $x\sp{c}y\in B$ and  
$\term(f)\subset B$  for all $f\in B$. And in the case that  $d$ coincides with $e-1$, 
starting from the fact that $x\sp{c}y\in B$ and  applying the \tqt{closure} condition that $\term(f)\subset B$  for all $f\in B$, one  easily shows (by induction on the $y$-degree) that $\sT\subset B$. 

\sectitem[\rmNb(ii)]{It holds  for every finite subset $E\subset \eexT$  that $\upsilon E\subset \sT$ for  some $\upsilon\in\T(x\sp{p\sp{r}},y\sp{p\sp{r}})$; thus,  $A$ is a strongly nicely essential subring in $\eextA$} 
Proposition~\ref{prop:expl2-ncom} grants that the subset 
$\T(x\sp{p\sp{r}},y\sp{p\sp{r}})=\set{x\sp{ap\sp{r}}y\sp{bp\sp{r}} \with a,b\in\N}$ is    contained in the center  of   $\RingS=\ring[x][y;\alpha]$.
Thus for all $\upsilon\in \T(x\sp{p\sp{r}},y\sp{p\sp{r}})$ and 
$\tau\in \T(x,y)$,  the product $\upsilon\tau$  (in $\RingS$) is till a term living in $\T(x,y)$.    Given every $\tau\in \eexT$, we will show that
$\upsilon[\tau]\tau\in\sT$  for some 
\mthbox{\upsilon[\tau]\in \sT\cap \T(x\sp{p\sp{r}},y\sp{p\sp{r}})}.  
By definition, 
$\tau=x\sp{cm+\varepsilon kd}y\sp{m}=x\sp{\varepsilon kd} \cdot x\sp{cm}y\sp{m}$ with  $m,k\in\N$ and $\varepsilon\in\set{1,-1}$ (such that $cm+\varepsilon kd\geq 0$), while  the term $x\sp{cm}y\sp{m}$ already lies in $\sT$. Thus we only need to find $\upsilon, \lambda\in \sT\cap \T(x\sp{p\sp{r}},y\sp{p\sp{r}})$   with $\upsilon x\sp{kd}, \lambda x\sp{-kd} \in  \sT$.
Choosing a positive natural number 
\mthbox{a\geq \lceil\frac{1}{p\sp{r}}\cdot (\frac{k}{r}+\mrk)\rceil} and setting
$\upsilon=x\sp{cap\sp{r}}y\sp{ap\sp{r}} \in\sT$, we will
check that the product-term $\upsilon x\sp{kd}=x\sp{cap\sp{r}+kd}y\sp{ap\sp{r}}$ still lies in $\sT$. Indeed, 
since $p\geq 2$ and $a\geq 1$, it holds that $ap\sp{r}  > r$, so that the bound $\bSn[ap\sp{r}]$, which  is the highest $x$-degree of terms with $y$-degree $ap\sp{r}$ in $\sT$, is given  according to the defining equation~\eqref{eq:Sn}  by: 
$\bSn[ap\sp{r}]=c\cdot (ap\sp{r}-\mrk+\Sn[\mrk]) +(ap\sp{r}-\mrk)r(e-1)$. So having in view that $\Sn[m]\geq m$ for all $m\in\N$ while $d$ is positive divisor of $e-1$, we get:
\begin{align*}
	\bSn[ap\sp{r}]-(cap\sp{r}+kd) &=
	c\cdot (-\mrk+\Sn[\mrk]) +(ap\sp{r}-\mrk)r(e-1)-kd \\
	&\geq [(ap\sp{r}-\mrk)r-k](e-1)\geq 0.
\end{align*}
Hence, \mthbox{\upsilon x\sp{kd}\in \sT} as expected. 
Next,   let us show that one can find \mthbox{b,s\in\N}  such that
both terms \mthbox{\lambda =x\sp{(cb+sd)p\sp{r}}y\sp{bp\sp{r}}}  and \mthbox{\lambda x\sp{-kd}=x\sp{(cb+sd)p\sp{r}-kd}y\sp{bp\sp{r}}} live in $\sT$.
Since \mthbox{(cb+sd)p\sp{r}-c bp\sp{r} \in d\N} and \mthbox{[(cb+sd)p\sp{r}-kd]-cbp\sp{r} =(sp\sp{r}-k)d}, 
we only need  to ensure that one can choose $b$  and $s$  such that  
\begin{align*}
	k\leq sp\sp{r} \text{ and } (cb+sd)p\sp{r} &\leq \bSn[bp\sp{r}] =
	c\cdot (bp\sp{r}-\mrk+\Sn[\mrk]) +(bp\sp{r}-\mrk)r(e-1). \\
	\intertext{Since \mthbox{\Sn[\mrk]-\mrk\geq 0} and \mthbox{1\leq d\leq e-1}, it suffices to find the natural numbers  \mth{b,s} such that}
	\frac{k}{p\sp{r}} \leq s &\leq (b-\frac{\mrk}{p\sp{r}})r.
\end{align*}
But since  $p\sp{r}\geq 2$, we may take $b=k +\mrk$ and $s=k$.  
This completes  the proof that for all  $\tau\in \eexT$   there exists some  $\upsilon[\tau] \in \sT\cap\T(x\sp{p\sp{r}},y\sp{p\sp{r}})$ with $\upsilon[\tau]\tau\in\sT$. 

Now since  $\T(x\sp{p\sp{r}},y\sp{p\sp{r}})$ is contained in the center of $\RingS$,  for every  finite subset  $E\subset \eexT$ of terms,  the term
$\upsilon=\prod\sb{\tau\in E}\upsilon[\tau]$   lies in
$\sT\cap \T(x\sp{p\sp{r}},y\sp{p\sp{r}})$  and $\upsilon E\subset \sT$. 
Since moreover every member of  $\T(x,y)$ is already a regular element in $\RingS$,  
we conclude that $A$ is a strongly nicely essential subring of  $\eextA$.

\sectitem[\rmNb(iii)]{$\eextA$  is a \special subextension in $\RingS$}
Regard $\RingS=\ring[x][y;\alpha]$ as univariate Ore extension over the polynomial ring $\ring[x]$. Recalling the definitions,  
\begin{align*}
	\eextA &=\ring\eexT \text{ with } \eexT=\set{x\sp{cm+dk}y\sp{m} \with m \in \N,\, -cm/d\leq  k\in\Z}; \text{ thus, } \\ 
	\eextA &=  \sum\sb{m\in\N} (\eextA\cap \ring[x]y\sp{m}) \text{ and } 
	\eextA\cap \ring[x]y\sp{m}=  \ring\cdot \set{x\sp{cm+dk} \with -cm/d\leq  k\in\Z} y\sp{m}.	
\end{align*}
With $m\in \N$  fixed, the terms $x\sp{cm}$ and $\splterm[m]=x\sp{cm}y\sp{m}$ are regular elements of $\RingS$  with
$\splterm[m] \in \eextA\cap \ring[x]y\sp{m}$; and setting $\isubext(\eextA)[m]=\ring[x\sp{d}]$, it is clear that for every $f\in \eextA\cap \ring[x]y\sp{m}$ there is a suitable $k\in\N$  such that $x\sp{dk}f$ lies in $\isubext(\eextA)[m]\splterm[m]=\ring[x\sp{d}]x\sp{cm}y\sp{m}$.  Moreover since $d$ divides $e-1$,  formula~\eqref{eq:power-conjmap} already grants that the conjugation   $\ring$-ring automorphism  \Seq{\alpha: \ring[x] \to \ring[x]} (with $\alpha(x)=x+px\sp{e}$) and its $m$th power $\alpha[,m]$ restrict   to automorphisms of $\ring[x\sp{d}]$, such that:
\begin{align*}
	\text{for all \mth{f\in \ring[x\sp{d}]}},\;\,  \splterm[m]f &=x\sp{cm}y\sp{m} f = 
	x\sp{cm}\alpha[,m](f)y\sp{m}=   \alpha[,m](f)\splterm[m].
\end{align*}
Hence we have proved that $\eextA$ 
is a \special subextension of $\RingS$ viewed as a univariate Ore extension of bijective type (and with vanishing derivation map) over the ordinary polynomial ring $\ring[x]$.  
\stoppr

For  $n\geq 3$ and leaving it to the reader,  one sees how to easily form  $n$-variate  \special subextensions $\eextA$ in an ambient $n$-variate skew polynomial rings $\RingS$ such that $\eextA$ is not an essential subring in $\RingS$. Likewise, one can also form $n$-variate strongly nicely essential subextensions   in  $\eextA$ such that $A$ is not \special in $\RingS$, 
though proving this last fact should  become very tedious (as illustrated by the  proof of Proposition~\ref{prop:expl2-ncom+}).   

\subsection*{Acknowledgements}
The author   would like to thank
the managing editor and the anonymous referee for their kind consideration and   valuable  comments.

%
%


\providecommand{\bysame}{\leavevmode\hbox to3em{\hrulefill}\thinspace}
\providecommand{\MR}{\relax\ifhmode\unskip\space\fi MR }
\providecommand{\MRhref}[2]{%
	\href{http://www.ams.org/mathscinet-getitem?mr=#1}{#2}
}
\providecommand{\href}[2]{#2}

 \end{document}